\documentclass[
	english
]{scrartcl}
\usepackage[T1]{fontenc}
\usepackage[utf8]{inputenc}

\usepackage{amsmath,amsfonts,amsthm,amssymb}
\usepackage{mathrsfs}
\usepackage[colorlinks,citecolor=blue]{hyperref}
\usepackage{cancel}
\usepackage{relsize}
\usepackage{algorithmic}

\usepackage{changes}

\usepackage{enumitem}
\setlist[enumerate]{label=(\alph*)}
\usepackage{pdfrender}

\usepackage[figure]{hypcap}
\usepackage{graphicx}
\graphicspath{{./img/}}
\usepackage{subcaption} 

\usepackage{preamble}

\usepackage{marginnote}

\usepackage[capitalise,nameinlink]{cleveref}
\allowdisplaybreaks

\numberwithin{equation}{section}

\crefname{assumption}{Assumption}{Assumptions}
\Crefname{ALC@unique}{Step}{Steps}

\newcommand\norm[1]{\left\Vert#1\right\Vert}
\newcommand\nnorm[1]{\Vert#1\Vert}

\newcommand\N{\mathbb{N}}
\newcommand\R{\mathbb{R}}

\newcommand{\JJ}{\mathscr{J}}

\newcommand\tto{\rightrightarrows}

\newcommand{\cl}{\operatorname{cl}}
\newcommand{\spa}{\operatorname{span}}

\newcommand{\sgn}{\operatorname{sgn}}
\newcommand{\intr}{\operatorname{int}}

\newcommand{\supp}{\operatorname{supp}}

\renewcommand{\subseteq}{\subset}

\DeclareMathAlphabet{\mathpzc}{OT1}{pzc}{m}{it}
\newcommand\oo{\mathpzc{o}}

\newtheorem{theorem}{Theorem}[section]
\newtheorem{lemma}[theorem]{Lemma}
\newtheorem{proposition}[theorem]{Proposition}

\newtheorem{corollary}[theorem]{Corollary}
\newtheorem{remark}[theorem]{Remark}
\newtheorem{definition}[theorem]{Definition}
\newtheorem{example}[theorem]{Example}

\makeatletter
\long\def\@firstoffiveparen#1#2#3#4#5{\textup{\tagform@{#1}}}
\def\eqref@nolink#1{\textup{\tagform@{\ref*{#1}}}}
\def\eqref@link#1{%
\Hy@safe@activestrue
\expandafter\@setref\csname r@#1\endcsname\@firstoffiveparen{#1}%
\Hy@safe@activesfalse
}
\protected\def\eqref{\@ifstar\eqref@nolink\eqref@link}
\makeatother

\usepackage{mathtools}

\begin{document}

\title{%
	Approximate stationarity in disjunctive optimization:
	concepts, qualification conditions, and application to MPCCs
	}%
\author{%
	Isabella K\"aming%
	\footnote{%
		Technische Universit\"at Dresden,
		Institute of Numerical Mathematics,
		01062 Dresden,
		Germany,
		\email{isabella.kaeming@tu-dresden.de},
		\orcid{0009-0001-8371-6025}
		}
	\and
	Patrick Mehlitz%
	\footnote{%
		Philipps-Universit\"at Marburg,
		Department of Mathematics and Computer Science,
		35032 Marburg,
		Germany,
		\email{mehlitz@uni-marburg.de},
		\orcid{0000-0002-9355-850X}%
		}
	}

% \date{\today}
\publishers{}
\maketitle

\begin{abstract}
	 In this paper, we are concerned with stationarity conditions and qualification conditions
	 for optimization problems with disjunctive constraints.
	 This class covers, among others, optimization problems with complementarity,
	 vanishing, or switching constraints, which are notoriously challenging due
	 to their highly combinatorial structure.
	 The focus of our study is twofold.
	 First, we investigate approximate stationarity conditions and
	 the associated strict constraint qualifications which can be used to infer
	 stationarity of local minimizers.
	 While such concepts are already known in the context of so-called Mordukhovich-stationarity,
	 we introduce suitable extensions associated with strong stationarity.
	 Second, a qualification condition is established which, 
	 based on an approximately Mordukhovich- or strongly stationary point, can be used to infer 
	 its Mordukhovich- or strong stationarity, respectively.
	 In contrast to the aforementioned strict constraint qualifications,
	 this condition depends on the involved sequences justifying approximate stationarity and, thus,
	 is not a constraint qualification in the narrower sense.
	 However, it is much easier to verify as it merely requires to check the (positive)
	 linear independence of a certain family of gradients.
	 In order to illustrate the obtained findings, 
	 they are applied to optimization problems with complementarity constraints,
	 where they can be naturally extended to the well-known concepts of weak and Clarke-stationarity.
\end{abstract}


\begin{msc}	
	\mscLink{49J53}, \mscLink{90C30}, \mscLink{90C33}, \mscLink{90C46}
\end{msc}

\section{Introduction}\label{sec:intro}

In this paper, 
we are concerned with optimality conditions and qualification conditions
for optimization problems with so-called disjunctive constraints of type
\begin{equation}\label{eq:disjunctive_constraints}
	F(x)\in\Gamma\coloneqq\bigcup\limits_{j=1}^t\Gamma_j,
\end{equation}
where $F\colon\R^n\to\R^\ell$ is continuously differentiable
and $\Gamma_1,\ldots,\Gamma_t\subset\R^\ell$ are convex polyhedral sets.
This type of constraints has been introduced in the seminal paper \cite{FlegelKanzowOutrata2007}
and covers, apart from standard inequality and equality constraints, so-called
complementarity constraints, 
	see e.g.\ \cite{ScheelScholtes2000,Ye2005},
vanishing constraints, 
	see e.g.\ \cite{AchtzigerKanzow2008,HoheiselPablosPooladianSchwartzSteverango2020},
switching constraints,
	see e.g.\ \cite{KanzowMehlitzSteck2021,Mehlitz2020},
cardinality constraints,
	see e.g.\ \cite{BeckHallak2016,XiaoYe2024},
reformulated cardinality constraints of complementarity type,
	see e.g.\ \cite{BurdakovKanzowSchwartz2016,CervinkaKanzowSchwartz2016},
or
relaxed probabilistic constraints,
	see e.g.\ \cite{AdamBranda2016}.
Optimization problems possessing these types of constraints 
possess several practically relevant applications,
while being notoriously difficult due to their inherent combinatorial structure.
Problem-tailored optimality conditions and qualification conditions for optimization problems
with constraints of type \eqref{eq:disjunctive_constraints} can be found, for example, in 
\cite{BenkoCervinkaHoheisel2019,BenkoGfrerer2018,Gfrerer2014,Mehlitz2020a}.
A subclass of disjunctive optimization problems which received special interest
in \cite{BenkoCervinkaHoheisel2019} is the one possessing so-called orthodisjunctive
constraints. In this situation, the sets $\Gamma_1,\ldots,\Gamma_t$ used in
\eqref{eq:disjunctive_constraints} are Cartesian products of closed, possibly unbounded intervals.
Let us mention that all example classes which have been mentioned earlier
possess this orthodisjunctive structure,
see \cite[Section~4]{BenkoCervinkaHoheisel2019} and \cite[Section~5]{Mehlitz2020a}
for suitable representations of $\Gamma$ associated with these problem classes.

Dual first-order stationarity concepts associated with optimization problems
possessing constraints of type \eqref{eq:disjunctive_constraints} involve the
first-order derivative of $F$ and a suitable normal cone to $\Gamma$.
Typically, the regular and limiting normal cone are exploited for that purpose,
see \cite{Mordukhovich2018,RockafellarWets1998}, resulting in the systems of
so-called strong stationarity and Mordukhovich-stationarity (S- and M-stationarity for short),
respectively. 
As the regular normal cone is a subset of the limiting normal cone,
S-stationarity implies M-stationarity, and the converse is typically not true
whenever $t\geq 2$ holds.
In order to guarantee that a local minimizer of the underlying problem 
is at least M-stationary, validity of a qualification condition is needed.
In contrast, it has been shown in \cite{Mehlitz2020b} that local minimizers
are always so-called approximately M-stationary (AM-stationary)
even in the absence of a qualification condition,
and that a rather mild constraint qualification, called AM-regularity,
guarantees that an AM-stationary point is already M-stationary.	
Note that AM-stationarity requires that, along a sequence converging to the reference point, 
perturbed M-stationarity conditions hold such that the involved residuals tend to zero.
AM-regularity ensures that, independently of the chosen objective function,
these sequences are somewhat well-behaved.
In fact, AM-regularity is the weakest among all so-called strict constraint qualifications,
i.e., constraint qualifications which ensure that a
given AM-stationary point is M-stationary.
In the context of geometrically-constrained or nonsmooth optimization, 
similar observations have been made e.g.\ in 
\cite{JiaKanzowMehlitzWachsmuth2023,Mehlitz2023}.
Furthermore, coming back to the aforementioned classes of (ortho)disjunctive problems,
related results have been obtained for mathematical problems with
complementarity constraints (MPCCs),
see e.g.\ \cite{AndreaniHaeserSeccinSilva2019,Mehlitz2023,Ramos2021},
switching-constrained problems,
see e.g.\ \cite{LiangYe2021},
and
problems with (reformulated) cardinality constraints (of complementarity type),
see e.g.\ \cite{KanzowRaharjaSchwartz2021}.
In order to tackle S-stationarity of local minimizers, one could be tempted to 
analogously define a concept of approximate S-stationarity (AS-stationarity)
by replacing the appearing 
limiting normal cone with the regular normal cone in the definition of 
AM-stationarity.
However, besides the fact that this approach could only be used to guarantee M-stationarity again
as the closure of the graph of the regular normal cone mapping is the graph of the limiting
normal cone mapping, it is not hard to see 
that the AS-stationarity concept obtained this way is equivalent to AM-stationarity
for problems with constraints of type \eqref{eq:disjunctive_constraints}.
Consequently, we introduce a novel, tightened concept of strict approximate S-stationarity 
(SAS-stationarity)
that is well-suited to obtain, together with the new condition AS-regularity, S-stationarity.
We show that this approach is reasonable as, given a local minimizer,
both SAS-stationarity and AS-regularity together are strictly weaker
than validity of a problem-tailored version 
of the linear independence constraint qualification,
which is the standard constraint qualification 
known from the literature to ensure S-stationarity of local minimizers for problems with 
disjunctive constraints, see e.g.\ \cite{Mehlitz2020a} for a general study. 

Turning to the special case of orthodisjunctive constraints, we exploit the structure of
the associated regular and limiting normal cone to derive a new 
problem-tailored qualification condition which, inspired by the recent study \cite{KaemingFischerZemkoho2024}
that considers Lipschitzian optimization problems with merely inequality constraints, we call
subset Mangasarian--Fromovitz condition (subMFC).
The latter further improves AM- and
AS-regularity from two main perspectives.
First, subMFC
may be fulfilled even if both AM- and AS-regularity are violated,
while still guaranteeing that an 
AM-stationary or SAS-stationary point is M-stationary 
or S-stationary, respectively. Moreover, in contrast to AM- and AS-regularity,
subMFC does not require to check the behavior of infinitely many sequences, which is
especially advantageous if one particular sequence is already provided, e.g., by the
execution of a numerical algorithm. 

To underline the applicability of our results obtained in the abstract setting of 
(ortho)disjunctive constraints, we apply these findings to the classes of purely 
inequality-constrained problems as
well as to MPCCs as 
special cases of problems with orthodisjunctive constraints. Concerning inequality-constrained problems, 
we find that our notions of approximate stationarity recover
the classical approximate Karush--Kuhn--Tucker conditions, see \cite{AndreaniHaeserMartinez2011,MovahedianPourahmad2024}, 
while the problem-tailored
version of subMFC recovers a related qualification condition from \cite{KaemingFischerZemkoho2024}.
For MPCCs, the application of our approximate
stationarity concepts not only recovers different (but equivalent) notions 
of AM-stationarity from the literature, but also allows for a
straightforward extension to restore further standard approximate
stationarity concepts, namely approximate weak stationarity and approximate Clarke-stationarity, 
see \cite{AndreaniHaeserSeccinSilva2019}.
Moreover, a problem-tailored version of subMFC provides a weak qualification 
condition that allows equivalent characterizations of
weak, Clarke-, Mordukhovich-, and strong stationarity via exploiting the corresponding
approximate stationarity concepts.
To the best of our knowledge, our notion of SAS-stationarity is a novel concept even for MPCCs.
As the tailored version of subMFC can be interpreted as a weakened version
of an MPCC-tailored version of the Mangasarian--Fromovitz constraint qualification, 
it is a very tractable condition that, again, can easily be used to verify exact 
stationarity if only the respective approximate stationarity condition holds, for example as a 
result of a numerical algorithm. 

The remainder of this paper is organized as follows. In \cref{sec:notation}, we
summarize the notation used in this paper and recall fundamental concepts from variational
analysis that are essential to our study. 
Thereafter, in \cref{sec:ODP}, we introduce stationarity 
concepts, approximate stationarity concepts, and regularity notions that address problems with 
constraints of disjunctive structure and investigate their relations.
The special case of problems with orthodisjunctive constraints 
is the focus of \cref{sec:ODP2}.
Therein, we first introduce a new qualification condition 
for orthodisjunctive problems and discuss its
relation to the regarded stationarity concepts. Then we apply our findings to purely
inequality-constrained problems and show that our approach recovers results from the 
literature available for this special case. In \cref{sec:MPCC}, 
we study MPCCs in the light of approximate stationarity and the new
qualification condition. 
In particular, we show that the latter serves as a regularity condition to
establish a connection between stationarity and approximate stationarity. 
Final remarks are given in \cref{sec:final_rem}.

\section{Notation and preliminaries}\label{sec:notation}

In this section, we comment on the notation used in this paper,
and some preliminary results are presented as well.

\subsection{Fundamental notation}

We use $\N$, $\R$, $\R_+$, and $\R_-$ 
to denote the positive integers
as well as the real, nonnegative real, and nonpositive real numbers, respectively.
For some $n\in\N$, we make use of $[n]\coloneqq\{1,\ldots,n\}$.
Let us mention that $0$, depending on the context,
is used to represent the scalar zero as well as the all-zero vector.
The standard sign function is denoted by $\sgn\colon\R\to\{-1,0,1\}$.
Throughout the paper, for natural numbers $n,m\in\N$,
it will be convenient to identify the Cartesian product $\R^n\times\R^m$
with $\R^{n+m}$ in a canonical way.
Given $x\in\R^n$, $\supp(x)\coloneqq\{i\in[n]\,|\,x_i\neq 0\}$ is the support of $x$.
Additionally, whenever $I\subset[n]$ is a nonempty index set,
$x_I\in\R^{|I|}$ denotes the vector which results from $x$ by deleting the
components with indices in $[n]\setminus I$,
and for notational convenience, 
the components of $x_I$ may be indexed with the elements of $I$. 
Let $\norm{\cdot}$ denote the Euclidean norm in $\R^n$ and, with a slight abuse of notation,
the associated induced matrix norm as well.
For $\varepsilon>0$, 
$\mathbb B_\varepsilon(x)\coloneqq\{y\in\R^n\,|\,\norm{y-x}\leq\varepsilon\}$
represents the closed ball around $x$ with radius $\varepsilon$.
For a set $\Omega\subset\R^n$, $\cl\Omega$ and $\spa\Omega$ denote the closure of $\Omega$
and the smallest linear subspace of $\R^n$ containing $\Omega$, respectively.
Whenever $\upsilon\colon\R^n\to\R^m$ is a differentiable function,
$\upsilon'\colon\R^n\to\R^{m\times n}$ is the derivative of $\upsilon$,
and $\upsilon'(x)\in\R^{m\times n}$ is the Jacobian of $\upsilon$ at $x$.
Furthermore, whenever $m\coloneqq 1$ holds,
we exploit $\nabla\upsilon(x)\coloneqq\upsilon'(x)^\top$ 
to represent the gradient of $\upsilon$ at $x$.

\subsection{Variational analysis}

In this paper, we exploit some standard concepts from variational analysis
taken, e.g., from \cite{Mordukhovich2018,RockafellarWets1998}.

For a set-valued mapping $\Upsilon\colon\R^n\tto\R^m$
and a point $\bar x\in\R^n$ such that $\Upsilon(\bar x)$ is nonempty, 
\[
	\limsup\limits_{x\to\bar x}\Upsilon(x)
	\coloneqq
	\left\{y\in\R^m\,\middle|\,
		\begin{aligned}
			&\exists\{x^k\}_{k=1}^\infty\subset\R^n,\,\exists\{y^k\}_{k=1}^\infty\subset\R^m\colon
			\\
			&\qquad
			x^k\to\bar x,\,y^k\to y,\,y^k\in\Upsilon(x^k),\,\forall k\in\N
		\end{aligned}
	\right\}
\]
is called the Painlev\'{e}--Kuratowski limit of $\Upsilon$ at $\bar x$.
By definition, $\limsup_{x\to\bar x}\Upsilon(x)$ is a closed set.

For a closed set $\Omega\subset\R^n$ and some point $\bar x\in\Omega$, 
we refer to
\begin{align*}
	\widehat N_\Omega(\bar x)
	&\coloneqq
	\{\xi\in\R^n\,|\,\forall x\in\Omega\colon\,\xi^\top(x-\bar x)\leq\oo(\norm{x-\bar x})\}
\end{align*}
as the regular normal cone to $\Omega$ at $\bar x$.
The latter is a closed, convex cone.
For the purpose of completeness, we set $\widehat N_\Omega(\tilde x)\coloneqq\emptyset$
for each $\tilde x\notin\Omega$.
The closed cone defined by
\[
	N_\Omega(\bar x)
	\coloneqq
	\limsup\limits_{x\to\bar x}\widehat N_\Omega(x)
\]
is called the limiting normal cone to $\Omega$ at $\bar x$.
By a diagonal sequence argument, it is easy to realize the relationship
\[
	N_\Omega(\bar x)
	=
	\limsup\limits_{x\to\bar x} N_\Omega(x),
\]
which is well known as the robustness of the limiting normal cone.
Again, we set $N_\Omega(\tilde x)\coloneqq\emptyset$ 
for each $\tilde x\notin\Omega$.
Let us note that both the regular and the limiting normal cone
equal the normal cone in the sense of convex analysis if $\Omega$ is convex.

Given some $\ell\in\N$ and, for each $i\in[\ell]$,
some $n_i\in\N$, a closed set $\Omega_i\subset\R^{n_i}$, and some point $x^i\in\Omega_i$,
the regular and limiting normal cone enjoy the product rule
\begin{align*}
	\widehat N_{\Omega_1\times\ldots\times\Omega_\ell}((x^1,\ldots,x^\ell))
	&=
	\widehat N_{\Omega_1}(x^1)\times\ldots\times\widehat N_{\Omega_\ell}(x^\ell),
	\\
	N_{\Omega_1\times\ldots\times\Omega_\ell}((x^1,\ldots,x^\ell))
	&=
	N_{\Omega_1}(x^1)\times\ldots\times N_{\Omega_\ell}(x^\ell),
\end{align*}
see e.g.\ \cite[Proposition~6.41]{RockafellarWets1998}.

In the following lemma, 
we comment on the behavior of the regular and limiting normal cone
to sets which possess a certain kind of disjunctive structure.

\begin{lemma}\label{lem:calculus_normals_disjunctive}
	Let $\Gamma_1,\ldots,\Gamma_t\subset\R^\ell$ be convex polyhedral sets,
	and consider $\Gamma\coloneqq\bigcup_{j=1}^t\Gamma_j$ 
	as well as some point $\bar y\in\Gamma$.
	Then the following assertions hold.
	\begin{enumerate}
		\item For $J(\bar{y}) \coloneqq \{ j\in[t] \mid \bar{y}\in \Gamma_j\}$,
			we have
			\[
				\widehat N_\Gamma(\bar y) 
				=
				\bigcap\limits_{j\in J(\bar{y})} \widehat N_{\Gamma_j}(\bar y),
				\qquad
				N_\Gamma(\bar y)
				\subset
				\bigcup\limits_{j\in J(\bar y)} \widehat N_{\Gamma_j}(\bar y).
			\]
		\item There is some $\varepsilon>0$ such that
			\[
				 N_\Gamma(y)\subset N_\Gamma(\bar y),
				 \quad
				 \forall y\in\Gamma\cap\mathbb B_\varepsilon(\bar y).
			\]
	\end{enumerate}
\end{lemma}
\begin{proof}
	The first assertion is immediately clear from \cite[Lemma~2.2]{Mehlitz2020a}.
	To prove the second assertion,
	we decompose $\Gamma$ into a disjoint partition
	$\bar\Gamma_1,\ldots,\bar\Gamma_{\bar t}\subset\R^\ell$
	whose components are nonempty, relatively open, and convex such that $\cl\bar \Gamma_j$ 
	is convex polyhedral for each $j\in[\bar t]$
	and
	\[
		\bar \Gamma_j\cap\cl\bar \Gamma_{j'}\neq\emptyset
		\quad\Longrightarrow\quad
		\bar \Gamma_j\subset\cl\bar \Gamma_{j'},
		\quad
		\forall j,j'\in[\bar t].
	\]
	This is possible due to \cite[Lemma~1]{AdamCervinkaPistek2016}.
	Let $\bar j\in[\bar t]$ 
	be the uniquely defined index such that $\bar y\in\bar \Gamma_{\bar j}$.
	For any $\tilde y\in\Gamma$ and the associated uniquely determined index $\tilde j$ 
	such that $\tilde y\in\bar \Gamma_{\tilde j}$,
	we define
	\[
		\mathcal I(\tilde j)
		\coloneqq
		\{j\in[\bar t]\,|\,\bar \Gamma_{\tilde j}\cap\cl\bar \Gamma_j\neq\emptyset\}.
	\]
	One can easily check that there is some $\varepsilon>0$ such that, for each
	$\tilde y\in\Gamma\cap\mathbb B_\varepsilon(\bar y)$,
	$\mathcal I(\tilde j)\subset\mathcal I(\bar j)$.
	Hence, \cite[Corollary~1]{AdamCervinkaPistek2016} yields
	$N_\Gamma(\tilde y)\subset N_\Gamma(\bar y)$.
\end{proof}

\section{Approximate stationarity conditions for disjunctive optimization problems}\label{sec:ODP}

Consider the disjunctive problem
\begin{equation}\label{eq:disjunctive_problem}\tag{DP}
		\min\limits_x\quad f(x)	\quad\text{s.t.}\quad 
		F(x)\in \Gamma 
		\coloneqq 
		\bigcup\limits_{j=1}^t\Gamma_j,
\end{equation}
where $f\colon\R^n\to\R$ and $F\colon\R^n\to\R^\ell$ are continuously differentiable 
and $\Gamma_1, \dotsc, \Gamma_t \subseteq \mathbb{R}^\ell$ are convex polyhedral sets.
Let us emphasize that some of the results in this section hold for arbitrary closed
sets $\Gamma$, and we will mention it explicitly whenever this is the case.

For $x\in\R^n$ and $\delta\in\R^\ell$, if $F(x)-\delta\in\Gamma$, we define
the index set of the active components as
\begin{align*}
	J(x,\delta) 
	&\coloneqq 
	\{ j\in [t] \mid F(x)-\delta \in \Gamma_j \}.
\end{align*}
Furthermore, whenever $\bar x\in\R^n$ is feasible for \eqref{eq:disjunctive_problem},
we use $J(\bar x)\coloneqq J(\bar x,0)$ for brevity.
For any sequence $\{ (x^k,\delta^k) \}_{k=1}^\infty\subset \R^{n+\ell}$ with
$(x^k,\delta^k) \to (\bar{x},0)$, 
the inclusion $J(x^k,\delta^k)\subset J(\bar{x})$ clearly holds for all large enough $k\in\N$ as 
$F$ is continuous and $\Gamma_1,\ldots,\Gamma_t$ are closed.

To start, let us define some stationarity concepts which are associated with 
\eqref{eq:disjunctive_problem}.
\begin{definition}\label{def:DP_exact_stationary}
	A point $\bar{x}\subseteq\mathbb{R}^n$ feasible for \eqref{eq:disjunctive_problem} 
	is called
	\begin{enumerate}
		\item\label{item:M_stat_DP}
			Mordukhovich-stationary (M-stationary) if there exists a multiplier
			$ \lambda \in \mathbb{R}^\ell$
			such that
			\begin{subequations}\label{eq:M_stat_DP}
				\begin{align}
					\label{eq:M_stat_DP_x}
					&0=\nabla f(\bar x)+F'(\bar x)^\top\lambda,
					\\
					\label{eq:M_stat_DP_lambda}
					&\lambda\in N_{\Gamma}(F(\bar x));
				\end{align}
			\end{subequations}
		\item\label{item:S_stat_DP}
			strongly stationary (S-stationary) if there exists a multiplier
			$ \lambda \in \mathbb{R}^\ell$
			such that \eqref{eq:M_stat_DP_x} and
			\begin{equation}\label{eq:S_stat_DP} 
					\lambda\in \widehat N_{\Gamma}(F(\bar x)).
			\end{equation}
	\end{enumerate}
\end{definition}

Let us note that the concepts from \cref{def:DP_exact_stationary} do not exploit
the particular structure of the set $\Gamma$ in \eqref{eq:disjunctive_problem}.
Using the fact that $\Gamma$ is the union of finitely many convex polyhedral sets,
condition \eqref{eq:S_stat_DP} is equivalent to
\begin{equation*}%\label{eq:S_stat_DP_intersection}
	\lambda\in\bigcap\limits_{j\in J(\bar x)}\widehat N_{\Gamma_j}(F(\bar x))
\end{equation*}
by means of the union rule for the regular normal cone, see \cref{lem:calculus_normals_disjunctive}.

Below, we present some approximate necessary optimality conditions 
for \eqref{eq:disjunctive_problem} which do not require validity of a constraint
qualification. 
The proof of this result is similar to the one of 
\cite[Proposition~2.5]{DeMarchiJiaKanzowMehlitz2023} and only drafted for convenience
of the reader.
\begin{lemma}\label{lem:sequential_necessary_condition}
	If $\bar x\in\R^n$ is a local minimizer of \eqref{eq:disjunctive_problem},
	then there exists a sequence 
	$\{(x^k,\lambda^k,\delta^k,\varepsilon^k)\}_{k=1}^\infty\subset\R^{n+\ell+\ell+n}$
	such that we have
	\begin{subequations}\label{eq:as_dp}
		\begin{align}
			\label{eq:as_dp_x}   
			\varepsilon^k&=\nabla f(x^k)+F'(x^k)^\top\lambda^k,
			\\
			\label{eq:as_dp_lambda}
			\lambda^k&\in 
			\widehat N_{\Gamma}(F(x^k)-\delta^k)	
		\end{align}
	\end{subequations}
	for all $k\in\N$ as well as the convergences
	\begin{subequations}\label{eq:as_dp_conv}
		\begin{align}
			\label{eq:as_dp_conv_standard}
			(x^k,\delta^k,\varepsilon^k)&\to(\bar x,0,0),
			\\
			\label{eq:as_dp_conv_extra}
			(\lambda^k)^\top\delta^k&\to 0.
		\end{align}
	\end{subequations}
\end{lemma}
\begin{proof}
	We consider a sequence of global minimizers 
	$\{(x^k,y^k)\}_{k=1}^\infty\subset\R^{n+\ell}$
	of the surrogate problem
	\begin{equation}\label{eq:DPk}\tag{DP$(k)$}
		\min\limits_{x,y}\{f(x)+\tfrac k2\norm{F(x)-y}^2+\tfrac12\norm{x-\bar x}^2\,|\,
			x\in\mathbb B_\epsilon(\bar x),\,y\in\mathbb B_1(F(\bar x))\cap \Gamma\},
	\end{equation}
	where $\epsilon>0$ is the radius of local optimality in \eqref{eq:disjunctive_problem}
	associated with $\bar x$.
	From the estimate
	\begin{equation}\label{eq:some_crucial_estimate}
		f(x^k)+\frac k2\nnorm{F(x^k)-y^k}^2+\frac12\nnorm{x^k-\bar x}^2\leq f(\bar x),
		\quad
		\forall k\in\N
	\end{equation}
	one can show $(x^k,y^k)\to(\bar x,F(\bar x))$.
	Hence, without loss of generality, we can omit the constraints
	$x\in\mathbb B_\epsilon(\bar x)$ and $y\in\mathbb B_1(F(\bar x))$
	in \eqref{eq:DPk} and still preserve local optimality of $(x^k,y^k)$ for each $k\in\N$.
	Noting that the remaining problem possesses a continuously differentiable objective function,
	we can rely on the sharpest necessary optimality condition from
	\cite[Theorem~6.12]{RockafellarWets1998}, which is stated in terms of the regular normal cone,
	and find
	\[
		\bar x-x^k = \nabla f(x^k)+k F'(x^k)^\top(F(x^k)-y^k),
		\qquad
		k(F(x^k)-y^k)\in\widehat N_\Gamma(y^k)
	\]
	for all $k\in\N$.
	This motivates to define
	$\delta^k\coloneqq F(x^k)-y^k$, $\varepsilon^k\coloneqq\bar x-x^k$, and
	$\lambda^k\coloneqq k\delta^k$ for each $k\in\N$.
	This shows
	\eqref{eq:as_dp} and \eqref{eq:as_dp_conv_standard}.
	Finally, \eqref{eq:some_crucial_estimate} yields
	\[
		(\lambda^k)^\top\delta^k
		=
		k\nnorm{F(x^k)-y^k}^2
		\to 0,
	\]
	i.e., \eqref{eq:as_dp_conv_extra} has been shown,
	and this completes the proof.
\end{proof}

Let us emphasize that \cref{lem:sequential_necessary_condition}
holds true for any closed set $\Gamma\subset\R^\ell$.
We note that, for the price of a slightly more technical proof
which exploits Ekeland's variational principle
and considers this more general setting,
one can even show that the approximate multipliers $\lambda^k$
can be chosen from the so-called proximal normal cone to
$\Gamma$ at $F(x^k)-\delta^k$ for each $k\in\N$,
see \cite[Example~6.16]{RockafellarWets1998} for the definition
of the proximal normal cone and
\cite[Theorem~8]{MovahedianPourahmad2024} for the mentioned result.
One can easily check that the proximal normal cone enjoys the
same union rule as the regular normal cone, 
and, for convex sets, the regular and proximal normal cone 
coincide by means of \cite[Proposition~6.17]{RockafellarWets1998}.
Thus, for the disjunctive setting from \eqref{eq:disjunctive_problem} 
we are interested in,
\cref{lem:sequential_necessary_condition} recovers \cite[Theorem~8]{MovahedianPourahmad2024}.
As outlined in \cite[Section~3.2]{MovahedianPourahmad2024},
in the setting of standard nonlinear or nonlinear conic optimization,
the necessary optimality conditions from \cref{lem:sequential_necessary_condition}
are slightly stronger than the so-called CAKKT conditions,
see \cite{AndreaniGomezHaeserMitoRamos2022,AndreaniMartinezSvaiter2010}
for a precise definition and further discussion.
Respecting the structural properties of $\Gamma$,
we can replace \eqref{eq:as_dp_lambda} equivalently by
\begin{equation*}%\label{eq:lambda_split_for_dp}
	\lambda^k\in \bigcap\limits_{j\in J(x^k,\delta^k)}
				\widehat N_{\Gamma_j}(F(x^k)-\delta^k),
\end{equation*}
see \cref{lem:calculus_normals_disjunctive} again. 

Let us emphasize that the sequence $\{\lambda^k\}_{k=1}^\infty$ of approximate
multipliers in \cref{lem:sequential_necessary_condition} does not need to be
bounded or even convergent in general.
If it is bounded, one could immediately take the limit $k\to\infty$
along a suitable subsequence
in \eqref{eq:as_dp} to deduce M-stationarity of $\bar x$,
and \eqref{eq:as_dp_conv_extra} would be meaningless.
However, as M-stationarity is not a necessary optimality condition for
\eqref{eq:disjunctive_problem} in general, it may happen that
$\{\lambda^k\}_{k=1}^\infty$ is unbounded.
In this situation, \eqref{eq:as_dp_conv_extra} delivers some extra information,
as this unboundedness, to some extent, 
is compensated by the residuals $\{\delta^k\}_{k=1}^\infty$, which tend to zero.
For simplicity of presentation, and to stay as close as possible to the
available results from the literature, we will neglect the information
from \eqref{eq:as_dp_conv_extra} in our upcoming analysis.
Indeed, the remaining conditions \eqref{eq:as_dp} and \eqref{eq:as_dp_conv_standard}
recover the famous AKKT conditions in the setting of standard nonlinear and
nonlinear conic optimization, as proven in 
\cite[Proposition~3]{MovahedianPourahmad2024}.

\cref{lem:sequential_necessary_condition} now motivates the following definition.
\begin{definition}\label{def:DP_special_sequence}
	Let $\bar x\in\R^n$ be feasible for \eqref{eq:disjunctive_problem}.
	A sequence 
	$\{(x^k,\lambda^k,\delta^k,\varepsilon^k)\}\subset\R^{n+\ell+\ell+n}$
	satisfying the convergences \eqref{eq:as_dp_conv_standard} is referred to as
	\begin{enumerate}
		\item an approximately S-stationary sequence w.r.t.\ $\bar x$
			whenever the conditions \eqref{eq:as_dp} hold for all
			$k\in\N$;
		\item an approximately M-stationary sequence w.r.t. $\bar x$
			whenever the conditions \eqref{eq:as_dp_x} and
			\begin{equation}\label{eq:am_dp_lambda}
				\lambda^k\in N_\Gamma(F(x^k)-\delta^k)
			\end{equation}
			hold for each $k\in\N$.
	\end{enumerate}
\end{definition}

Furthermore, we will make use of the following concept.
\begin{definition}\label{def:DP_special_sequence_SAS}
	Let $\bar x\in\R^n$ be feasible for \eqref{eq:disjunctive_problem}.
	A sequence 
	$\{(x^k,\lambda^k,\varepsilon^k)\}\subset\R^{n+\ell+n}$
	satisfying the convergences
	\begin{equation}\label{eq:as_dp_conv_reduced}
		(x^k,\varepsilon^k)\to(\bar x,0)
	\end{equation}
	is referred to as
	a strictly approximately S-stationary sequence w.r.t.\ $\bar x$
	whenever the conditions \eqref{eq:as_dp_x} and
	\begin{equation}\label{eq:as_dp_lambda_reg}
		\lambda^k\in\widehat N_{\Gamma}(F(\bar x))
	\end{equation}
	hold for each $k\in\N$.
\end{definition}

In the next proposition,
we present an alternative characterization of strictly
approximately S-stationary sequences
that underlines some structural similarities to the notions of
approximately M- and S-stationary sequences.

\begin{proposition}\label{prop:equivalence_of_sas_stationarity}
	Let $\bar x\in\R^n$ be a feasible point for \eqref{eq:disjunctive_problem},
	and let $\{(x^k,\lambda^k,\varepsilon^k)\}\subset\R^{n+\ell+n}$ be a sequence.
	Then the following statements are equivalent.
	\begin{enumerate}
		\item\label{item:sas_via_def}
			The sequence
			$\{(x^k,\lambda^k,\varepsilon^k)\}$ 
			is strictly approximately S-stationary w.r.t.\ $\bar x$.
		\item\label{item:sas_via_extra_delta}
			There exists a sequence $\{\delta^k\}_{k=1}^\infty\subset\R^\ell$
			such that 
			$
				\{ (x^k, \lambda^k, \delta^k, \varepsilon^k) \}_{k=1}^\infty
			$
			satisfies the convergences
			\eqref{eq:as_dp_conv_standard}
			while conditions 
			\eqref{eq:as_dp_x}, $F(x^k)-\delta^k\in\Gamma$, and \eqref{eq:as_dp_lambda_reg}
			hold for all $k\in\N$.
	\end{enumerate}
\end{proposition}
\begin{proof}
	$[\ref{item:sas_via_extra_delta} \Longrightarrow \ref{item:sas_via_def}]$:
	This relation is clear as \eqref{eq:as_dp_conv_reduced} is implied by 
	\eqref{eq:as_dp_conv_standard}.
	
	$[\ref{item:sas_via_def} \Longrightarrow \ref{item:sas_via_extra_delta}]$:
	Let 
	$
		\{ (x^k, \lambda^k, \varepsilon^k) \}_{k=1}^\infty
	$
	be a strictly approximately S-stationary sequence w.r.t.\ $\bar x$.
	It suffices to introduce  
	a sequence $\{ \delta^k \}_{k=1}^\infty\subset\R^\ell$ such that $\delta^k\to 0$
	as $k\to\infty$	and $F(x^k)-\delta^k\in\Gamma$ for all $k\in\N$.
	To this end, let us define 
	$ \delta^k \coloneqq F(x^k)-F(\bar{x})$
	for each $k\in\N$.
	By continuity of $F$ and \eqref{eq:as_dp_conv_reduced}, $\delta^k\to 0$
	is clearly satisfied. Moreover, as $\bar{x}$ is feasible for 
	\eqref{eq:disjunctive_problem},	$F(x^k)-\delta^k=F(\bar{x})\in\Gamma$ 
	holds as well.
\end{proof}

\cref{prop:equivalence_of_sas_stationarity} shows that using a residual sequence
$\{\delta^k\}_{k=1}^\infty\subset\R^\ell$ in the definition of strictly approximately
S-stationary sequences is not necessary, which is why it is dropped in 
\cref{def:DP_special_sequence_SAS} for brevity of presentation.

Based on \cref{def:DP_special_sequence,def:DP_special_sequence_SAS},
the following concepts are reasonable.

\begin{definition}\label{def:DP_stationary}
	A point $\bar{x}\subseteq\mathbb{R}^n$ feasible for \eqref{eq:disjunctive_problem} 
	is called
	\begin{enumerate}
		\item\label{item:AMstat_DP}
			approximately M-stationary (AM-stationary) if there exists an
			approximately M-stationary sequence w.r.t.\ $\bar x$;
		\item\label{item:SASstat_DP}
			strictly approximately S-stationary (SAS-stationary) if there exists a
			strictly approximately S-stationary sequence w.r.t.\ $\bar x$.
	\end{enumerate}
\end{definition}

It is obvious by definition that each M-stationary point is 
AM-stationary and that each S-stationary point is SAS-stationary.
The converse relations, in general, require a qualification condition to hold.
Let us mention that related concepts of AM-stationarity can be found, for example, in
\cite{AndreaniHaeserSeccinSilva2019,DeMarchiJiaKanzowMehlitz2023,JiaKanzowMehlitzWachsmuth2023,KanzowRaharjaSchwartz2021,KrugerMehlitz2022,Mehlitz2023,Ramos2021}.
To the best of our knowledge, the concept of SAS-stationarity is new.

In the following result,
we prove that it does not make a difference whether to claim the existence of an
approximately M- or approximately S-stationary sequence w.r.t.\ a feasible 
point of \eqref{eq:disjunctive_problem} in order to define AM-stationarity.
This result merely requires $\Gamma$ to be a closed set
and does not rely on its disjunctive structure.
However, we also show that, equivalently, one can claim
\begin{equation}\label{eq:multipliers_in_normal_cone}
	\lambda^k\in  N_{\Gamma}(F(\bar x))
\end{equation}
instead of \eqref{eq:am_dp_lambda} for all $k\in\N$ in the definition
of an approximately M-stationary sequence w.r.t.\ $\bar x$,
and this equivalence only holds for sets $\Gamma$ possessing the
disjunctive structure which is claimed in \eqref{eq:disjunctive_problem}.
Let us note that both \eqref{eq:am_dp_lambda} and \eqref{eq:multipliers_in_normal_cone}
might be difficult to evaluate as there does not exist
a straightforward union rule for the limiting normal cone.
However, as $\Gamma$ is the union of finitely many convex polyhedral sets, 
one still has a chance to
compute  \eqref{eq:am_dp_lambda} and \eqref{eq:multipliers_in_normal_cone} 
explicitly via a so-called normally admissible stratification of $\Gamma$,
see \cite{AdamCervinkaPistek2016} and the proof of \cref{lem:calculus_normals_disjunctive}.
Other strategies to compute the limiting normal cone to such unions can be found
in \cite{HenrionOutrata2008}.

\begin{proposition}\label{prop:equivalence_of_stationarities}
	Let $\bar x\in\R^n$ be a feasible point for \eqref{eq:disjunctive_problem}.
	Then the following statements are equivalent.
	\begin{enumerate}
		\item\label{item:am_via_reg_normals}
			There exists an approximately S-stationary sequence w.r.t.\ $\bar x$.
		\item\label{item:am_def}
			The point $\bar x$ is AM-stationary.
		\item\label{item:am_via_fixed_normal_cone}
			There exists a sequence
			$
				\{ (x^k, \lambda^k, \varepsilon^k) \}_{k=1}^\infty
				\subseteq \mathbb{R}^{n+\ell+n}
			$
			satisfying the convergences
			\eqref{eq:as_dp_conv_reduced}
			such that conditions 
			\eqref{eq:as_dp_x} and \eqref{eq:multipliers_in_normal_cone}
			hold for all $k\in\N$. 
	\end{enumerate}
\end{proposition}
\begin{proof}
	$[\ref{item:am_via_reg_normals} \Longrightarrow \ref{item:am_def}]$:
	This relation is clear as \eqref{eq:am_dp_lambda} is implied by \eqref{eq:as_dp_lambda}
	for each $k\in\N$.
	
	$[\ref{item:am_def} \Longrightarrow \ref{item:am_via_fixed_normal_cone}]$:
	This relation is clear as \eqref{eq:multipliers_in_normal_cone}
	is implied by \eqref{eq:am_dp_lambda} for large enough $k\in\N$
	by \cref{lem:calculus_normals_disjunctive}.
	
	$[\ref{item:am_via_fixed_normal_cone} \Longrightarrow \ref{item:am_via_reg_normals}]$:
	Let a sequence
		$
			\{ (x^k, \lambda^k, \varepsilon^k) \}_{k=1}^\infty
			\subseteq \mathbb{R}^{n+\ell+n}
		$
	be chosen satisfying the convergences \eqref{eq:as_dp_conv_reduced}
	such that conditions \eqref{eq:as_dp_x} and \eqref{eq:multipliers_in_normal_cone}
	hold for all $k\in\N$.
	By definition of the limiting normal cone,
	for each $k\in\N$, we find sequences 
	$\{y^{k,\iota}\}_{\iota=1}^\infty\subset\Gamma$ and
	$\{\lambda^{k,\iota}\}_{\iota=1}^\infty\subset\R^\ell$
	such that 
	$y^{k,\iota}\to F(\bar x)$ and $\lambda^{k,\iota}\to\lambda^k$
	as $\iota\to\infty$, and $\lambda^{k,\iota}\in\widehat N_\Gamma(y^{k,\iota})$
	for all $\iota\in\N$.
	For each $k\in\N$, pick $\iota(k)\in\N$ such that
	\begin{equation}\label{eq:diag_sequence_quality}
		\nnorm{y^{k,\iota(k)}-F(\bar x)}\leq\frac1k,
		\qquad
		\nnorm{\lambda^{k,\iota(k)}-\lambda^k}\leq\frac1k.
	\end{equation}
	Next, we define a sequence 
	$
		\{ (x^k, \tilde\lambda^k, \tilde\delta^k, \tilde\varepsilon^k) \}_{k=1}^\infty
		\subseteq \mathbb{R}^{n+\ell+\ell+n}
	$
	by means of
	\begin{equation}\label{eq:some_surrogate_sequences}
		\tilde\lambda^k
		\coloneqq
		\lambda^{k,\iota(k)},
		\quad
		\tilde\delta^k
		\coloneqq
		F(x^k)-y^{k,\iota(k)},
		\quad
		\tilde\varepsilon^k
		\coloneqq
		\varepsilon^k+F'(x^k)^\top(\tilde\lambda^k-\lambda^k),
		\quad
		\forall k\in\N.
	\end{equation}
	By construction, we have
	\begin{equation}\label{eq:properties_surrogate_sequences}
		\tilde\varepsilon^k
		=
		\nabla f(x^k)+F'(x^k)^\top\tilde\lambda^k,
		\quad
		\tilde\lambda^k
		\in 
		\widehat N_\Gamma(F(x^k)-\tilde\delta^k),
		\quad
		\forall k\in\N.
	\end{equation}
	Due to \eqref{eq:diag_sequence_quality} and continuity of $F$, we find
	\[
		\nnorm{\tilde\delta^k}
		\leq
		\nnorm{F(x^k)-F(\bar x)}+\nnorm{F(\bar x)-y^{k,\iota(k)}}
		\leq
		\nnorm{F(x^k)-F(\bar x)} + \frac1k
		\to 
		0.
	\]
	Furthermore, boundedness of $\{x^k\}_{k=1}^\infty$
	and continuous differentiability of $F$ yield
	boundedness of $\{F'(x^k)\}_{k=1}^\infty$, 
	and \eqref{eq:diag_sequence_quality} can be used to obtain
	\[
		\nnorm{\tilde\varepsilon^k}
		\leq
		\nnorm{\varepsilon^k}
		+
		\nnorm{F'(x^k)}\nnorm{\tilde\lambda^k-\lambda^k}
		\leq
		\nnorm{\varepsilon^k}
		+
		\frac1k\nnorm{F'(x^k)}
		\to
		0.
	\]
	Hence, $(x^k,\tilde\delta^k,\tilde\varepsilon^k)\to(\bar x,0,0)$,
	showing that $\{(x^k,\tilde\lambda^k,\tilde\delta^k,\tilde\varepsilon^k)\}_{k=1}^\infty$
	is an approximately S-stationary sequence w.r.t. $\bar{x}$.
\end{proof}

As already mentioned before,
one can easily show the equivalence of~\ref{item:am_via_reg_normals}
and~\ref{item:am_def} in \cref{prop:equivalence_of_stationarities} for any
closed set $\Gamma$ (i.e., in the absence of any disjunctive structure)
by a standard diagonal sequence argument.

The next result is now an immediate corollary of \cref{prop:equivalence_of_stationarities}.
\begin{corollary}\label{cor:SAS_stat_vs_AM_stat}
	If $\bar x\in\R^n$ is an SAS-stationary point of \eqref{eq:disjunctive_problem},
	then it is AM-stationary.
	If $t=1$, the converse holds true as well. 
\end{corollary}

As a corollary of \cref{lem:sequential_necessary_condition} and
\cref{prop:equivalence_of_stationarities}, we find the following result.
\begin{corollary}\label{cor:local_minimizers_AM_stat}
	If $\bar x\in\R^n$ is a local minimizer of \eqref{eq:disjunctive_problem},
	then it is AM-stationary.
\end{corollary}

The subsequently stated example shows that a local minimizer of \eqref{eq:disjunctive_problem}
does neither have to be SAS- nor M-stationary.
\begin{example}\label{ex:non_SAS_stat}
	We consider the disjunctive problem
	\[
		\min\limits_x\quad -x_1-x_2\quad\textup{s.t.}\quad
		(-x_1+x_2^3,x_1,x_2, x_1^3)\in\Gamma\coloneqq\Gamma_1\cup\Gamma_2
	\]
	with
	\[
		\Gamma_1\coloneqq\R_-\times\{0\}\times\R_+\times\R,
		\qquad
		\Gamma_2\coloneqq\R_-\times\R_+\times\{0\}\times\R_-.
	\]
	The origin $\bar x\coloneqq(0,0)$ is the only feasible point and,
	thus, the global minimizer of the problem.
	We also observe that $J(\bar x)=\{1,2\}$.
	Let us pick a sequence 
	$\{(x^k,\lambda^k,\varepsilon^k)\}_{k=1}^\infty\subset\R^{2+4+2}$
	satisfying the convergences \eqref{eq:as_dp_conv_reduced}
	as well as condition \eqref{eq:as_dp_x}, amounting to
	\begin{equation}\label{eq:station_ex}
		\varepsilon^k
		=
		\begin{pmatrix}
			-1 \\ -1
		\end{pmatrix}
		+
		\lambda^k_1
		\begin{pmatrix}
			-1 \\ 3(x_2^k)^2
		\end{pmatrix}
		+
		\lambda^k_2
		\begin{pmatrix}
			1\\0
		\end{pmatrix}
		+
		\lambda^k_3
		\begin{pmatrix}
			0\\1
		\end{pmatrix}
		+
		\lambda^k_4
		\begin{pmatrix}
			3(x_1^k)^2\\0
		\end{pmatrix}
		,
	\end{equation}
	for each $k\in\N$.
	To ensure SAS-stationarity of $\bar x$, 
	it would be necessary to guarantee
	\begin{align*}
		(\lambda^k_1,\lambda^k_2,\lambda^k_3,\lambda^k_4)
		&\in 
		\widehat N_\Gamma((0,0,0,0))
		\\
		&=
		\widehat N_{\Gamma_1}((0,0,0,0))\cap\widehat N_{\Gamma_2}((0,0,0,0))
		\\
		&=
		(\R_+\times\R\times\R_- \times\{0\})\cap(\R_+\times\R_-\times\R \times\R_+)
		=
		\R_+\times\R_-\times\R_- \times\{0\}
	\end{align*}
	for each $k\in\N$, additionally.
	Due to $\varepsilon^k\to (0,0)$ and $\lambda_4^k = 0$,
	this requires $-\lambda^k_1+\lambda^k_2\to 1$,
	which is impossible as $\lambda_1^k \geq 0$ and $\lambda_2^k \leq 0$.
	Hence, we conclude that $\bar x$ is not SAS-stationary.
	
	Concerning AM-stationarity, it suffices by \cref{prop:equivalence_of_stationarities} to 
	fulfill condition \eqref{eq:station_ex} for each $k\in\N$ using 
	\begin{align*}
		(\lambda_1^k,\lambda_2^k,\lambda_3^k,\lambda^k_4)
		\in 
		N_\Gamma((0,0,0,0))
		=
		(\R_+\times\R\times\{0\}\times\{0\})
		&\cup
		(\R_+\times\{0\}\times\R\times\R_+)
		\\
		&\cup
		(\R_+\times\R_-\times\R_-\times\{0\}).		
	\end{align*}
	Choosing 
	\[
		(x_1^k, x_2^k) \coloneqq \left(0,\frac{1}{\sqrt{3k}}\right),
		\quad
		(\lambda_1^k,\lambda_2^k,\lambda_3^k,\lambda_4^k) \coloneqq (k,k+1,0,0), 
		\quad
		\varepsilon^k \coloneqq (0,0)
	\]
	for each $k\in\N$ confirms that $\bar{x}$ is AM-stationary, 
	as expected due to \cref{cor:local_minimizers_AM_stat}.

	However, investigating M-stationarity of $\bar{x}$, condition \eqref{eq:M_stat_DP_x} reads as
	\[
		0
		=
		\begin{pmatrix}
			-1\\ -1
		\end{pmatrix}
		+
		\lambda_1
		\begin{pmatrix}
			-1\\ 0
		\end{pmatrix}
		+
		\lambda_2
		\begin{pmatrix}
			1\\0
		\end{pmatrix}
		+
		\lambda_3
		\begin{pmatrix}
			0\\1
		\end{pmatrix}
		+
		\lambda_4
		\begin{pmatrix}
			0\\0
		\end{pmatrix}
		,
	\]
	which needs to be fulfilled for some
	$
		(\lambda_1,\lambda_2,\lambda_3,\lambda_4)
		\in 
		N_\Gamma((0,0,0,0)).
	$
	Clearly, the above system requires $\lambda_2 > 0$ and $\lambda_3 > 0$ at the same time,
	which is not possible. Thus, $\bar{x}$ is not M-stationary.
\end{example}

To proceed, we aim to explore which kind of condition is required to ensure
that an AM- or SAS-stationary point is already M- or S-stationary, respectively.
Here, we follow the approach from \cite{Mehlitz2020b} and introduce 
so-called approximate (or asymptotic) constraint qualifications.

\begin{definition}\label{def:asymptotic_regularity}
	A point $\bar{x}\in\mathbb{R}^n$ feasible for \eqref{eq:disjunctive_problem} 
	is called
	\begin{enumerate}
		\item
			approximately M-regular (AM-regular) if 
			\[
				\limsup\limits_{x\to\bar x,\,\delta\to 0}
					F'(x)^\top N_\Gamma(F(x)-\delta)
				\subset
				F'(\bar x)^\top N_\Gamma(F(\bar x));
			\]
		\item
			approximately S-regular (AS-regular) if
			\[
				\limsup\limits_{x\to\bar x} 
					F'(x)^\top\widehat N_\Gamma(F(\bar x))
				\subset
				F'(\bar x)^\top \widehat N_\Gamma(F(\bar x)).
			\]
	\end{enumerate}
\end{definition}

Let us start to present equivalent alternatives to define AM-regularity.
The following result, which in parts can already be found in 
\cite[Theorem~5.3]{Mehlitz2020b}, 
parallels \cref{prop:equivalence_of_stationarities}.

\begin{proposition}\label{prop:equivalence_of_regularities}
	Let $\bar x\in\R^n$ be a feasible point for \eqref{eq:disjunctive_problem}.
	Then the following statements are equivalent.
	\begin{enumerate}
		\item\label{item:am_reg_via_reg_normals}
			It holds that
			\[
				\limsup\limits_{x\to\bar x,\,\delta\to 0}
					F'(x)^\top\widehat N_\Gamma(F(x)-\delta)
				\subset
				F'(\bar x)^\top N_\Gamma(F(\bar x)).
			\]
		\item\label{item:am_reg_def}
			The point $\bar x$ is AM-regular.
		\item\label{item:am_reg_via_fixed_normal_cone}
			It holds that
			\[
				\limsup\limits_{x\to\bar x}
					F'(x)^\top N_\Gamma(F(\bar x))
				\subset
				F'(\bar x)^\top N_\Gamma(F(\bar x)).
			\]
	\end{enumerate}
\end{proposition}

\begin{proof}
	$[\ref{item:am_reg_def} \Longleftrightarrow \ref{item:am_reg_via_fixed_normal_cone}]$:
	See \cite[Theorem 5.3]{Mehlitz2020b}. \\
	$[\ref{item:am_reg_def} \Longrightarrow \ref{item:am_reg_via_reg_normals}]$:
	This relation is clear as 
	\[
		\limsup\limits_{x\to\bar x,\,\delta\to 0}
			F'(x)^\top\widehat N_\Gamma(F(x)-\delta)
		\subset			
		\limsup\limits_{x\to\bar x,\,\delta\to 0}
			F'(x)^\top N_\Gamma(F(x)-\delta).
	\]
	$[\ref{item:am_reg_via_reg_normals} \Longrightarrow \ref{item:am_reg_via_fixed_normal_cone}]$:
	Consider $\xi \in \limsup_{x\to\bar x} F'(x)^\top N_\Gamma(F(\bar x))$. By definition of the 
	Painlev\'{e}--Kuratowski limit, there exist sequences $\{ x^k \}_{k=1}^\infty$ and 
	$\{ \xi^k \}_{k=1}^\infty$ with $\xi^k \in F'(x^k)^\top N_\Gamma(F(\bar x))$ for all $k\in\N$ and
	$x^k \to \bar{x}$, $\xi^k \to \xi$ as $k\to\infty$. In particular, for each $k\in\N$, 
	there exists 
	$\lambda^k\in N_\Gamma(F(\bar x))$ such that $\xi^k = F'(x^k)^\top \lambda^k$.
	Following the proof of \cref{prop:equivalence_of_stationarities}, we now use the definition
	of the limiting normal cone to find, for each $k\in\N$, sequences 	
	$\{y^{k,\iota}\}_{\iota=1}^\infty\subset\Gamma$ and
	$\{\lambda^{k,\iota}\}_{\iota=1}^\infty\subset\R^\ell$
	such that 
	$y^{k,\iota}\to F(\bar x)$ and $\lambda^{k,\iota}\to\lambda^k$
	as $\iota\to\infty$, and $\lambda^{k,\iota}\in\widehat N_\Gamma(y^{k,\iota})$
	for all $\iota\in\N$.
	Note that the convergence $\lambda^{k,\iota}\to\lambda^k$ implies 
	that $F'(x^k)^\top \lambda^{k,\iota} \to \xi^k$ as $\iota\to\infty$.
	For each $k\in\N$, pick $\iota(k)\in\N$ such that
	\begin{equation}\label{eq:diag_sequence_quality_reg}
		\nnorm{y^{k,\iota(k)}-F(\bar x)}\leq\frac1k,
		\qquad
		\nnorm{F'(x^k)^\top \lambda^{k,\iota(k)}-\xi^k}\leq\frac1k.
	\end{equation}
	Next, we define 
	\begin{equation*}
		\tilde\lambda^k
		\coloneqq
		\lambda^{k,\iota(k)},
		\quad
		\delta^k
		\coloneqq
		F(x^k)-y^{k,\iota(k)},
		\quad
		\tilde \xi^k
		\coloneqq
		F'(x^k)^\top \lambda^{k,\iota(k)},
		\quad
		\forall k\in\N.
	\end{equation*}
	By construction, we have
	\begin{equation*}
		\tilde\lambda^k
		\in 
		\widehat N_\Gamma(F(x^k)-\delta^k),
		\quad
		\tilde \xi^k
		\in 
		F'(x^k)^\top \widehat N_\Gamma(F(x^k)-\delta^k),
		\quad
		\forall k\in\N.
	\end{equation*}
	Due to \eqref{eq:diag_sequence_quality_reg} and continuity of $F$, we find
	\[
		\nnorm{\delta^k}
		\leq
		\nnorm{F(x^k)-F(\bar x)}+\nnorm{F(\bar x)-y^{k,\iota(k)}}
		\leq
		\nnorm{F(x^k)-F(\bar x)} + \frac1k
		\to 
		0
	\]
	and
	\[
		\nnorm{\tilde \xi^k - \xi} 
		\leq
		\nnorm{\tilde \xi^k - \xi^k} + \nnorm{\xi^k-\xi}
		\leq 
		\frac{1}{k} +  \nnorm{\xi^k-\xi}
		\to 
		0
	\]
	as $k\to\infty$.
	Thus, the sequence $\{ \tilde \xi^k \}_{k=1}^\infty$ fulfills $\tilde \xi^k \to \xi$ 
	as $k\to\infty$ and can, together with the sequence 
	$\{ (x^k,\delta^k) \}_{k=1}^\infty$, be used to verify that 
	\[
		\xi \in \limsup\limits_{x\to\bar x,\,\delta\to 0}
					F'(x)^\top\widehat N_\Gamma(F(x)-\delta)
	\]
	holds. The claim now follows from \ref{item:am_reg_via_reg_normals}.
\end{proof}

We emphasize that the equivalence of assertions 
\ref{item:am_reg_via_reg_normals} and \ref{item:am_reg_def}
in \cref{prop:equivalence_of_regularities},
which has been pointed out in \cite[Lemma~3.10]{Ramos2021} 
in the particular setting of complementarity-constrained optimization, 
holds for each closed set $\Gamma$
and does not require any disjunctive structure.

\cref{prop:equivalence_of_regularities} shows that AM- and AS-regularity are the same
whenever $t=1$ is valid,
and this observation, obviously, extends to situations where $|J(\bar x)|=1$ holds.
Particularly, this is the case whenever $\Gamma$ represents the setting of standard
nonlinear programming, i.e., $\Gamma \coloneqq \R^m_-\times\{0\}^p$,
and in this situation, these two conditions recover the
so-called cone-continuity-property from \cite{AndreaniMartinezRamosSilva2016}.
In more general situations, however, AM- and AS-regularity are independent conditions,
as the following example illustrates.

\begin{example}\label{ex:AM_vs_AS_regularity}
	Let us consider the disjunctive constraint systems $F(x)\in\Gamma^{i}$, $i=1,2$,
	where $F\colon\R\to\R^2$ is given by $F(x)\coloneqq(x,-x^2)$, $x\in\R$, and
	\[
		\Gamma^1\coloneqq(\R_+\times\R)\cup(\R\times\R_+),
		\qquad
		\Gamma^2\coloneqq(\R_+\times\{0\})\cup(\{0\}\times\R_+).
	\]
	We note that $\bar x\coloneqq 0$ is feasible for both systems.
	\\
	Noting that $\widehat N_{\Gamma^1}(F(\bar x))=\{0\}\times\{0\}$, 
	$\bar x$ is AS-regular for the system $F(x)\in\Gamma^1$.
	Furthermore, we find
	\[
		N_{\Gamma^1}(F(\bar x))=(\R_-\times\{0\})\cup(\{0\}\times\R_-)
	\]
	and, thus, $F'(\bar x)^\top N_{\Gamma^1}(F(\bar x))=\R_-$.
	Defining $\{(x^k,\delta^k)\}_{k=1}^\infty\subset\R\times\R^2$ by
	\[
		x^k\coloneqq \frac{1}{2k},\qquad\delta^k\coloneqq F(x^k),\qquad\forall k\in\N,
	\]
	we have $1\in F'(x^k)^\top N_{\Gamma^1}(F(x^k)-\delta^k)$
	for each $k\in\N$, showing that $\bar x$ is not AM-regular
	for the system $F(x)\in\Gamma^1$.
	\\
	Observing that $\R\times\{0\}\subset N_{\Gamma^2}(F(\bar x))$ holds,
	we find $F'(\bar x)^\top N_{\Gamma^2}(F(\bar x))=\R$, 
	and this trivially guarantees that $\bar x$ is an AM-regular point
	of the system $F(x)\in\Gamma^2$.
	From $\widehat N_{\Gamma^2}(F(\bar x))=\R_-\times\R_-$
	we obtain $F'(\bar x)^\top \widehat N_{\Gamma^2}(F(\bar x))=\R_-$.
	However, exploiting the sequence $\{x^k\}_{k=1}^\infty$ from above,
	we find $1\in F'(x^k)^\top\widehat N_{\Gamma^2}(F(\bar x))$
	for each $k\in\N$,
	verifying that $\bar x$ is not AS-regular for the system $F(x)\in\Gamma^2$.
\end{example}

Let us fix a feasible point $\bar x\in\R^n$ for \eqref{eq:disjunctive_problem}.
Obviously, since $f$ is continuously differentiable, $\bar x$ is AM-stationary if and only if
\[
	-\nabla f(\bar x)
	\in
	\limsup\limits_{x\to\bar x,\,\delta\to 0} F'(x)^\top N_\Gamma(F(x)-\delta),
\]
see e.g.\ \cite[Lemma~3.6]{Mehlitz2020b},
while $\bar x$ is SAS-stationary if and only if
\[
	-\nabla f(\bar x)
	\in 
	\limsup\limits_{x\to\bar x} F'(x)^\top\widehat N_\Gamma(F(\bar x)).
\]
Our interest in the regularity concepts from \cref{def:asymptotic_regularity} 
is, thus, based on the following obvious observations.

\begin{proposition}\label{prop:AM_and_AS_reg_CQs}
	Let $\bar x\in\R^n$ be a feasible point for \eqref{eq:disjunctive_problem}.
	Then the following assertions hold.
	\begin{enumerate}
		\item\label{item:am_stat_am_reg}
			If $\bar x$ is AM-stationary and AM-regular,
			then it is M-stationary.
		\item\label{item:sas_stat_as_reg}
			If $\bar x$ is SAS-stationary and AS-regular,
			then it is S-stationary.
	\end{enumerate}
\end{proposition}

\cref{cor:local_minimizers_AM_stat} and \cref{prop:AM_and_AS_reg_CQs} together show that
AM-regularity is a constraint qualification implying M-stationarity of local
minimizers. As illustrated in \cite{Mehlitz2020b}, AM-regularity is a comparatively
mild constraint qualification.
Related concepts of regularity have been studied, for example, in 
\cite{AndreaniHaeserSeccinSilva2019,DeMarchiJiaKanzowMehlitz2023,JiaKanzowMehlitzWachsmuth2023,KanzowRaharjaSchwartz2021,KrugerMehlitz2022,Mehlitz2023,Ramos2021},
and several results similar to \cref{prop:AM_and_AS_reg_CQs}\,\ref{item:am_stat_am_reg}
can be found in the literature.
Recently, it has been shown that AM-regularity can be refined with the 
aid of directional data, see \cite{BenkoMehlitz2025}.
The presence of AM-regularity turned out to guarantee M-stationary of feasible
accumulation points associated with sequences generated by augmented Lagrangian methods,
as such accumulation points are typically AM-stationary,
see e.g.\ \cite{JiaKanzowMehlitzWachsmuth2023,KrugerMehlitz2022,Ramos2021}.

To the best of our knowledge, AS-regularity is a novel concept.
Although \cref{prop:AM_and_AS_reg_CQs}\,\ref{item:sas_stat_as_reg}
seems to be appealing at the first glance,
we emphasize that it might be of limited practical use since
SAS-stationarity is, in general, not a necessary optimality condition
for problem \eqref{eq:disjunctive_problem} if $t\geq 2$,
see \cref{ex:non_SAS_stat} again.
However, we know from \cref{cor:local_minimizers_AM_stat} that
AM-stationarity serves as a necessary optimality condition.
Whenever the associated sequence is already a strictly approximately
S-stationary sequence, see \cref{prop:equivalence_of_sas_stationarity} again, 
then validity of AS-regularity is enough to
guarantee S-stationarity of the underlying local minimizer.

Typically, S-stationarity of local minimizers of
\eqref{eq:disjunctive_problem} can be guaranteed in the presence of 
problem-tailored versions of the linear independence constraint
qualification, as visualized in \cite{Mehlitz2020a}.
The following concept originates from \cite[Definition~3.1]{Mehlitz2020a}.
\begin{definition}\label{def:MPDC_LICQ}
	Let $\bar x\in\R^n$ be a feasible point for \eqref{eq:disjunctive_problem}.
	We say that the linear independence constraint qualification
	(LICQ) holds at $\bar x$ whenever 
	\[
		F'(\bar x)^\top\lambda =0,\,
		\lambda\in\sum_{j\in J(\bar x)}\spa\widehat N_{\Gamma_j}(F(\bar x))
		\quad\Longrightarrow\quad
		\lambda=0
	\]
	is valid.
\end{definition}

Below, we briefly outline that inferring S-stationarity of
a given local minimizer of \eqref{eq:disjunctive_problem}
via its SAS-stationarity and AS-regularity is strictly
milder than just postulating validity of LICQ.

\begin{proposition}\label{lem:LICQ_stronger}
	Let $\bar x\in\R^n$ be a local minimizer for \eqref{eq:disjunctive_problem}
	where LICQ holds.
	Then $\bar x$ is SAS-stationary and AS-regular.
	Particularly, $\bar x$ is S-stationary.
\end{proposition}
\begin{proof}
	Due to \cite[Theorem~3.8]{Mehlitz2020a}, $\bar x$ is S- and, thus, SAS-stationary.
	Hence, it remains to prove AS-regularity of $\bar x$.
	
	To this end, take sequences $\{x^k\}_{k=1}^\infty\subset\R^n$ and
	$\{\lambda^k\}_{k=1}^\infty\subset\widehat N_\Gamma(F(\bar x))$
	as well as some point $\xi\in\R^n$ such that $F'(x^k)^\top\lambda^k\to\xi$.
	Let us show that $\{\lambda^k\}_{k=1}^\infty$ is bounded.
	Assume, on the contrary, that $\{\lambda^k\}_{k=1}^\infty$ possesses
	an unbounded subsequence (without relabeling).
	Set $\tilde\lambda^k\coloneqq\lambda^k/\nnorm{\lambda^k}$ for each $k\in\N$,
	and observe that along a subsequence (without relabeling)
	we have $\tilde\lambda^k\to\tilde\lambda$ for some nonvanishing
	$\tilde\lambda\in\R^\ell$.
	As we have $F'(x^k)^\top\tilde\lambda^k\to 0$ by construction,
	$F'(\bar x)^\top\tilde\lambda=0$ follows.
	Furthermore, since $\widehat N_\Gamma(F(\bar x))$ is a closed cone,
	$\tilde\lambda\in\widehat N_\Gamma(F(\bar x))$ is inherent.
	Noting that we have
	\[
		\widehat N_\Gamma(F(\bar x))
		=
		\bigcap\limits_{j\in J(\bar x)}\widehat N_{\Gamma_j}(F(\bar x))
		\subset
		\sum_{j\in J(\bar x)}\widehat N_{\Gamma_j}(F(\bar x))
		\subset
		\sum_{j\in J(\bar x)}\spa\widehat N_{\Gamma_j}(F(\bar x))
	\]
	from \cref{lem:calculus_normals_disjunctive},
	a contradiction to LICQ has been obtained.
	Hence, $\{\lambda^k\}_{k=1}^\infty$ is bounded and, thus,
	converges along a subsequence (without relabeling)
	to some $\lambda\in\widehat N_{\Gamma}(F(\bar x))$
	by closedness of $\widehat N_{\Gamma}(F(\bar x))$.
	This shows $\xi\in F'(\bar x)^\top\widehat N_{\Gamma}(F(\bar x))$
	and, thus, validity of AS-regularity.
\end{proof}

Our next example shows that LICQ can be more restrictive than
both SAS-stationarity and AS-regularity together,
i.e., the converse relation in \cref{lem:LICQ_stronger} does not hold.

\begin{example}\label{ex:LICQ_stronger}
	Let us consider the disjunctive problem
	\[
		\min\limits_{x}\quad x_1\quad\textup{s.t.}\quad
		(x_1+x_2,x_1,x_2)\in\Gamma\coloneqq\Gamma_1\cup\Gamma_2
	\]
	with 
	\[
		\Gamma_1\coloneqq\R_-\times\{0\}\times\R_+,
		\qquad
		\Gamma_2\coloneqq\R_-\times\R_+\times\{0\}.
	\]
	As the feasible set of this problem reduces to the origin $\bar x\coloneqq(0,0)$,
	the latter is the uniquely determined minimizer of the problem.
	For later use, let us note that $J(\bar x)=\{1,2\}$ and
	\[
		\widehat N_{\Gamma_1}((0,0,0))
		=
		\R_+\times\R\times\R_-,
		\qquad
		\widehat N_{\Gamma_2}((0,0,0))
		=
		\R_+\times\R_-\times\R.
	\]
	From \cref{lem:calculus_normals_disjunctive} we find
	\[
		\widehat N_\Gamma((0,0,0))
		=
		\widehat N_{\Gamma_1}((0,0,0))\cap\widehat N_{\Gamma_2}((0,0,0))
		=
		\R_+\times\R_-\times\R_-,
	\]
	so choosing $\lambda\coloneqq(0,-1,0)$ yields
	\[
		\begin{pmatrix}
			0\\0
		\end{pmatrix}
		=
		\begin{pmatrix}
			1\\0
		\end{pmatrix}
		+
		\lambda_1
		\begin{pmatrix}
			1\\1			
		\end{pmatrix}
		+
		\lambda_2
		\begin{pmatrix}
			1\\0
		\end{pmatrix}
		+
		\lambda_3
		\begin{pmatrix}
			0\\1
		\end{pmatrix},
	\]
	and this proves that $\bar x$ is S- and, thus, SAS-stationary.
	Noting that the mapping $F$ is affine in this example,
	AS-regularity is trivially valid.
	
	Observing that
	\[
		\spa\widehat N_{\Gamma_1}((0,0,0))
		+
		\spa\widehat N_{\Gamma_2}((0,0,0))
		=
		\R^3,
	\]
	LICQ at $\bar x$ claims the linear independence of three distinct
	vectors in $\R^2$, which is impossible.
	Hence, LICQ is violated at $\bar x$.
\end{example}

It remains an interesting and important question for future research to
find more practically relevant conditions which ensure
SAS-stationarity of local minimizers of
\eqref{eq:disjunctive_problem}.

Let us mention that, based on
\cref{prop:equivalence_of_stationarities} and
\cref{cor:local_minimizers_AM_stat},
\[
	\limsup\limits_{x\to\bar x} F'(x)^\top N_\Gamma(F(\bar x))
	\subset
	F'(\bar x)^\top \widehat N_\Gamma(F(\bar x))
\]
provides a constraint qualification ensuring 
S-stationarity of a local minimizer $\bar x\in\R^n$ for \eqref{eq:disjunctive_problem}.
This condition is, in general, too strong to be applicable.
Indeed, it particularly requires 
$F'(\bar x)^\top N_\Gamma(F(\bar x))\subset F'(\bar x)^\top\widehat N_\Gamma(F(\bar x))$,
which is unlikely to hold whenever $|J(\bar x)|\geq~2$.

Let us close the section with a last remark.

\begin{remark}\label{rem:strict_complementarity}
	Let $\bar x\in\R^n$ be a feasible point of \eqref{eq:disjunctive_problem}
	such that $J(\bar x)=\{j_0\}$ for some $j_0\in[t]$.
	Particularly, for all $x\in\R^n$ sufficiently close to $\bar x$ 
	and all $\delta\in\R^\ell$ sufficiently close to $0$,
	we have $J(x,\delta)\subset\{j_0\}$.
	Hence, locally around $\bar x$,
	the disjunctive structure is redundant,
	and the constraints of \eqref{eq:disjunctive_problem}
	can be replaced by $F(x)\in\Gamma_{j_0}$ for any kind of local analysis.
	Noting that $\Gamma_{j_0}$ is convex polyhedral,
	for the localized model, all concepts related to (approximate)
	M- and S-stationarity or the respective regularities
	coincide. 
	In the particular setting of complementarity-constrained
	optimization, this situation is referred to as \emph{strict complementarity},
	and we may adapt this terminology to \eqref{eq:disjunctive_problem} as well.
	Clearly, whenever strict complementarity holds, 
	the main difficulties encapsulated in \eqref{eq:disjunctive_problem}
	are not present anymore.
	In the analysis and examples in this paper,
	we, thus, mainly consider reference points $\bar x$ such that $|J(\bar x)|\geq 2$
	(or, at least, have such points in mind).
\end{remark}

\section{A qualification condition for orthodisjunctive optimization problems}\label{sec:ODP2}

Throughout this section,
we consider the case where the sets $\Gamma_1,\ldots,\Gamma_t$ in \eqref{eq:disjunctive_problem} 
are products of closed intervals, i.e., 
the so-called orthodisjunctive problem
\begin{equation}\label{eq:orthodisjunctive_problem}\tag{ODP}
	\min\limits_x \quad f(x)	
		\quad\text{s.t.}\quad F(x)\in \Gamma 
		\coloneqq \bigcup\limits_{j=1}^t\Gamma_j	\,\,\text{with}\,\,
		\Gamma_j \coloneqq  \prod\limits_{i=1}^\ell [a_i^j, b_i^j], 
		\quad 
		\forall j\in[t]
\end{equation}
with $-\infty\leq a_i^j\leq b_i^j\leq\infty$ 
for all $i\in[\ell]$ and $j\in[t]$.
Here, we used
\[
	[-\infty,r]\coloneqq(-\infty,r],
	\qquad
	[r,\infty]\coloneqq [r,\infty),
	\qquad
	[-\infty,\infty]\coloneqq \R
\]
for $r\in \R$.
To the best of our knowledge, 
the notion of orthodisjunctive optimization problems dates back to 
\cite{BenkoCervinkaHoheisel2019}, where the authors exploited the particular
structure of $\Gamma$ in \eqref{eq:orthodisjunctive_problem} to derive
tractable sufficient conditions for the so-called metric subregularity 
constraint qualification.
Let us recall that most of the popular problem classes in disjunctive optimization, such as
mathematical problems with complementarity, vanishing, cardinality, or switching
constraints, see e.g.\ \cite[Section~5]{Mehlitz2020a} for an overview,
possess the orthodisjunctive structure.

\subsection{A qualification condition for M- and S-stationarity in orthodisjunctive optimization}\label{sec:ODP_stuff}

For $x\in\R^n$ and $\delta\in\R^\ell$ satisfying $F(x)-\delta\in\Gamma$, 
let us define the index sets
\begin{align*}
	I^\forall (x,\delta) 
	&\coloneqq 
	\{ i\in [\ell] \mid \forall j\in J(x,\delta)\colon\,F_i(x)-\delta_i \in \{a_i^j,b_i^j \}\},
	\\
	I^\exists (x,\delta)
	&\coloneqq
	\{ i\in[\ell] \mid \exists j\in J(x,\delta)\colon\,
		F_i(x)-\delta_i \in \{a^j_i,b^j_i\}\}.
\end{align*}
Furthermore, to simplify notation,
for each point $\bar x\in\R^n$ feasible for \eqref{eq:orthodisjunctive_problem},
we will make use of	$I^\exists(\bar x) \coloneqq I^\exists(\bar x,0)$.
Let us note that, for any $x\in\R^n$ sufficiently close to $\bar x$
and any $\delta\in\R^\ell$ sufficiently small in norm satisfying $F(x)-\delta\in\Gamma$,
we find $I^\exists(x,\delta)\subset I^\exists(\bar x)$
from the continuity of $F$.

To proceed, let us fix $x\in\R^n$ and $\delta\in\R^\ell$ such that $F(x)-\delta\in\Gamma$.
From \cref{lem:calculus_normals_disjunctive} we find
\begin{equation}\label{eq:regular_normals_ODP}
	\begin{aligned}
	\widehat N_\Gamma(F(x)-\delta)
	&=
	\bigcap\limits_{j\in J(x,\delta)}\widehat N_{\Gamma_j}(F(x)-\delta)
	\\
	&=
	\bigcap\limits_{j\in J(x,\delta)}\prod\limits_{i=1}^\ell
		\widehat N_{[a^j_i,b^j_i]}(F_i(x)-\delta_i)
	=
	\prod\limits_{i=1}^\ell\bigcap\limits_{j\in J(x,\delta)}
		\widehat N_{[a^j_i,b^j_i]}(F_i(x)-\delta_i)
	\end{aligned}
\end{equation}
and 
\begin{equation}\label{eq:limiting_normals_ODP}
	\begin{aligned}
	N_\Gamma(F(x)-\delta)
	&\subset
	\bigcup\limits_{j\in J(x,\delta)}\widehat N_{\Gamma_j}(F(x)-\delta)
	\\
	&=
	\bigcup\limits_{j\in J(x,\delta)}\prod\limits_{i=1}^\ell
		\widehat N_{[a^j_i,b^j_i]}(F_i(x)-\delta_i)
	\subset
	\prod\limits_{i=1}^\ell\bigcup\limits_{j\in J(x,\delta)}
		\widehat N_{[a^j_i,b^j_i]}(F_i(x)-\delta_i).
	\end{aligned}
\end{equation}

Some consequences of these formulas are summarized in the following remark.

\begin{remark}\label{rem:lambda_AS_ODP}
	Pick $x\in\R^n$ and $\delta\in\R^\ell$ such that $F(x)-\delta\in\Gamma$.
	\begin{enumerate}
	\item\label{item:reg_normals_ODP}
		Using the index sets
		\begin{align*}
			I_=^\forall(x,\delta)
			&\coloneqq
			\{i\in I^\forall(x,\delta) \mid\forall j\in J(x,\delta)\colon\, a_i^j=b_i^j\},
			\\
			I_{\neq}^\forall(x,\delta)
			&\coloneqq
			I^\forall(x,\delta)\setminus I^\forall_=(x,\delta),
			\\
			I_\downarrow^\forall (x,\delta)
			&\coloneqq
			\{i\in I_{\neq}^\forall (x,\delta) \mid 
				\forall j\in J(x,\delta)\colon\,F_i(x)-\delta_i=a^j_i\},
			\\
			I_{\uparrow}^\forall(x,\delta)
			&\coloneqq
			\{i\in I_{\neq}^\forall(x,\delta) \mid 
				\forall j\in J(x,\delta)\colon\,F_i(x)-\delta_i=b^j_i\},
			\\
			I_{\updownarrow}^\forall(x,\delta)
			&\coloneqq
			I_{\neq}^\forall(x,\delta)\setminus
				(I_\downarrow^\forall(x,\delta)\cup I_\uparrow^\forall(x,\delta)),
		\end{align*}
		it can be derived from \eqref{eq:regular_normals_ODP} that
		\[
			\widehat N_\Gamma(F(x)-\delta)
			=
			\left\{
				\lambda\in\R^\ell\,\middle|\,
				\begin{aligned}	
					&\forall i\in ([\ell]\setminus I^\forall(x,\delta))
						\cup I_\updownarrow^\forall(x,\delta)\colon\,\lambda_i=0,
						\\
					&\forall i\in I_\downarrow^\forall(x,\delta)\colon\,\lambda_i\leq 0,
						\\
					&\forall i\in I_\uparrow^\forall(x,\delta)\colon\,\lambda_i\geq 0
				\end{aligned}
			\right\}.
		\]
		Indeed, if, for example, we have $i\in[\ell]\setminus I^\forall(x,\delta)$, then
		$F_i(x)-\delta_i\in(a^{j_0}_i,b^{j_0}_i)$ holds for some $j_0\in J(x,\delta)$,
		and $\widehat N_{[a^{j_0}_i,b^{j_0}_i]}(F_i(x)-\delta_i)=\{0\}$ follows,
		which yields $\lambda_i=0$. The restrictions concerning the other index sets can be
		obtained likewise.
	\item\label{item:lim_normals_trivial_ODP}
		From \eqref{eq:limiting_normals_ODP} we obtain
		\begin{equation}\label{eq:trivial_components_limiting_normal_cone}
			\supp(\lambda)\subset I^\exists(x,\delta),
			\quad
			\forall\lambda\in N_\Gamma(F(x)-\delta).
		\end{equation}
		Indeed, whenever $i\in[\ell]\setminus I^\exists(x,\delta)$, then
		$F_i(x)-\delta_i\in(a^{j}_i,b^{j}_i)$, and, thus,
		$\widehat N_{[a^{j}_i,b^{j}_i]}(F_i(x)-\delta_i)=\{0\}$ holds for all 
		$j\in J(x,\delta)$, 
		which yields $\lambda_i=0$ via \eqref{eq:limiting_normals_ODP}.
	\end{enumerate}
\end{remark}

We note that the formula from \cref{rem:lambda_AS_ODP}\,\ref{item:reg_normals_ODP}
can, particularly, be useful to specify the precise meaning of SAS- and, 
taking \cref{prop:equivalence_of_stationarities} into account, AM-stationarity 
for \eqref{eq:orthodisjunctive_problem}.
\Cref{rem:lambda_AS_ODP}\,\ref{item:lim_normals_trivial_ODP} will become handy later on.

\begin{remark}\label{rem:dp_seq_const}
	Whenever we consider a sequence 
	$\{ (x^k,\lambda^k,\delta^k, \varepsilon^k) \}_{k=1}^\infty\subset\R^{n+\ell+\ell+n}$ 
	as in \cref{def:DP_special_sequence} for \eqref{eq:orthodisjunctive_problem} 
	in the following, 
	we may assume without loss of generality that the sets $J(x^k,\delta^k)$ 
	are the same for all $k\in\N$, and likewise, the sets
	$I^\exists(x^k,\delta^k)$ are the same 
	for all $k\in\N$, which is reasonable as it is always possible 
	to consider a suitable subsequence by pigeonhole principle.
	Similarly, we may also assume that 
	each component of $\lambda^k$ possesses
	a constant sign for all $k\in\N$.
\end{remark}

Now, we are well-prepared to tailor the so-called 
subset Mangasarian--Fromovitz condition from \cite{KaemingFischerZemkoho2024},
where it has been introduced for nonsmooth Lipschitzian optimization problems with
inequality constraints, to \eqref{eq:orthodisjunctive_problem}.
\begin{definition}\label{def:subMFC_ODP}
	Let $\bar{x}\in\R^n$ be an AM-stationary point of \eqref{eq:orthodisjunctive_problem}.
	We say that the \textit{ODP Subset Mangasarian--Fromovitz Condition (ODP-subMFC)}
	holds at $\bar{x}$ if there 
	exist $I \subseteq I^\exists(\bar{x})$ and a sequence 
	$
		\{(x^k,\lambda^k, \delta^k, \varepsilon^k) \}_{k=1}^\infty
		\subset 
		\mathbb{R}^{n+\ell+\ell+n}
	$ 
	such that the following conditions are satisfied.
	\begin{enumerate}[label=(\roman*)]
		\item \label{subMFC_ODP_I} Either $I=\emptyset$, 
			or
			it holds for all 
			$u\in\R^{\ell}\setminus\{0\}$ with $u\geq 0$ and $u_{[\ell]\setminus I} = 0$ that
			\begin{equation}\label{eq:ODP-subMFC_ODP_lin_indep}
				0  
				\neq 
				\sum_{i\in I} \sgn(\lambda_i^k) u_i \nabla F_i(\bar x).
			\end{equation}
		\item\label{subMFC_ODP_II} The sequence 
			$\{ (x^k,\lambda^k,\delta^k,\varepsilon^k) \}_{k=1}^\infty$ 
			is approximately M-stationary w.r.t.\ $\bar x$,
			and $I=I^\exists(x^k,\delta^k)$ is valid for all $k\in\N$.
	\end{enumerate}
\end{definition}

As we will see in \cref{sec:NLPs}, \cref{def:subMFC_ODP} recovers the original
definition of the subset Mangasarian--Fromovitz condition from 
\cite[Definition~3.7]{KaemingFischerZemkoho2024}
in the setting of smooth inequality-constrained problems.

From \cref{sec:ODP}, 
see \cref{prop:equivalence_of_stationarities,prop:equivalence_of_regularities} in particular,
one may anticipate that it does not play a role if the 
sequence $\{(x^k,\lambda^k,\delta^k,\varepsilon^k)\}_{k=1}^\infty$ in \cref{def:subMFC_ODP}
is required to be an approximately M- or approximately S-stationary sequence.
In the subsequently stated result, we formally prove this conjecture.

\begin{proposition}\label{prop:ODPsubMFC_via_AS_stat}
	Let $\bar x\in\R^n$ be an AM-stationary point of \eqref{eq:orthodisjunctive_problem}.
	Then ODP-subMFC holds at $\bar x$ if and only if there exist
	$I\subset I^\exists(\bar x)$ and a sequence 
	$\{(x^k,\lambda^k,\delta^k,\varepsilon^k)\}_{k=1}^\infty\subset\R^{n+\ell+\ell+n}$
	being approximately S-stationary w.r.t.\ $\bar x$
	such that item \ref{subMFC_ODP_I} of \cref{def:subMFC_ODP}
	and $I=I^\exists(x^k,\delta^k)$ for all $k\in\N$ are satisfied.
\end{proposition}
\begin{proof}
	$[\Longleftarrow]$:
		This implication is trivial
		as, for each $k\in\N$,
		we have $\widehat N_\Gamma(F(x^k)-\delta^k) \subset N_\Gamma(F(x^k)-\delta^k)$,
		such that the index set $I$ and the sequence 
		$\{ (x^k,\lambda^k,\delta^k,\varepsilon^k) \}_{k=1}^\infty$,
		which is approximately S-stationary w.r.t.\ $\bar x$,
		satisfy all requirements of \cref{def:subMFC_ODP}.
	
	$[\Longrightarrow]$:
		Pick an index set $I\subset I^\exists(\bar{x})$ and a sequence 
		$\{ (x^k,\lambda^k,\delta^k,\varepsilon^k) \}_{k=1}^\infty$ 
		which satisfy all requirements in \cref{def:subMFC_ODP}, 
		implying in particular $\sgn(\lambda_i^k) \neq 0$ for all $i\in I$
		and $k\in\N$, where we recall \cref{rem:dp_seq_const}.
		By definition of the limiting normal cone,
		for each $k\in\N$, we find sequences 
		$\{y^{k,\iota}\}_{\iota=1}^\infty\subset\Gamma$ and
		$\{\lambda^{k,\iota}\}_{\iota=1}^\infty\subset\R^\ell$
		such that 
		$y^{k,\iota}\to F(x^k)-\delta^k$ as well as $\lambda^{k,\iota}\to\lambda^k$
		as $\iota\to\infty$, 
		and $\lambda^{k,\iota}\in\widehat N_\Gamma(y^{k,\iota})$
		for all $\iota\in\N$.
	
		For each $k\in\N$, pick $\iota(k)\in\N$ large enough such that
		\begin{subequations}\label{eq:diag_sequence_quality_new}
			\begin{align}
			\label{eq:diag_sequence_quality1}
			&\nnorm{y^{k,\iota(k)}-(F(x^k)-\delta^k)}\leq\frac1k,
			\qquad
			\nnorm{\lambda^{k,\iota(k)}-\lambda^k}\leq\frac1k, 
			\\
			\label{eq:diag_sequence_quality1.5}
			&\sgn(\lambda_i^{k,\iota(k)}) = \sgn(\lambda_i^k), \quad \forall i\in I,
			\\
			\label{eq:diag_sequence_quality2}
			&F(x^k)-\delta^k \notin \Gamma_j \Longrightarrow y^{k,\iota(k)} \notin \Gamma_j,
			\quad  \forall j\in [t],
			\\
			\label{eq:diag_sequence_quality3}
			&F_i(x^k)-\delta_i^k \in ( a_i^j,b_i^j) 
			\Longrightarrow y_i^{k,\iota(k)} \in ( a_i^j,b_i^j),
			\quad  \forall j\in J(x^k,\delta^k),\, \forall i\in[\ell],
			\end{align}
		\end{subequations}
		where \eqref{eq:diag_sequence_quality2} and \eqref{eq:diag_sequence_quality3}
		can be fulfilled as the complement of $\Gamma_j$ and $(a_i^j, b_i^j)$,
		$j\in[t]$, $i\in[\ell]$, are open, respectively.
		As in the proof of \cref{prop:equivalence_of_stationarities}, 
		we now define a sequence 
		$
			\{ (x^k, \tilde\lambda^k, \tilde\delta^k, \tilde\varepsilon^k) \}_{k=1}^\infty
			\subseteq \mathbb{R}^{n+\ell+\ell+n}
		$
		by means of \eqref{eq:some_surrogate_sequences}.
		This sequence fulfills \eqref{eq:properties_surrogate_sequences}
		and, by construction, $(x^k,\tilde\delta^k,\tilde\varepsilon^k)\to(\bar x,0,0)$.
		Hence, it is approximately S-stationary w.r.t. $\bar x$.
		
		We will now show that $I^\exists(x^k,\tilde\delta^k)\subset I^\exists(x^k,\delta^k)$
		holds for each $k\in\N$.
		Therefore, fix $i\in I^\exists(x^k,\tilde\delta^k)$. 
		By definition, there is some
		$j\in J(x^k,\tilde\delta^k)$, giving $y^{k,\iota(k)}\in\Gamma_j$, 
		such that $y_i^{k,\iota(k)}\in \{a_i^j,b_i^j\}$.
		Due to \eqref{eq:diag_sequence_quality2}, $j\in J(x^k,\delta^k)$
		follows, and \eqref{eq:diag_sequence_quality3} yields
		$F_i(x^k)-\delta_i^k\in \{a_i^j,b_i^j\}$.
		Hence, $i\in I^\exists(x^k,\delta^k)$ follows.
		
		Picking a subsequence if necessary, we may assume that
		the index sets $I^\exists(x^k,\tilde\delta^k)$ are the same for all $k\in\N$,
		see \cref{rem:dp_seq_const}.
		Thus, let us set $\tilde I\coloneqq I^\exists(x^k,\tilde\delta^k)$ for any $k\in\N$.
		Then we have $\tilde I\subset I$, and, due to \eqref{eq:diag_sequence_quality1.5},
		$\lambda^k_i$ and $\tilde\lambda^k_i$ have the same sign for each $i\in \tilde I$
		and $k\in\N$.	
		Hence, it finally follows that 
		item~\ref{subMFC_ODP_I} of \cref{def:subMFC_ODP} 
		also holds when using $\tilde I$
		instead of $I$.
\end{proof}

Let us point out that condition ODP-subMFC from \cref{def:subMFC_ODP}
is not a constraint qualification in the narrower sense as it depends,
via the required existence of an approximately M-stationary sequence,
also on the objective function of \eqref{eq:orthodisjunctive_problem}.
We, thus, refer to it as a qualification condition.
The most popular qualification condition in the literature we are aware of
is the so-called strict Mangasarian--Fromovitz condition,
which guarantees uniqueness of multipliers in situations where
stationarity of the underlying feasible point is already known,
see \cite{Wachsmuth2013} for a discussion in standard nonlinear programming
and \cite[Theorem~2.2]{Shapiro1997} for the more general setting of
conically-constrained optimization in arbitrary Banach spaces.

The subsequently stated result justifies our interest in ODP-subMFC.

\begin{theorem}\label{thm:ODP_subMFC}
	Let $\bar x\in\R^n$ be an AM-stationary point of \eqref{eq:orthodisjunctive_problem}, and let
	ODP-subMFC be satisfied at $\bar{x}$ with
	$\{(x^k,\lambda^k,\delta^k,\varepsilon^k)\}_{k=1}^\infty\subset\R^{n+\ell+\ell+n}$
	being the involved approximately M-stationary sequence w.r.t $\bar{x}$.
	Then the following assertions hold.
	\begin{enumerate}
		\item\label{item:ODP_subMFC_AM} The point $\bar x$ is M-stationary.
		\item\label{item:ODP_subMFC_AS} 
			If $\{(x^k,\lambda^k,\varepsilon^k)\}_{k=1}^\infty$ is a
			strictly approximately S-stationary sequence w.r.t. $\bar{x}$,
			then the point $\bar x$ is S-stationary.
	\end{enumerate}
\end{theorem}
\begin{proof}
	Let us start with the proof of the first assertion.
		By item \ref{subMFC_ODP_II} of \cref{def:subMFC_ODP}, 
		we find \eqref{eq:am_dp_lambda},
		the convergences \eqref{eq:as_dp_conv_standard},
		and
		\begin{equation}\label{eq:odp_proof1}
			\varepsilon^k=\nabla f(x^k)+ \sum_{i\in I} \lambda_i^k \nabla F_i(x^k)
		\end{equation}
		for all $k\in\N$, 
		where $I\subset I^\exists(\bar x)$ is the index set whose existence
		is claimed in \cref{def:subMFC_ODP}, fulfilling
		$I=I^\exists(x^k,\delta^k)$ for all $k\in\N$.
		Above, we used \eqref{eq:trivial_components_limiting_normal_cone} 
		to obtain $\lambda_{[\ell]\setminus I}^k = 0$, while
		taking into account \cref{rem:dp_seq_const} if necessary. 
		If $I=\emptyset$, the continuity of $\nabla f$ implies 
		$0 = \nabla f(\bar{x})$	when taking the limit $x^k\to\bar{x}$ in \eqref{eq:odp_proof1}.
		As $\bar{x}$ is feasible to \eqref{eq:orthodisjunctive_problem}, it is already
		M-stationary with multiplier	$\lambda \coloneqq 0$.
		Let us note that this also proves the second assertion in this trivial case,
		which, thus, is omitted later on.			

		In the case where $I \neq \emptyset$, 
		we show that the sequence 
		$\{ \lambda^k_{I} \}_{k=1}^\infty$ is
		bounded. To this end, let us assume the contrary, 
		i.e., that this sequence is unbounded.
		Then we have 
		$\nnorm{ \lambda^k_{I} } \to \infty$
		along a subsequence (without relabeling)
		and, hence, without loss of generality, the convergence of the bounded 
		sequence
		\begin{equation*}%\label{eq:proof3.5}  
			\frac{ \lambda^k_{I} }
				{\nnorm{ \lambda^k_{I} }} 
			\to 
			\bar{\lambda}_{I}
		\end{equation*} 
		for some $\bar\lambda\in\R^\ell$ such that $\bar\lambda_I\neq 0$.
		Dividing \eqref{eq:odp_proof1} by 
		$\nnorm{ \lambda^k_{I} }$ 
		and passing to the limit now yields
		\begin{equation}\label{eq:odp_proof4}
			0 = \sum_{i\in I} \bar{\lambda}_i \nabla F_i(\bar{x}).
		\end{equation}  
		Defining $u\in\R^{\ell}$ by means of
		$u_i \coloneqq |\bar{\lambda}_i|\geq 0$ for each $i\in I$
		and $u_{[\ell]\setminus I}\coloneqq 0$,
		and exploiting the set $I_\pm \coloneqq\{i\in I\,|\,u_i\neq 0\}$,
		\eqref{eq:odp_proof4} is equivalent to
		\begin{equation*}
			0 
			= 
			\sum_{i\in I} \sgn(\bar{\lambda}_i) u_i \nabla F_i(\bar{x})
			=
			\sum_{i\in I_\pm} \sgn(\lambda^k_i) u_i \nabla F_i(\bar x)
			=
			\sum_{i\in I} \sgn(\lambda^k_i) u_i \nabla F_i(\bar x).
		\end{equation*}
		This is a contradiction to item~\ref{subMFC_ODP_I} of \cref{def:subMFC_ODP} 
		as $u\geq 0$, $u_{[\ell]\setminus I}=0$, and $u_I\neq 0$ due to $\bar\lambda_{I}\neq 0$.
		Hence, the sequence 
		$\{ \lambda^k_{I} \}_{k=1}^\infty$ 
		has to be bounded, and without loss of 
		generality, 
		we can pick $\widetilde\lambda\in\R^\ell$ 
		such that $\widetilde\lambda_{[\ell]\setminus I} = 0$ and 
		$ \lambda^k_I \to \widetilde{\lambda}_{I} $. 
		Taking the limit $k\to\infty$, it follows from 
		\eqref{eq:odp_proof1}, $(x^k,\delta^k,\varepsilon^k)\to(\bar x,0,0)$,
		and continuous differentiability of all involved functions that
		\[
			0
			=
			\nabla f(\bar x)
			+
			\sum_{i\in I}\widetilde\lambda_i\nabla F_i(\bar x).
		\] 
		Due to \eqref{eq:am_dp_lambda} and 
		the robustness of the limiting normal cone,
		we obtain $\widetilde\lambda \in N_\Gamma(F(\bar{x}))$.
		Hence, we have shown that $\bar x$ is an M-stationary point of 
		\eqref{eq:orthodisjunctive_problem}.

		Let us now prove the second assertion.
		As mentioned earlier, let us directly address the nontrivial case,
		i.e., $I\neq\emptyset$.
		From the first assertion we already know that
		$\bar x$ is M-stationary with multiplier $\tilde\lambda$,
		which is the limit of $\{\lambda^k\}_{k=1}^\infty$ 
		(along a subsequence if necessary).
		As we have $\{\lambda^k\}_{k=1}^\infty\subset\widehat N_\Gamma(F(\bar x))$
		by definition of a strictly approximately S-stationary sequence 
		w.r.t. $\bar{x}$,
		the closedness of $\widehat N_\Gamma(F(\bar x))$ already yields
		$\tilde\lambda\in\widehat N_\Gamma(F(\bar x))$,
		i.e., $\bar x$ is already S-stationary with multiplier $\tilde\lambda$.
\end{proof}

The above proof provides the following observation that parallels 
\cite[Remark~3.10]{KaemingFischerZemkoho2024}.
\begin{remark}
	Let	ODP-subMFC be satisfied at $\bar{x}$ with
	$\{(x^k,\lambda^k,\delta^k,\varepsilon^k)\}_{k=1}^\infty\subset\R^{n+\ell+\ell+n}$
	being the involved approximately M-stationary sequence w.r.t $\bar{x}$.
	The proof of \cref{thm:ODP_subMFC} shows that, if ODP-subMFC holds at $\bar{x}$ 
	with $I\neq\emptyset$, the sequence $\{ \lambda_I^k \}_{k=1}^\infty$ is bounded.
	\\
	In fact, we even have that $\{\lambda^k\}_{k=1}^\infty$ must be bounded
	in the present situation.
	Indeed, supposing that $\{\lambda^k_i\}_{k=1}^\infty$ is unbounded for some
	$i\in[\ell]$, we find $i\in I^\exists(x^k,\delta)=I$ for
	all sufficiently large $k\in\N$ from \eqref{eq:trivial_components_limiting_normal_cone}
	and condition \ref{subMFC_ODP_II} of \cref{def:subMFC_ODP},
	and, thus, we end up with a contradiction.\\
	Conversely, it can be checked by means of simple examples that not every sequence 
	$\{(x^k,\lambda^k,\delta^k,\varepsilon^k)\}_{k=1}^\infty\subset\R^{n+\ell+\ell+n}$
	that is approximately M-stationary w.r.t $\bar{x}$ such that 
	$\{ \lambda_I^k \}_{k=1}^\infty$ is bounded 
	necessarily constitutes a suitable candidate for verifying ODP-subMFC at $\bar{x}$, 
	see e.g.\ \cref{ex:ODP-subMFC_Mstat_not_equivalent} below.
\end{remark}

The combination of \cref{cor:local_minimizers_AM_stat} and \cref{thm:ODP_subMFC} 
particularly yields the following.

\begin{corollary}
	If $\bar{x}$ is a local minimizer of \eqref{eq:orthodisjunctive_problem} and
	ODP-subMFC holds at $\bar{x}$, then $\bar{x}$ is M-stationary.
\end{corollary}

As ODP-subMFC depends on the objective function, naturally the question 
arises whether it may be weak enough to simply be an equivalent characterization 
of M-stationarity. The following example demonstrates that this is not generally 
the case as the presence of equality constraints, for example, may result in M-stationarity not 
necessarily implying ODP-subMFC.

\begin{example}\label{ex:ODP-subMFC_Mstat_not_equivalent}
	Let us consider the optimization problem
	\[
		\min\limits_x\quad x\quad\text{s.t.}\quad 
		\left(\frac{1}{2}x^2,-x\right)\in\Gamma\coloneqq\Gamma_1\coloneqq \{ 0 \} \times \{ 0 \}.
	\]
	The only feasible point is $\bar x\coloneqq 0$, which, thus, is the uniquely determined
	global minimizer of the problem.\\
	To verify M-stationarity of $\bar{x}$, condition \eqref{eq:M_stat_DP_x} is given by 
	$0 = 1 + \lambda_1 \bar{x} - \lambda_2  = 1 - \lambda_2$				
	and needs to be fulfilled for some 
	$(\lambda_1, \lambda_2) \in N_\Gamma((0,0))= \R\times \R$ 
	due to \eqref{eq:M_stat_DP_lambda}.
	Thus, $\bar{x}$ is M-stationary, 
	and, due to $\widehat N_\Gamma((0,0))= N_\Gamma((0,0))$,
	even S-stationary with multiplier 
	$\lambda = (\lambda_1,1)$ for all $\lambda_1\in\R$. \\
	Concerning ODP-subMFC at $\bar{x}$, consider any sequence
	$\{ (x^k, \lambda^k, \delta^k, \varepsilon^k) \}_{k=1}^\infty$ 
	that is approximately M-stationary w.r.t. $\bar{x}$, i.e., that 
	fulfills the convergences \eqref{eq:as_dp_conv_standard} as well 
	as the conditions \eqref{eq:as_dp_x} and \eqref{eq:am_dp_lambda}, amounting to 
	$\varepsilon^k = 1 + \lambda_1^kx^k - \lambda_2^k$ and 
	\[
		(\lambda_1^k,\lambda_2^k) 
		\in 
		\begin{cases}
		 	\R \times \R & (\delta_1^k,\delta_2^k) = \left( \frac{1}{2} (x^k)^2, -x^k \right), \\
		 	\emptyset & (\delta_1^k,\delta_2^k) \neq \left( \frac{1}{2} (x^k)^2, -x^k \right),
		\end{cases}			
	\]
	for all $k\in\N$. Thus, any approximately M-stationary sequence w.r.t. $\bar x$ 
	that could be used in item \ref{subMFC_ODP_II} of \cref{def:subMFC_ODP} fulfills 
	$(\delta_1^k,\delta_2^k) = \left( \frac{1}{2} (x^k)^2, -x^k \right)$, which 
	always implies $J(x^k,\delta^k) = \{1\}$ and $I^\exists(x^k,\delta^k)=\{1,2\}$. 
	Consequently, to verify ODP-subMFC at $\bar{x}$, we necessarily have to choose
	$I=I^\exists(x^k,\delta^k)=\{1,2\}$. Accordingly, item \ref{subMFC_ODP_I} of \cref{def:subMFC_ODP}
	requires for all $u\in\R^2_+$ with $u\neq (0,0)$ that
	\[
		0 \neq \sum_{i\in\{1,2\}} \sgn(\lambda_i^k) u_i \nabla F_i(\bar{x}) 
		= \sgn(\lambda_1^k) u_1 \bar{x} - \sgn(\lambda_2^k) u_2
		= -\sgn(\lambda_2^k) u_2,
	\]
	which is clearly violated as can be confirmed by taking, for example, $u=(1,0)$. 
	As the above observations hold for 
	any approximately M-stationary sequence w.r.t. $\bar x$ that could be used
	in item \ref{subMFC_ODP_II} of \cref{def:subMFC_ODP}, it is clear that ODP-subMFC
	at $\bar{x}$ cannot be satisfied.
\end{example}

However, for certain subclasses of \eqref{eq:orthodisjunctive_problem}, 
ODP-subMFC does in fact provide an equivalent characterization of M-stationarity, as shown in 
the next theorem. In \cref{sec:NLPs}, a specific subclass fulfilling the assumptions of this 
theorem will be examined.

\begin{theorem}\label{thm:ODP-subMFC_for_ODP-equi}
Let $\bar x\in\R^n$ be an AM-stationary point of \eqref{eq:orthodisjunctive_problem}
such that $J(\bar x) = \{ j_0 \}$ for some $j_0 \in [t]$.
If $\intr(\Gamma_{j_0}) \neq \emptyset$,
then ODP-subMFC holds at $\bar{x}$ if and only if $\bar{x}$ is M-stationary.
\end{theorem}

\begin{proof}
	$[\Longrightarrow]$:
		This implication follows directly from \cref{thm:ODP_subMFC}~\ref{item:ODP_subMFC_AM}.\\
	$[\Longleftarrow]$:
		Suppose that $\bar{x}$ is M-stationary, i.e., there exists a multiplier $\lambda \in\R^\ell$ such that
		\[ 
			-\nabla f(\bar{x}) = \sum_{i\in I^\exists(\bar{x})} \lambda_i \nabla F_i(\bar{x})
		\]
		and $\lambda\in N_\Gamma(F(\bar{x}))$ hold, where we used in the above equation that 
		$\lambda\in N_\Gamma(F(\bar{x}))$ implies 
		$\lambda_{[\ell]\setminus I^\exists(\bar{x})} = 0$ by \eqref{eq:trivial_components_limiting_normal_cone}.
		By \cite[Lemma~1]{AndreaniHaeserSchuverdtSilva2012}, we know that there exist
		$I \subset I^\exists(\bar{x})$ and $\widetilde{\lambda}_i\in\R\setminus\{0\}$ 
		with 
		$\sgn(\widetilde{\lambda}_i)=\sgn(\lambda_i)$
		for all $i\in I$ 
		such that the vectors $\{ \nabla F_i(\bar{x}) \}_{i\in I}$ are linearly 
		independent and
		\begin{equation}\label{eq:ODP-subMFC_Mstat_caratheodory}
			-\nabla f(\bar{x}) = \sum_{i\in I} \widetilde\lambda_i \nabla F_i(\bar{x})
		\end{equation}
		holds.
		Due to $J(\bar{x}) = \{ j_0\}$ for some $j_0 \in [t]$, we know that 
		\[	
			F(\bar{x}) \in \Gamma_{j_0} =  \prod\limits_{i=1}^\ell [a_i^{j_0}, b_i^{j_0}].
		\]
		For all $i\in [\ell]$, set $\eta_i \coloneqq \min \{ 1, b_i^{j_0} - a_i^{j_0}\}/2$, where
		we have $\eta_i > 0$ for all $i\in [\ell]$ due to $\intr(\Gamma_{j_0}) \neq \emptyset$.
		We now define the sequence 
		$\{ (x^k, \lambda^k, \delta^k, \varepsilon^k) \}_{k=1}^\infty $
		given by
		\[
			x^k\coloneqq \bar{x}, \quad
			\lambda_i^k\coloneqq\begin{cases}
				\widetilde\lambda_i & i\in I,\\
				0 & i\in[\ell]\setminus I,
			\end{cases} \quad
			\delta_i^k \coloneqq\begin{cases}
				0 & i\in I \cup ([\ell]\setminus I^\exists(\bar x)),\\
				\frac{\eta_i}{k} & i\in I^\exists(\bar x)\setminus I\colon F_i(\bar x) = b_i^{j_0}, \\
				-\frac{\eta_i}{k} & i\in I^\exists(\bar x)\setminus I\colon F_i(\bar x) = a_i^{j_0},
			\end{cases}\quad
			\varepsilon^k \coloneqq 0
		\]
		for all $k\in\N$.
		For this sequence, it can be checked that $J(x^k,\delta^k) = \{ j_0 \}$ and $I=I^\exists(x^k,\delta^k)$
		are valid for $k\in\N$ large enough. Without loss of generality, in the following,
		we refer by $\{ (x^k, \lambda^k, \delta^k, \varepsilon^k) \}_{k=1}^\infty $ to the corresponding 
		subsequence for which the latter properties hold for all $k\in\N$.
		As the vectors $\{ \nabla F_i(\bar{x}) \}_{i\in I}$ are linearly 
		independent and $\sgn(\lambda_i^k) = \sgn(\widetilde \lambda_i) \neq 0$ 
		for all $i\in I$ and $k\in\N$, 
		we clearly see that either $I = \emptyset$ holds or 
		\eqref{eq:ODP-subMFC_ODP_lin_indep} is satisfied for all 
		$u\in\R^{\ell}\setminus\{0\}$ with $u\geq 0$ and $u_{[\ell]\setminus I}=0$, 
		meaning that item \ref{subMFC_ODP_I} of \cref{def:subMFC_ODP} is satisfied. 
		Thus, applying \cref{prop:ODPsubMFC_via_AS_stat} shows that
		it remains to prove the approximate S-stationarity of 
		$\{ (x^k, \lambda^k, \delta^k, \varepsilon^k) \}_{k=1}^\infty$ 
		w.r.t. $\bar{x}$ to verify ODP-subMFC at $\bar{x}$.
		To this end, we first see that \eqref{eq:ODP-subMFC_Mstat_caratheodory} 
		implies \eqref{eq:as_dp_x}. Moreover, we need to show \eqref{eq:as_dp_lambda}, 
		i.e., $\lambda^k \in \widehat N_\Gamma (F(x^k)-\delta^k)$ for all $k\in \N$.
		It follows from \eqref{eq:limiting_normals_ODP}, $J(\bar x)= \{j_0\}$, and 
		$\lambda \in N_\Gamma (F(\bar{x}))$ that 
		\[
			\lambda \in %N_\Gamma (F(\bar{x}))
			%\subset 
			\prod_{i=1}^\ell \widehat{N}_{[a_i^{j_0},b_i^{j_0}]}(F_i(\bar{x})),
			\quad
			\widehat{N}_{[a_i^{j_0},b_i^{j_0}]}(F_i(\bar{x}))
			= 
			\begin{cases}
				\R_+ & i\in[\ell]\colon F_i(\bar{x}) = b_i^{j_0},\\
				\R_- & i\in[\ell]\colon F_i(\bar{x}) = a_i^{j_0},\\
				\{ 0 \} & i\in [\ell]\colon F_i(\bar{x}) \in (a_i^{j_0}, b_i^{j_0}).
			\end{cases}
		\]
		Due to $\sgn(\lambda_i^k) = \sgn(\widetilde{\lambda}_i) = \sgn(\lambda_i)$ 
		for all $i\in I$ and $\lambda_i^k = 0$ for all $i\in [\ell] \setminus I$,
		the above clearly implies 
		\[
			\lambda_i^k \in 
			\begin{cases}
				\R_+ & i \in I\colon F_i(\bar{x}) = b_i^{j_0},\\
				\R_- & i \in I\colon F_i(\bar{x}) = a_i^{j_0},\\
				\{ 0 \} & i\in [\ell] \setminus I
			\end{cases}
		\]
		for all $i\in [\ell]$.
		Moreover, for all $k\in\N$, it follows from $I=I^\exists(x^k,\delta^k)$ and the definition
		of $\delta^k$ that $F_i(\bar{x}) =  F_i(x^k)-\delta_i^k \in \{a_i^{j_0},b_i^{j_0}\}$ 
		for all $i\in I$
		and $F_i(x^k) - \delta_i^k \in (a_i^{j_0}, b_i^{j_0})$ 
		for all $i\in [\ell] \setminus I$.
		This yields
		\begin{align*}
			\lambda_i^k 
			&\in 
			\begin{cases}
				\R_+ & i\in[\ell]\colon F_i(x^k) - \delta_i^k = b_i^{j_0},\\
				\R_- & i\in[\ell]\colon F_i(x^k) - \delta_i^k = a_i^{j_0},\\
				\{ 0 \} & i\in[\ell]\colon F_i(x^k)-\delta_i^k \in (a_i^{j_0}, b_i^{j_0})
			\end{cases}
			\\
			&= 
			\widehat{N}_{[a_i^{j_0},b_i^{j_0}]}(F_i(\bar{x})-\delta_i^k)
		\end{align*}
		for all $i\in [\ell]$.
		Using \eqref{eq:regular_normals_ODP} and $J(x^k,\delta^k) = \{ j_0\}$,
		we obtain
		\[
			\lambda^k
			\in
			\prod_{i=1}^\ell 
			\widehat {N}_{[a_i^{j_0},b_i^{j_0}]}(F_i(\bar{x})-\delta_i^k)
			=
			\widehat N_\Gamma(F(x^k)-\delta^k)
		\]
		for all $k\in \N$. Thus, it is shown that ODP-subMFC holds at $\bar{x}$.
\end{proof}

Recall from \cref{ex:AM_vs_AS_regularity} that AM- and AS-regularity are independent 
conditions. 
The next proposition shows that ODP-subMFC is independent of both of these 
independent conditions and, thus, can particularly be useful in situations where
both AM- and AS-regularity are violated. 

\begin{proposition}\label{rem:ODP-subMFC_indep_AM-AS-reg}
	Let $\bar{x}$ be feasible for \eqref{eq:orthodisjunctive_problem}.
	Then ODP-subMFC at $\bar{x}$ is independent of AM-regularity and of
	AS-regularity at $\bar{x}$.
\end{proposition}	

\begin{proof}
	As both AM- and AS-regularity reduce to the cone-continuity property from
	\cite{AndreaniMartinezRamosSilva2016}
	in the special case $\Gamma\coloneqq\R^\ell_-$,
	we can use
	\cite[Example~A.3]{KaemingFischerZemkoho2024} to see that ODP-subMFC 
	can hold at a point that is neither AM- nor AS-regular.

	On the opposite, let us again consider the optimization problem from
	\cref{ex:ODP-subMFC_Mstat_not_equivalent}, for which we have shown that
	ODP-subMFC does not hold at the global minimizer $\bar x = 0$.
	Moreover, using
	$\widehat N_\Gamma(F(\bar{x})) =  N_\Gamma(F(\bar{x})) = \R \times \R$,
	it follows that
	$F'(\bar{x})^\top \widehat N_\Gamma(F(\bar{x})) 
	= F'(\bar{x})^\top N_\Gamma(F(\bar{x})) = \R$,
	which trivially guarantees that $\bar{x}$ is both AM- and AS-regular. 	
\end{proof}

Let us comment on the value of assertion~\ref{item:ODP_subMFC_AS} in \cref{thm:ODP_subMFC}
in more detail.
For that purpose, pick a local minimizer $\bar x\in\R^n$ of \eqref{eq:orthodisjunctive_problem}.
As explained in \cref{sec:ODP},
a potential way in order to obtain S-stationarity of $\bar x$ is to guarantee SAS-stationary
and AS-regularity of $\bar x$.
However, this approach has two drawbacks.
First, SAS-stationarity is not a necessary optimality condition for \eqref{eq:orthodisjunctive_problem}
without additional assumptions. Second, AS-regularity requires to check the behavior of
infinitely many sequences.
\cref{thm:ODP_subMFC}\,\ref{item:ODP_subMFC_AS} offers another way to justify S-stationarity of $\bar x$.
Indeed, from \cref{cor:local_minimizers_AM_stat} we know that $\bar x$ is AM-stationary,
and we can check validity of ODP-subMFC when using the associated approximately M-stationary sequence
w.r.t.\ $\bar x$, which would ensure M-stationarity of $\bar x$.
If, by chance, the involved sequence is already a strictly approximately S-stationary sequence w.r.t.\ $\bar x$,
then we have already confirmed S-stationarity of $\bar{x}$. 
It follows from \cref{rem:ODP-subMFC_indep_AM-AS-reg} that this approach is not only beneficial
whenever the verification of AM- and AS-regularity is too difficult with regard to the computation itself, 
but especially in situations where AM- and AS-regularity are actually violated.

Thus, starting from a given approximately M-stationary sequence w.r.t.\ a
given feasible point for \eqref{eq:orthodisjunctive_problem},
ODP-subMFC provides a simple criterion to infer M-stationarity (or even S-stationarity) of this point. 
However, one has to note the obvious fact that ODP-subMFC may hold for some but not
necessarily all approximately M-stationary sequences w.r.t.\ the reference point.
Indeed, item~\ref{subMFC_ODP_II} of \cref{def:subMFC_ODP} merely requires the existence
of such a sequence, see \cref{ex:AM-reg_fails_at_M_stat_point} below.

Therefore, to close the paragraph, 
let us draft a practically relevant situation where the observations from
\cref{thm:ODP_subMFC} can be applied beneficially.
It is well known from the literature, see e.g.\ 
\cite{AndreaniHaeserSeccinSilva2019,JiaKanzowMehlitzWachsmuth2023,Ramos2021},
that several solution methods when applied to (subclasses of) problem \eqref{eq:orthodisjunctive_problem}
generate approximately M-stationary sequences w.r.t.\ accumulation points
(if not terminating finitely).
In order to infer an accumulation point's M-stationarity, it is typically claimed in the literature 
that AM-regularity or a stronger constraint qualification holds.
However, as AM-regularity requires control over infinitely many sequences,
this might still be too restrictive to be applicable.
Condition ODP-subMFC now offers a distinct approach to recognize M-stationarity.
Indeed, in the present situation, item~\ref{subMFC_ODP_II} of \cref{def:subMFC_ODP} is already valid,
so it remains to check item~\ref{subMFC_ODP_I} of \cref{def:subMFC_ODP}.
If it holds, M-stationarity of the accumulation point under consideration follows.
Note, however, that it may happen that the considered accumulation point is M-stationary
while the generated approximately M-stationary sequence does not enable a confirmation via
\cref{def:subMFC_ODP}~\ref{subMFC_ODP_I}.

Let us visualize this by means of a simple example
which has been inspired by \cite[Example~5.5]{BenkoMehlitz2025}.
\begin{example}\label{ex:AM-reg_fails_at_M_stat_point}
	Let us consider the optimization problem
	\[
		\min\limits_x\quad x\quad\text{s.t.}\quad (x,-x^2)\in\Gamma\coloneqq\Gamma_1\cup\Gamma_2
	\]
	with
	\[
		\Gamma_1\coloneqq \R_+\times\R,\qquad
		\Gamma_2\coloneqq \R\times\R_+.
	\]
	Its feasible set equals $\R_+$, such that $\bar x\coloneqq 0$ is the uniquely determined
	global minimizer of the problem.
	\\
	Let us assume that, for the numerical solution of the problem,
	a multiplier-penalty method has been applied which generated the sequence
	$\{(x^k,\lambda^k,\delta^k,\varepsilon^k)\}_{k=1}^\infty\subset\R^{1+2+2+1}$
	given by
	\[
		x^k\coloneqq -\frac{1}{k},
		\quad
		\lambda^k\coloneqq (-1,0),
		\quad
		\delta^k\coloneqq\left(-\frac1k,0\right),
		\quad
		\varepsilon^k\coloneqq 0,
		\quad
		\forall k\in\N.
	\]
	As we have
	\[
		N_\Gamma\left(\left(0,-\frac{1}{k^2}\right)\right)=\R_-\times\{0\},
	\]
	one can easily check that $\{(x^k,\lambda^k,\delta^k,\varepsilon^k)\}_{k=1}^\infty$ 
	is an approximately M-stationary sequence w.r.t.\ $\bar x$.
	Note that $\{x^k\}_{k=1}^\infty$ converges to $\bar x$ from outside the feasible set.
	It is not hard to check that AM-regularity does not hold for the given problem
	at $\bar x$, see \cref{ex:AM_vs_AS_regularity}.
	However, for the above sequence, we find $I^\exists(x^k,\delta^k) = \{1\}$ for all $k\in\N$, 
	such that, using this sequence in ODP-subMFC and, thus, $I\coloneqq \{1\}$,
	ODP-subMFC is valid if $	0 \neq -u(1,0)^\top$ holds for all $u>0$.
	As the latter is clearly fulfilled, ODP-subMFC holds, and it follows 
	that $\bar{x}$ must be M-stationary. 

	To close, suppose that another run of the multiplier-penalty method produces the sequence
	$\{(\tilde x^k,\tilde \lambda^k,\tilde \delta^k,\tilde\varepsilon^k)\}_{k=1}^\infty\subset\R^{1+2+2+1}$
	given by
	\[
		\tilde x^k\coloneqq -\frac{1}{k},
		\quad
		\tilde\lambda^k\coloneqq \left(0,-\frac{k}{2}\right),
		\quad
		\tilde\delta^k\coloneqq\left(-\frac1k,-\frac{1}{k^2}\right),
		\quad
		\tilde\varepsilon^k\coloneqq 0,
		\quad
		\forall k\in\N.
	\]
	Again, it is not hard to see that 
	$\{(\tilde x^k,\tilde \lambda^k,\tilde \delta^k,\tilde\varepsilon^k)\}_{k=1}^\infty$
	is an approximately M-stationary sequence w.r.t.\ $\bar x$.
	We have $I^\exists(x^k,\delta^k)=\{1,2\}$
	and, thus, need to choose $I\coloneqq \{1,2\}$ when checking ODP-subMFC
	for this particular sequence. 
	Then condition \ref{subMFC_ODP_I} of \cref{def:subMFC_ODP}
	is clearly violated as $\sgn(\tilde\lambda_1^k) = 0$ and $\{1\}\subset I$ are valid.
	Hence, ODP-subMFC does not hold when considering 
	$\{(\tilde x^k,\tilde \lambda^k,\tilde \delta^k,\tilde\varepsilon^k)\}_{k=1}^\infty$
	as the approximately M-stationary sequence w.r.t.\ $\bar x$.
	As AM-regularity also fails, see above,
	M-stationarity of $\bar x$ remains undetected.
\end{example}

\begin{remark}\label{rem:algorithmic_consequences}
	We already mentioned that some prominent
	solution algorithms for problems of type \eqref{eq:orthodisjunctive_problem}
	turn out to produce approximately M-stationary sequences w.r.t.\ accumulation points,
	see e.g.\ \cite{AndreaniHaeserSeccinSilva2019,JiaKanzowMehlitzWachsmuth2023,Ramos2021}.
	Now, in order to infer that accumulation points of the primal iterates are
	already M-stationary, a qualification condition has to be exploited which
	controls the behavior of \emph{all} approximately M-stationary sequences
	w.r.t.\ the accumulation point under consideration,
	as one does not know a priori which of them will be produced by the algorithm.
	On the one hand, a typical example of such a qualification condition 
	is AM-regularity from \cref{def:asymptotic_regularity}.
	On the other hand, ODP-subMFC from \cref{def:subMFC_ODP} is not suitable for that
	purpose as it makes an assertion about \emph{one} particular
	approximately M-stationary sequence w.r.t.\ the accumulation point.
	\\
	In order to overcome this issue, one may formulate a stronger condition which,
	given some feasible point $\bar x\in\R^n$ of \eqref{eq:orthodisjunctive_problem},
	for each index set $I\subset I(\bar x)$ and each approximately M-stationary sequence w.r.t. $\bar{x}$
	satisfying condition~\ref{subMFC_ODP_II} of \cref{def:subMFC_ODP}, 
	demands that we also have condition~\ref{subMFC_ODP_I} of \cref{def:subMFC_ODP}.
	Clearly, this condition is much stronger than ODP-subMFC.
	In fact, 
	it can be shown by straightforward arguments that it is equivalent to the 
	\emph{generalized Mangasarian--Fromovitz constraint qualification}
	\[
		F'(\bar x)^\top\lambda=0,\,\lambda\in N_\Gamma(F(\bar x))
		\quad\Longrightarrow\quad\lambda=0,
	\]
	which is a standard constraint qualification in geometrically-constrained
	optimization and a comparatively strong condition.
	\\
	Depending on the algorithm,
	it might be possible to reasonably strengthen ODP-subMFC by postulating its validity
	for a certain, smaller subclass of approximately M-stationary sequences w.r.t.\ an accumulation point,
	namely the subclass of those sequences the algorithm generates.
	A precise analysis of this kind of condition and the investigation of
	its algorithmic consequences is, however, beyond the scope of this paper
	but should be addressed in some future research.
\end{remark}

\subsection{The subset Mangasarian--Fromovitz condition for problems with inequality constraints}\label{sec:NLPs}

In this subsection, we embed our findings from \cref{sec:ODP_stuff} into the results
obtained in \cite{KaemingFischerZemkoho2024}, where the subset
Mangasarian--Fromovitz condition has been introduced for
(merely Lipschitzian) optimization problems with pure inequality constraints.

Therefore, we apply \cref{def:subMFC_ODP} to the inequality-constrained problem
\begin{equation}\label{eq:inequality_problem}\tag{P}
	\min\limits_x \quad f(x) \quad \text{s.t.}\quad g(x)\leq 0
\end{equation}
with continuously differentiable functions $f\colon\R^n\to\R$
and $g\colon\R^n\to\R^m$ 
to show that the results in
\cite{KaemingFischerZemkoho2024},
at least in the continuously differentiable setting,
can be recovered.
Clearly, \eqref{eq:inequality_problem} can be modeled as a special case of
\eqref{eq:orthodisjunctive_problem} with
$t\coloneqq 1$, $\ell\coloneqq m$, and 
\[
	\Gamma = \Gamma_1 \coloneqq \prod_{i=1}^\ell [-\infty,0].
\]
As in \cite{KaemingFischerZemkoho2024}, we define 
\[
	I^g(x, \delta) \coloneqq \{ i \in [\ell] \mid g_i(x) = \delta_i\}
\]
for $x\in\R^n$ and $\delta\in\R^\ell$ such that $g(x)-\delta\in\Gamma$.
Then we obtain $J(x,\delta)=\{1\}$, $I^\exists(x,\delta)=I^g(x,\delta)$,
and
\begin{align*}
	\widehat{N}_\Gamma(g(x)-\delta) 
	=
	N_\Gamma(g(x)-\delta)
	=
	\left\{
		\lambda\in\R^\ell\,\middle|\,
		\begin{aligned}	
			&\forall i\in I^g(x,\delta)\colon\,\lambda_i\geq 0,				
			\\
			&\forall i\in[\ell]\setminus I^g(x,\delta)\colon\,\lambda_i = 0
		\end{aligned}
	\right\}.
\end{align*}
Recall from \cref{cor:SAS_stat_vs_AM_stat} that AM- and SAS-stationarity are the same
for \eqref{eq:inequality_problem}. Moreover, they recover the classical AKKT conditions for
inequality-constrained problems, see \cite[Proposition~3]{MovahedianPourahmad2024}, that are considered
in \cite{KaemingFischerZemkoho2024}.

Applying \cref{def:subMFC_ODP} to the situation at hand yields the following lemma.

\begin{lemma}\label{lem:ODP-subMFC_for_P}
Let $\bar x\in\R^n$ be an AM-stationary point of \eqref{eq:inequality_problem}. Then ODP-subMFC
holds at $\bar{x}$ if and only if there exist $I\subset I^g(\bar x)$,
where $I^g(\bar x)\coloneqq I^g(\bar x,0)$,
and a sequence $\{(x^k,\lambda^k,\delta^k,\varepsilon^k)\}_{k=1}^\infty\subset\R^{n+\ell+\ell+n}$
such that the following conditions hold.
\begin{enumerate}[label=(\roman*)]
	\item\label{item:classical_subMFC_i} Either $I=\emptyset$, 
		or
		it holds for all 
		$u\in\R^{\ell}\setminus\{0\}$ with $u\geq 0$ and $u_{[\ell]\setminus I}=0$ that
		\begin{equation}\label{eq:classical_subMFC_improved}
			0  
			\neq 
			\sum_{i\in I} u_i \nabla g_i(\bar x).
		\end{equation}
	\item\label{item:classical_subMFC_ii} The sequence 
		$\{ (x^k,\lambda^k,\delta^k,\varepsilon^k) \}_{k=1}^\infty$ 
		satisfies the convergences \eqref{eq:as_dp_conv_standard}
		and
		\begin{align*}
			&\varepsilon^k = \nabla f(x^k) + \sum_{i=1}^\ell \lambda^k_i\nabla g_i(x^k),
			\\
			&\lambda^k_i\geq 0,\quad \forall i\in I^g(x^k,\delta^k),
			\\
			&\lambda^k_i=0,\quad \forall i\in[\ell]\setminus I^g(x^k,\delta^k),
			\\
			&g(x^k) \leq \delta^k,
			\\
			&I=I^g(x^k,\delta^k)
		\end{align*}
		for all $k\in\N$.
\end{enumerate}
\end{lemma}

Note that in a straightforward application of ODP-subMFC to \eqref{eq:inequality_problem},
condition \eqref{eq:ODP-subMFC_ODP_lin_indep}
reads as 
\begin{equation}\label{eq:classical_subMFC}
	0  
	\neq 
	\sum_{i\in I} \sgn(\lambda_i^k) u_i \nabla g_i(\bar x),
\end{equation}	
which clearly implies \eqref{eq:classical_subMFC_improved}
as
\cref{lem:ODP-subMFC_for_P}~\ref{item:classical_subMFC_ii} yields
$\lambda^k_i\geq 0$ for each index $i\in I$ and for all $k\in\N$. 
To prove the opposite, i.e., that ODP-subMFC with \eqref{eq:classical_subMFC_improved} also 
implies ODP-subMFC with \eqref{eq:classical_subMFC}, we again find 
$\lambda^k_i\geq 0$ for each index $i\in I$ and for all $k\in\N$ from 
\cref{lem:ODP-subMFC_for_P}~\ref{item:classical_subMFC_ii}.
Define $I_1 \coloneqq \{i\in I \mid \lambda_i^k = 0 \text{ for all } k\in\N\}$, where
we recall \cref{rem:dp_seq_const}, and set 
$\tilde\delta^k_i \coloneqq \delta^k_i + \frac{1}{k}$ for all $i\in I_1$
and $\tilde\delta^k_i \coloneqq \delta^k_i$
for all $i\in I \setminus I_1$.
Then all conditions in \cref{lem:ODP-subMFC_for_P}
are fulfilled for
$\tilde{I} \coloneqq I^g(x^k,\tilde\delta^k) = I \setminus I_1$
and the sequence 
$\{ (x^k,\lambda^k,\tilde\delta^k,\varepsilon^k) \}$, which additionally fulfills 
$\sgn(\lambda_i^k) = 1$ for
all $i\in \tilde{I}$.
Hence, ODP-subMFC holds when replacing
\eqref{eq:classical_subMFC_improved} by \eqref{eq:classical_subMFC}.
 
We have, thus, shown that ODP-subMFC in the setting of \eqref{eq:inequality_problem}
precisely recovers \cite[Definition~3.7]{KaemingFischerZemkoho2024}.
This also reveals that \cref{thm:ODP_subMFC} applied to \eqref{eq:inequality_problem}
recovers \cite[Theorem~3.9]{KaemingFischerZemkoho2024} as the M-stationary (and S-stationary) 
points of \eqref{eq:inequality_problem} correspond to the Karush--Kuhn--Tucker (KKT) 
stationary points considered in \cite{KaemingFischerZemkoho2024}.

Moreover, we find the following additional result as a corollary of
\cref{thm:ODP-subMFC_for_ODP-equi}.

\begin{corollary}\label{thm:ODP-subMFC_for_P-equi}
Let $\bar x\in\R^n$ be an AM-stationary point of \eqref{eq:inequality_problem}.
Then ODP-subMFC holds at $\bar{x}$ if and only if $\bar{x}$ is M-stationary.
\end{corollary}

Thus, we obtain from \cref{thm:ODP-subMFC_for_P-equi}
that ODP-subMFC and the validity of the KKT conditions are 
even equivalent for \eqref{eq:inequality_problem}. 
Note that this is only valid because we restricted the functions 
used in \eqref{eq:inequality_problem} to be continuously differentiable. Indeed, 
if the functions in \eqref{eq:inequality_problem} were merely Lipschitz continuous, 
it is known that this equivalence between ODP-subMFC and the KKT conditions 
does \textit{not} persist, see \cite[Example~3.12]{KaemingFischerZemkoho2024}.

\section{Application to complementarity-constrained optimization}\label{sec:MPCC}

Let us consider so-called mathematical problems with complementarity constraints,
which are optimization problems of the form
\begin{equation}\label{eq:MPCC} \tag{MPCC}
	\min_{x} \quad f(x) \quad\text{s.t.}\quad g(x)\leq 0,\,\, h(x)=0,\,\ 0 \leq G(x) \perp H(x) \geq 0 
\end{equation} 
with continuously differentiable functions $f\colon\mathbb{R}^n \to \mathbb{R}$, 
$g\colon\mathbb{R}^n\to \mathbb{R}^m$, $h\colon\mathbb{R}^n\to\mathbb{R}^p$, and 
$G,H\colon\mathbb{R}^n\to \mathbb{R}^q$.
In order to avoid the standard situation of nonlinear programming,
we will assume that $q\geq 1$ holds throughout.

To start our considerations,
we embed \eqref{eq:MPCC} into the class of (ortho)disjunctive optimization problems.
Therefore, setting $\ell\coloneqq m+p+2q$, we may define a function
$F\colon\R^n\to\R^\ell$ and a set $\Gamma\subset\R^\ell$ by means of 
\begin{subequations}\label{eq:dp_for_mpcc}
	\begin{align}
	\label{eq:dp_for_mpcc1}   
	F &\coloneqq (g_1, \dotsc, g_m, h_1, \dotsc, h_p, -G_1, -H_1, \dotsc, -G_q, -H_q), 
	\\
	\label{eq:dp_for_mpcc2}   
	\Gamma &\coloneqq \mathbb{R}^m_- \times \{ 0 \}^p \times C^q,
	\end{align} 
\end{subequations}
where the closed set $C\subset\R^2$ is given as
\begin{equation}\label{eq:compl_angle}
	C \coloneqq \{ (a,b) \in\mathbb{R}^2 \mid a \leq 0, b \leq 0, ab = 0 \}.
\end{equation}
Observe that $C$ is the union of the two convex polyhedral sets $C_1\coloneqq\{0\}\times\R_-$
and $C_2\coloneqq \R_-\times\{0\}$, so $\Gamma$ from \eqref{eq:dp_for_mpcc2} can be represented
as the union of $2^q$ convex polyhedral sets.
Hence, with the setting \eqref{eq:dp_for_mpcc}, 
\eqref{eq:disjunctive_problem} corresponds to \eqref{eq:MPCC}.
To be more precise, by defining, for some multi-index $j\in\JJ\coloneqq\{1,2\}^q$, 
the convex polyhedral set
\begin{equation}\label{eq:dp_for_mpcc3}
	\Gamma_j \coloneqq \R^m_-\times\{0\}^p\times\prod\limits_{i=1}^q C_{j_i},
\end{equation}
we find $\Gamma=\bigcup_{j\in\JJ}\Gamma_j$.
We also observe that each of these components $\Gamma_j$, $j\in\JJ$,
is the Cartesian product of closed intervals, such that \eqref{eq:MPCC}
is indeed an orthodisjunctive problem of type \eqref{eq:orthodisjunctive_problem}.

For a point $\bar{x}\in\mathbb{R}^n$ that is feasible for \eqref{eq:MPCC}, we define the 
following index sets of active constraints:
\begin{equation*}
	\begin{aligned}
		I_{0+}(\bar{x}) 
		&\coloneqq 
		\{ i\in [q] \mid G_i(\bar{x}) = 0, H_i(\bar{x}) > 0 \}, \\
		I_{+0}(\bar{x}) 
		&\coloneqq 
		\{ i\in [q] \mid G_i(\bar{x}) > 0, H_i(\bar{x}) = 0 \}, \\
		I_{00}(\bar{x}) 
		&\coloneqq \{ i\in [q] \mid G_i(\bar{x}) = 0, H_i(\bar{x}) = 0 \}.
	\end{aligned}
\end{equation*} 
To shorten notation, we further define 
\[
	I^h \coloneqq [p],
	\qquad 
	I^G(\bar{x}) \coloneqq I_{00}(\bar{x}) \cup I_{0+}(\bar{x}), 
	\qquad 
	I^H(\bar{x}) \coloneqq I_{00}(\bar{x}) \cup I_{+0}(\bar{x}).
\]
For brevity of notation,
we are going to exploit the so-called MPCC-Lagrangian function
$\mathcal L^\textup{cc}\colon\R^n\times\R^\ell\to\R$ given by
\[
	\mathcal L^\textup{cc}(x,\lambda)
	\coloneqq
	f(x)+\lambda_g^\top g(x)+\lambda_h^\top h(x)-\lambda_G^\top G(x)-\lambda_H^\top H(x),
\]
where $\lambda=(\lambda_g,\lambda_h,\lambda_G,\lambda_H)$.
We are going to work with the subsequently defined stationary concepts for \eqref{eq:MPCC} 
that are frequently used in the literature, see e.g.\ \cite{Ye2005} for an overview.

\begin{definition}\label{def:stat_MPCC}
	A point $\bar{x}\in\mathbb{R}^n$ feasible for \eqref{eq:MPCC} is called 
	\begin{enumerate}
		\item \textit{weakly stationary (W-stationary)} if there exist multipliers 
		$\lambda=(\lambda_g,\lambda_h,\lambda_G,\lambda_H)\in\mathbb{R}^\ell$ such that
				\begin{align*} 
					&0 = \nabla_x\mathcal L^\textup{cc}(\bar x,\lambda), \\
					&(\lambda_g)_i \geq 0,\,(\lambda_g)_ig_i(\bar x)=0,
						\quad \forall i\in[m], \\
					&(\lambda_G)_i = 0,\quad \forall i\in I_{+0}(\bar{x}), \\
					&(\lambda_H)_i = 0,\quad \forall i\in I_{0+}(\bar{x}); 
				\end{align*}
		\item \textit{Clarke-stationary (C-stationary)} 
			if it is W-stationary with multipliers 
			$\lambda\in\mathbb{R}^\ell$ and, for all $i\in I_{00}(\bar{x})$,
			condition $(\lambda_G)_i(\lambda_H)_i \geq 0$ holds;
		\item \textit{Mordukhovich-stationary (M-stationary)} if it is W-stationary with 
			multipliers $\lambda\in\mathbb{R}^\ell$ and, for all $i\in I_{00}(\bar{x})$, 
			condition $\min\{(\lambda_G)_i,(\lambda_H)_i\}>0$ or 
			$(\lambda_G)_i(\lambda_H)_i=0$ holds;
		\item \textit{strongly stationary (S-stationary)} if it is W-stationary with 
			multipliers $\lambda\in\mathbb{R}^\ell$ and, for all $i\in I_{00}(\bar{x})$, condition 
			$\min\{(\lambda_G)_i,(\lambda_H)_i\} \geq 0$ holds.
	\end{enumerate}
\end{definition}

Let us emphasize that the notions of M- and S-stationarity 
for \eqref{eq:MPCC} from \cref{def:stat_MPCC}
precisely correspond to the notions of M- and S-stationarity for \eqref{eq:disjunctive_problem} 
from \cref{def:DP_exact_stationary} when applied to the setting \eqref{eq:dp_for_mpcc}, respectively.
In order to compute the appearing regular and limiting normal cone to $\Gamma$,
one can directly rely on the representation \eqref{eq:dp_for_mpcc2} and exploit the
product rule as normal cones to $C$ from \eqref{eq:compl_angle} are easy to compute.

The notion of C-stationarity for \eqref{eq:MPCC} is induced 
by reformulating the complementarity constraints
by means of the nonsmooth equations $\min\{G_i(x),H_i(x)\}=0$, $i\in[q]$, 
and applying Clarke's subdifferential construction for its generalized differentiation,
see e.g.\ \cite[Section~2.2]{ScheelScholtes2000} for an illustration.
The concept of W-stationarity corresponds to the KKT conditions of a
so-called tightened standard nonlinear optimization problem associated with \eqref{eq:MPCC}
at the point of interest.

It is common knowledge, and clear from the definitions, that the implications
\[
	\text{S-stationary}
	\quad\Longrightarrow\quad
	\text{M-stationary}
	\quad\Longrightarrow\quad
	\text{C-stationary}
	\quad\Longrightarrow\quad
	\text{W-stationary}
\]
hold w.r.t.\ a given feasible point for \eqref{eq:MPCC}.

\subsection{Revisiting approximate stationarity for complementarity-constrained problems}

Approximate stationarity concepts for \eqref{eq:MPCC} are known from the literature, 
see e.g.\ \cite{AndreaniHaeserSeccinSilva2019,Mehlitz2023,Ramos2021}, 
and of our concern in this subsection.
For the introduction of suitable concepts, we 
use a notation that is inspired by \cite{KaemingFischerZemkoho2024,QiWei2000} and well-suited 
for the derivation of a problem-tailored version of ODP-subMFC.
Later on, we show that our definitions are actually equivalent to 
the ones which can be found in \cite{AndreaniHaeserSeccinSilva2019,Mehlitz2023}.

For $x\in\mathbb{R}^n$ and $\delta=(\delta_g,\delta_h,\delta_G,\delta_H)\in\R^\ell$, we define the following index sets for later use:
\begin{equation*}
	\begin{aligned}
		I_{0+}(x,\delta) 
		& \coloneqq 
		\{ i\in [q]\mid G_i(x) = -(\delta_G)_i,\, H_i(x) > -(\delta_H)_i\},\\
		I_{+0}(x,\delta) 
		& \coloneqq 
		\{ i\in [q]\mid G_i(x) > -(\delta_G)_i,\, H_i(x) = -(\delta_H)_i\},\\
		I_{00}(x,\delta) 
		&\coloneqq
		\{ i\in [q]\mid G_i(x) = -(\delta_G)_i,\, H_i(x) = -(\delta_H)_i\},\\
		I^G(x,\delta) 
		&\coloneqq 
		I_{0+}(x,\delta)\cup  I_{00}(x,\delta),\\
		I^H(x,\delta) 
		&\coloneqq 
		I_{+0}(x,\delta)\cup I_{00}(x,\delta).
	\end{aligned}
\end{equation*}

\begin{definition}\label{def:aw_ac_am_as_stationary}
	A point $\bar{x}\in\mathbb{R}^n$ feasible for \eqref{eq:MPCC} is called 
	\begin{enumerate}
		\item\label{item:AWstat} \textit{approximately weakly stationary (AW-stationary)} 
			if there exists an approximately W-stationary sequence 
			$\{(x^k,\lambda^k,\delta^k, \varepsilon^k) \}_{k=1}^\infty
			\subset 
			\mathbb{R}^{n+\ell+\ell+n}$ 
			w.r.t.\ $\bar x$, where
			\begin{equation}\label{eq:decoupling_lambda_and_delta}
				\lambda^k = (\lambda^k_g,\lambda^k_h,\lambda^k_G,\lambda^k_H),
				\qquad 
				\delta^k = (\delta_g^k, \delta_h^k, \delta_G^k, \delta_H^k),
			\end{equation}
			i.e., $\{(x^k,\lambda^k,\delta^k,\varepsilon^k)\}_{k=1}^\infty$
			is such that the conditions 
			\begin{subequations}\label{eq:aw}
				\begin{align}
					\label{eq:aw1}   
					&\varepsilon^k 
					= 
					\nabla_x\mathcal L^{\textup{cc}}(x^k,\lambda^k),  
					\\
					\label{eq:aw2}
					&(\lambda_g^k)_i \geq 0,\quad g_i(x^k) \leq (\delta_g^k)_i,
						\quad (\lambda^k_g)_i (g_i(x^k)-(\delta_g^k)_i) = 0,
						\quad
						\forall i\in[m], 
					\\
					\label{eq:aw3}
					&h(x^k)=\delta^k_h,
					\\
					\label{eq:aw4}
					&G_i(x^k)\geq -(\delta^k_G)_i,\quad 
					(\lambda^k_G)_i(G_i(x^k)+(\delta^k_G)_i) = 0,
					\quad \forall i\in[q],
					\\
					\label{eq:aw5}
					&H_i(x^k)\geq -(\delta^k_H)_i,\quad 
					(\lambda^k_H)_i(H_i(x^k)+(\delta^k_H)_i) = 0,
					\quad \forall i\in[q],
					\\
					\label{eq:aw6}
					&(G_i(x^k)+(\delta^k_G)_i)(H_i(x^k)+(\delta^k_H)_i)=0,
					\quad \forall i\in[q]
				\end{align}
			\end{subequations} 
			are fulfilled for all $k\in\mathbb{N}$ and the convergences
			\eqref{eq:as_dp_conv_standard} hold;
		\item \textit{approximately Clarke-stationary (AC-stationary)}
			if there exists an approximately C-stationary sequence 
			$\{(x^k,\lambda^k,\delta^k, \varepsilon^k) \}_{k=1}^\infty
			\subset 
			\mathbb{R}^{n+\ell+\ell+n}$ 
			w.r.t.\ $\bar x$, where \eqref{eq:decoupling_lambda_and_delta},
			i.e., $\{(x^k,\lambda^k,\delta^k,\varepsilon^k)\}_{k=1}^\infty$
			is an approximately W-stationary sequence w.r.t.\ $\bar x$ such that
			\begin{equation}\label{eq:ac}
				(\lambda^k_G)_i(\lambda^k_H)_i\geq 0,
				\quad
				\forall i\in I_{00}(\bar x)
			\end{equation}
			holds for all $k\in\N$;
		\item\label{item:AMstat} \textit{approximately Mordukhovich-stationary (AM-stationary)}
			if there exists an approximately M-stationary sequence 
			$\{(x^k,\lambda^k,\delta^k, \varepsilon^k) \}_{k=1}^\infty
			\subset 
			\mathbb{R}^{n+\ell+\ell+n}$ 
			w.r.t.\ $\bar x$, where \eqref{eq:decoupling_lambda_and_delta},
			i.e., $\{(x^k,\lambda^k,\delta^k,\varepsilon^k)\}_{k=1}^\infty$
			is an approximately W-stationary sequence w.r.t.\ $\bar x$ such that
			\begin{equation}\label{eq:am}
				\min\{(\lambda^k_G)_i,(\lambda^k_H)_i\}>0 
				\text{ or } 
				(\lambda^k_G)_i(\lambda^k_H)_i=0 ,
				\quad
				\forall i\in I_{00}(\bar x)
			\end{equation}
			holds for all $k\in\N$;
		\item\label{item:ASstat} \textit{strictly approximately strongly stationary 
			(SAS-stationary)} 
			if there exists a strictly approximately S-stationary sequence 
			$\{(x^k,\lambda^k,\delta^k, \varepsilon^k) \}_{k=1}^\infty
			\subset 
			\mathbb{R}^{n+\ell+\ell+n}$ 
			w.r.t.\ $\bar x$, where \eqref{eq:decoupling_lambda_and_delta},
			i.e., $\{(x^k,\lambda^k,\delta^k,\varepsilon^k)\}_{k=1}^\infty$
			is an approximately W-stationary sequence w.r.t.\ $\bar x$ such that
			\begin{equation}\label{eq:as}
					\min\{(\lambda^k_G)_i,(\lambda^k_H)_i\}\geq 0,
					\quad
					\forall i\in I_{00}(\bar x)
			\end{equation}
			holds for all $k\in\N$.
	\end{enumerate}
\end{definition}

Given a feasible point $\bar x\in\R^n$ of \eqref{eq:MPCC}, we obtain the relations
\begin{align*}
	\text{SAS-stationary}
	\quad&\Longrightarrow\quad
	\text{AM-stationary}
	\quad\Longrightarrow\quad
	\text{AC-stationary}
	\quad\Longrightarrow\quad
	\text{AW-stationary}
\end{align*}
regarding the different types of approximate stationarity of $\bar x$
directly from \cref{def:aw_ac_am_as_stationary}.
Our first result in this section justifies the terminology in \cref{def:aw_ac_am_as_stationary}
in the light of \cref{def:DP_special_sequence,def:DP_stationary}.

\begin{lemma}\label{lem:stat_DP_vs_MPCC}
	Let $\bar x\in\R^n$ be a feasible point for \eqref{eq:MPCC}.
	Then the following assertions hold.
	\begin{enumerate}
		\item\label{lem:stat_DP_vs_MPCC_AM}
			A sequence $\{(x^k,\lambda^k,\delta^k,\varepsilon^k)\}_{k=1}^\infty\subset\R^{n+\ell+\ell+n}$
			is an approximately M-stationary sequence w.r.t.\ $\bar x$ in the sense of
			\cref{def:aw_ac_am_as_stationary} if and only if it is an approximately M-stationary sequence
			w.r.t.\ $\bar x$ in the sense of \cref{def:DP_special_sequence} when applied to
			\eqref{eq:disjunctive_problem} in the setting \eqref{eq:dp_for_mpcc}.\\
			Particularly, $\bar x$ is AM-stationary in the sense of \cref{def:aw_ac_am_as_stationary}
			if and only if $\bar x$ is AM-stationary in the sense of \cref{def:DP_stationary}
			when applied to \eqref{eq:disjunctive_problem} in the setting \eqref{eq:dp_for_mpcc}.
		\item\label{lem:stat_DP_vs_MPCC_SAS}
			The point $\bar x$ is SAS-stationary in the sense of \cref{def:aw_ac_am_as_stationary}
			if and only if $\bar x$ is SAS-stationary in the sense of \cref{def:DP_stationary}
			when applied to \eqref{eq:disjunctive_problem} in the setting \eqref{eq:dp_for_mpcc}.
	\end{enumerate}
\end{lemma}
\begin{proof}
Throughout the proof, let $F\colon\R^n\to\R^\ell$ and $\Gamma\subset\R^\ell$ be
defined as in \eqref{eq:dp_for_mpcc}.
\begin{enumerate}
\item
We just verify the equivalence of the notion of an approximately M-stationary sequence
w.r.t.\ $\bar x$. The second statement is an immediate consequence of this observation.\\
Let $\{(x^k,\lambda^k,\delta^k,\varepsilon^k)\}_{k=1}^\infty$, 
where \eqref{eq:decoupling_lambda_and_delta}, be an approximately
M-stationary sequence w.r.t.\ $\bar x$ in the sense of \cref{def:aw_ac_am_as_stationary}.
Due to the 
clear equivalence of \eqref{eq:as_dp_x} and \eqref{eq:aw1},
it remains to show that \eqref{eq:aw2} - \eqref{eq:aw6} and \eqref{eq:am}
are equivalent to $\lambda^k\in N_\Gamma(F(x^k)-\delta^k)$.
Obviously, we can rewrite \eqref{eq:aw2} by means of
$\lambda^k_g\in N_{\R^m_-}(g(x^k)-\delta_g^k)$,
and \eqref{eq:aw3} is the same as
$\lambda^k_h\in N_{\{0\}^p}(h(x^k)-\delta_h^k)$.
Hence, having the product rule for the limiting normal cone in mind,
in order to show the claim, it remains to prove that 
\eqref{eq:aw4} - \eqref{eq:aw6} and \eqref{eq:am} are equivalent to
\begin{equation}\label{eq:componentwise_normal_cone_to_C}
	\bigl((\lambda^k_G)_i,(\lambda^k_H)_i\bigr)
	\in 
	N_C\bigl(\bigl(-G_i(x^k)-(\delta^k_G)_i,-H_i(x^k)-(\delta^k_H)_i\bigr)\bigr),
	\quad
	\forall i\in [q].
\end{equation}
Clearly, \eqref{eq:componentwise_normal_cone_to_C} is the same as
\[
	\bigl((\lambda^k_G)_i,(\lambda^k_H)_i\bigr)
	\in 
	\begin{cases}
		\R\times\{0\}	&	i\in I_{0+}(x^k,\delta^k),\\
		\{0\}\times\R	&	i\in I_{+0}(x^k,\delta^k),\\
		\{\xi\in\R^2\,|\,\min\{\xi_1,\xi_2\}>0\text{ or }\xi_1\xi_2=0\}	
						&	i\in I_{00}(x^k,\delta^k)
	\end{cases}
\]
together with $\bigl(-G_i(x^k)-(\delta^k_G)_i,-H_i(x^k)-(\delta^k_H)_i\bigr)\in C$
for all $i\in[q]$.
Hence, as we have $I_{00}(x^k,\delta^k)\subset I_{00}(\bar x)$
for all large enough $k\in\N$ by continuity of $G$ and $H$,
conditions \eqref{eq:aw4} - \eqref{eq:aw6} and \eqref{eq:am} imply
\eqref{eq:componentwise_normal_cone_to_C}.
Vice versa, \eqref{eq:componentwise_normal_cone_to_C} clearly
implies \eqref{eq:aw4} - \eqref{eq:aw6}.
For $i\in I_{00}(x^k,\delta^k)$, we find
$\min\{(\lambda^k_G)_i,(\lambda^k_H)_i\}>0$ or $(\lambda^k_G)_i(\lambda^k_H)_i=0$.
For $i\in I_{00}(\bar x)\setminus I_{00}(x^k,\delta^k)$,
\eqref{eq:aw4} and \eqref{eq:aw5} yield $(\lambda^k_G)_i(\lambda^k_H)_i=0$.
Hence, \eqref{eq:am} is valid as well.
\item 
$[\Longrightarrow]$: Whenever $\bar x$ is assumed to be SAS-stationary in the sense
of \cref{def:aw_ac_am_as_stationary}, we find a sequence 
$\{(x^k,\lambda^k,\delta^k,\varepsilon^k)\}_{k=1}^\infty$, 
where \eqref{eq:decoupling_lambda_and_delta}, being strictly approximately
S-stationary w.r.t.\ $\bar x$ in the sense of \cref{def:aw_ac_am_as_stationary}.
Following the proof of the first assertion, 
it suffices to prove $\lambda^k\in\widehat N_\Gamma(F(\bar x))$,
which reduces to $\lambda^k_g\in\widehat N_{\R^m_-}(g(\bar x))$, $\lambda^k_h\in\widehat N_{\{0\}^p}(h(\bar x))$,
and $\bigl((\lambda^k_G)_i,(\lambda^k_H)_i\bigr)\in\widehat N_{C}\bigl((-G_i(\bar x),-H_i(\bar x)\bigr)$
for all $i\in[q]$.
We note that \eqref{eq:aw2} implies $\lambda^k_g\geq 0$.
Furthermore, for each $i\in[m]\setminus I_g(\bar x)$,
we find $g_i(\bar x)<0$, which yields $g_i(x^k)-(\delta^k_g)_i<0$ for all large enough $k\in\N$
due to \eqref{eq:as_dp_conv_standard}. Hence, in that case,
$(\lambda^k_g)_i=0$ follows from \eqref{eq:aw2} for
large enough $k\in\N$. Consequently, $\lambda^k_g\in\widehat N_{\R^m_-}(g(\bar x))$ holds for all
large enough $k\in\N$.
As we have $h(\bar x)=0$ and, thus, $\widehat N_{\{0\}^p}(h(\bar x))=\R^p$,
$\lambda^k_h\in\widehat N_{\{0\}^p}(h(\bar x))$ is trivially satisfied.
Finally, the continuity of $G$ and $H$,
the convergences \eqref{eq:as_dp_conv_standard},
and \eqref{eq:aw4} - \eqref{eq:aw6} imply the inclusions
\[
	I_{0+}(\bar x) \subset I_{0+}(x^k,\delta^k),
	\qquad
	I_{+0}(\bar x) \subset I_{+0}(x^k,\delta^k),
\] 
and therefore, together with \eqref{eq:as}, the relations
\begin{align*}
	\bigl((\lambda^k_G)_i,(\lambda^k_H)_i\bigr)
	&\in 
	\begin{cases}
		\R\times\{0\}	&	i\in I_{0+}(\bar x),\\
		\{0\}\times\R	&	i\in I_{+0}(\bar x),\\
		\{\xi\in\R^2\,|\,\min\{\xi_1,\xi_2\}\geq 0\} &	i\in I_{00}(\bar x)
	\end{cases}
	\\
	&= \widehat N_{C}\bigl(\bigl(-G_i(\bar x),-H_i(\bar x)\bigr)\bigr)
\end{align*}
for all $i\in[q]$.
Thus, $\{(x^k,\lambda^k,\varepsilon^k)\}_{k=1}^\infty$ is a strictly
approximately S-stationary sequence w.r.t.\ $\bar x$ 
in the sense of \cref{def:DP_special_sequence_SAS}
when applied to \eqref{eq:disjunctive_problem} in the setting \eqref{eq:dp_for_mpcc}.
Particularly, $\bar x$ is SAS-stationary in the sense of \cref{def:DP_stationary} when 
applied to \eqref{eq:disjunctive_problem} in the setting \eqref{eq:dp_for_mpcc}.

$[\Longleftarrow]$: Let $\bar x$ be SAS-stationary in the sense of \cref{def:DP_stationary}
when applied to \eqref{eq:disjunctive_problem} in the setting \eqref{eq:dp_for_mpcc}.
Then we find a sequence 
$\{(x^k,\lambda^k,\varepsilon^k)\}_{k=1}^\infty$, 
where $\lambda^k=(\lambda^k_g,\lambda^k_h,\lambda^k_G,\lambda^k_H)$,
satisfying the convergences \eqref{eq:as_dp_conv_reduced}
as well as relations \eqref{eq:aw1}
and $\lambda^k\in \widehat{N}_\Gamma(F(\bar{x}))$ for each $k\in\N$, 
where the latter amounts to 
$\lambda^k_h\in\R^p$ and
\begin{subequations}\label{eq:actual_SAS_stat}
	\begin{align}
	\label{eq:actual_SAS_stat1}
	& (\lambda^k_g)_i\geq 0,\quad (\lambda^k_g)_i\,g_i(\bar x)=0,\quad 
		\forall i\in[m],
	\\	
	\label{eq:actual_SAS_stat2}
	& (\lambda^k_G)_i G_i(\bar x) = 0,\quad \forall i\in[q],
	\\
	\label{eq:actual_SAS_stat3}
	& (\lambda^k_H)_i H_i(\bar x) = 0,\quad \forall i\in[q],
	\\
	\label{eq:actual_SAS_stat4}
	& \min\{(\lambda^k_G)_i,(\lambda^k_H)_i\} \geq 0,\quad \forall i\in I_{00}(\bar{x})
	\end{align}
\end{subequations}
for all $k\in\N$, showing immediately that \eqref{eq:as} is satisfied.
Let us now define a sequence $\{\tilde\delta^k\}_{k=1}^\infty\subset\R^\ell$,
where $\tilde\delta^k=(\tilde\delta^k_g,\tilde\delta^k_h,\tilde\delta^k_G,\tilde\delta^k_H)$,
via
\[
	\tilde\delta_g^k \coloneqq g(x^k)-g(\bar{x}), \quad
	\tilde\delta_h^k\coloneqq h(x^k), \quad
	\tilde\delta_G^k \coloneqq G(\bar{x})-G(x^k), \quad
	\tilde\delta_H^k \coloneqq H(\bar{x})-H(x^k)
\]
for all $k\in\N$.
Using \eqref{eq:actual_SAS_stat} and the feasibility of $\bar{x}$ for \eqref{eq:MPCC}, 
we immediately see that, for all $k\in\N$, 
all remaining conditions from \eqref{eq:aw} are fulfilled for the sequence 
$\{(x^k,\lambda^k,\tilde\delta^k,\varepsilon^k)\}_{k=1}^\infty$.
Thus, $\{(x^k,\lambda^k,\tilde\delta^k,\varepsilon^k)\}_{k=1}^\infty$
is a strictly approximately S-stationary sequence w.r.t.\ $\bar x$ in the sense of
\cref{def:aw_ac_am_as_stationary}, i.e., $\bar x$ is SAS-stationary in the sense of
\cref{def:aw_ac_am_as_stationary}. \qedhere
\end{enumerate} 
\end{proof}

As we will prove in the upcoming two lemmas, the concepts 
of AW-, AC-, and AM-stationarity from \cref{def:aw_ac_am_as_stationary}
are equivalent to the corresponding ones defined in
\cite[Definitions~3.2,~3.3]{AndreaniHaeserSeccinSilva2019}.
To the best of our knowledge, the concept of SAS-stationarity 
from \cref{def:aw_ac_am_as_stationary} is new.

To start, we show that AW-stationarity from \cref{def:aw_ac_am_as_stationary}
is the same as AW-stationarity from \cite[Definition~3.2]{AndreaniHaeserSeccinSilva2019}.

\begin{lemma}\label{lem:def_equi}
	A point $\bar{x}\in\mathbb{R}^n$ feasible for \eqref{eq:MPCC} is AW-stationary 
	if and only if there exists a sequence 
	$\{(x^k,\lambda^k) \}_{k=1}^\infty\subset\mathbb{R}^{n+\ell}$,
	where $\lambda^k=(\lambda^k_g,\lambda^k_h,\lambda^k_G,\lambda^k_H)$, 
	such that the following conditions hold:
	\begin{subequations}\label{eq:def_andreani}
		\begin{align}
			\label{eq:def_andreani4} 
			& x^k \to \bar{x},
			\\
			\label{eq:def_andreani1}   
			&\nabla_x\mathcal L^{\textup{cc}}(x^k,\lambda^k)\to 0,  
			\\
			\label{eq:def_andreani2}
			&\{\lambda^k_g\}_{k=1}^\infty\subset\R^m_+, \quad \min\{-g_i(x^k),(\lambda^k_g)_i\}\to 0, \
			\quad \forall i \in [m],
			\\
			\label{eq:def_andreani3}
			&\min \bigl\{ G_i(x^k), |(\lambda^k_G)_i| \bigr\}\to 0,
			\quad
				\min \bigl\{ H_i(x^k), |(\lambda^k_H)_i| \bigr\} \to 0,\quad 
				\forall i \in [q].
		\end{align}
	\end{subequations} 
\end{lemma}
\begin{proof}
	Clearly, \eqref{eq:def_andreani1} is satisfied if and only if 
	there exists a sequence
	$\{ \varepsilon^k \}_{k=1}^\infty$ with $\varepsilon^k \to 0$ as $k\to \infty$ such that
	\eqref{eq:aw1} holds for all $k\in\N$.
	Thus, we now show that $\bar{x}$ is AW-stationary if and
	only if, for each $k\in\N$, \eqref{eq:aw1}, \eqref{eq:def_andreani4}, \eqref{eq:def_andreani2}, and \eqref{eq:def_andreani3}
	hold for a sequence 
	$\{(x^k,\lambda^k,\varepsilon^k) \}_{k=1}^\infty\subset\R^{n+\ell+n}$,
	where $\lambda^k=(\lambda^k_g,\lambda^k_h,\lambda^k_G,\lambda^k_H)$,
	such that $\varepsilon^k\to 0$.
	
	$ [ \Longrightarrow ]$:
	First, let $\{(x^k,\lambda^k,\delta^k, \varepsilon^k) \}_{k=1}^\infty$,
	where \eqref{eq:decoupling_lambda_and_delta},
	be an approximately W-stationary sequence w.r.t.\ $\bar x$.
	We immediately find \eqref{eq:def_andreani4} and $\varepsilon^k\to 0$ from \eqref{eq:as_dp_conv_standard}. 
	It remains to show validity of the convergences \eqref{eq:def_andreani2} and \eqref{eq:def_andreani3}.
		
	Let us start to prove \eqref{eq:def_andreani2}.
	As \eqref{eq:aw2} holds for each $k\in\N$, we find $\{\lambda^k_g\}_{k=1}^\infty\subset\R^m_+$.
	Pick $i\in[m]$ arbitrarily.
	In the case $g_i(\bar{x}) < 0$, by continuity of $g_i$ and \eqref{eq:as_dp_conv_standard}, we find
	$g_i(x^k)<0$ for sufficiently large $k\in\mathbb{N}$. In turn, as $\delta_g^k \to 0$ by 
	\eqref{eq:as_dp_conv_standard}, it follows that $g_i(x^k) \neq (\delta_g^k)_i$ for large enough
	$k\in\mathbb{N}$. Then \eqref{eq:aw2} implies $(\lambda^k_g)_i = 0$ and, thus, 
	$\min \{ -g_i(x^k), (\lambda^k_g)_i \} = 0$ for sufficiently large $k\in\mathbb{N}$. In the case 
	$g_i(\bar{x}) = 0$, we know that $g_i(x^k) \to 0$. As $\lambda^k_g \geq 0$ holds for all $k\in\N$, 
	this yields
	$\min\{-g_i(x^k),(\lambda^k_g)_i\} \to 0$, such that \eqref{eq:def_andreani2} is satisfied.
	
	Next, we will show validity of \eqref{eq:def_andreani3}.
	First, pick $i\in I_{+0}(\bar x)$. Then, 
	repeating the arguments from above, we find $G_i(x^k)>0$ and, thus, as $\delta_G^k \to 0$, 
	 $-G_i(x^k) \neq (\delta_G^k)_i$ for sufficiently large $k\in\mathbb{N}$. By 
	\eqref{eq:aw4}, this implies $(\lambda^k_G)_i = 0$ and, hence, $\min\{ G_i(x^k), |(\lambda^k_G)_i| \} = 0$ 
	for $k\in\mathbb{N}$ large enough.
	For the case $i\in I^G(\bar x)$, we know that $G_i(x^k) \to 0$, which clearly implies 
	$\min\{ G_i(x^k), |(\lambda^k_G)_i| \} \to 0$. 
	Analogously, $\min\{ H_i(x^k), |(\lambda^k_H)_i| \} \to 0$ can be proven
	for each $i\in[q]$. 
	Thus, \eqref{eq:def_andreani3} is shown as well.
	
	$ [ \Longleftarrow ] $:
	Assume that there is a sequence 
	$\{(x^k,\lambda^k,\varepsilon^k) \}_{k=1}^\infty$,
	where $\lambda^k=(\lambda^k_g,\lambda^k_h,\lambda^k_G,\lambda^k_H)$,
	that satisfies, for each $k\in\N$, \eqref{eq:aw1}, \eqref{eq:def_andreani4}, 
	\eqref{eq:def_andreani2}, \eqref{eq:def_andreani3}, and $\varepsilon^k\to 0$. 
	Subsequently, we construct a residual sequence $\{\delta^k\}_{k=1}^\infty$,
	where $\delta^k=(\delta^k_g,\delta^k_h,\delta^k_G,\delta^k_H)$,
	with $\delta^k\to 0$ to verify AW-stationarity of $\bar x$.
	To proceed, let us address the different types of constraints separately.
	
	Pick an arbitrary index $i\in I^g(\bar x)$.	
	By setting $(\delta_g^k)_i \coloneqq g_i(x^k)$ for each $k\in\N$ and using 
	\eqref{eq:def_andreani2}, all parts of \eqref{eq:aw2} are clearly fulfilled.
	Continuity of $g_i$ yields $(\delta_g^k)_i = g_i(x^k)\to g_i(\bar{x})=0$ as desired.
	For later use, let us set $(\bar\lambda^k_g)_i\coloneqq(\lambda^k_g)_i$ for each $k\in\N$ in this case.
 	In the case $i\in[m]\setminus I^g(\bar x)$, 
 	\eqref{eq:def_andreani2} yields $(\lambda^k_g)_i\to 0$,
 	and we set $(\bar\lambda^k_g)_i\coloneqq 0$ and 
 	$(\delta_g^k)_i \coloneqq \max\{g_i(x^k),1/k\}$
 	for each $k\in\N$.
 	By construction, we have $g_i(x^k)\leq (\delta^k_g)_i$ for each $k\in\N$,
 	and $g_i(x^k) \to g_i(\bar{x})<0$ yields $(\delta_g^k)_i\to 0$.
 	Furthermore, the complementarity slackness condition $(\bar\lambda^k_g)_i(g_i(x^k)-(\delta^k_g)_i)=0$ 
 	holds trivially for each $k\in\N$.
 	Thus, \eqref{eq:aw2} with $(\lambda^k_g)_i\coloneqq(\bar\lambda^k_g)_i$ 
 	holds for all $i\in[m]$ and $k\in\N$, 
 	and $\delta_g^k \to 0$ is fulfilled. 
	
 	Next, we are going to show that \eqref{eq:aw3} holds. As $h(x^k)\to h(\bar{x})$ by continuity
 	of $h$ and \eqref{eq:def_andreani4}, we know that $h(x^k) \to 0$ holds as $\bar{x}$ is 
 	feasible. This allows us to set $\delta^k_h \coloneqq h(x^k)$ to fulfill \eqref{eq:aw3} for each $k\in\N$
 	while $\delta_h^k \to 0$ is valid. For later use, we set 
 	$\bar\lambda^k_h\coloneqq \lambda^k_h$ for each $k\in\N$.
 	
 	Let us investigate an index from $i\in I^G(\bar x)$. 
 	We set $(\delta_G^k)_i \coloneqq -G_i(x^k)$, which implies 
 	that \eqref{eq:aw4} and \eqref{eq:aw6} are trivially fulfilled in component $i$. 
 	Moreover, this choice ensures 
 	$(\delta_G^k)_i = -G_i(x^k) \to -G_i(\bar{x})=0$ by continuity of $G_i$. 
	For later use, we set $(\bar\lambda^k_G)_i\coloneqq(\lambda^k_G)_i$ for each $k\in\N$.	
 	If $i\in I_{+0}(\bar x)$, \eqref{eq:def_andreani3} yields $(\lambda^k_G)_i\to 0$. 
	Let us set $(\bar\lambda^k_G)_i\coloneqq 0$ and
	$(\delta^k_G)_i\coloneqq\max\{-G_i(x^k),1/k\}$ for each $k\in\N$.
	By construction, we have $G_i(x^k)\geq -(\delta^k_G)_i$ for each $k\in\N$,
	and $G_i(x^k) \to G_i(\bar{x}) > 0$ yields $(\delta_G^k)_i \to 0$. 
	Furthermore, the complementarity slackness condition 
	$(\bar\lambda^k_G)_i(G_i(x^k)+(\delta^k_G)_i)=0$ 
	holds trivially for each $k\in\N$. 
	Thus, \eqref{eq:aw4} with $(\lambda^k_G)_i\coloneqq(\bar\lambda^k_G)_i$ holds for all $i\in[q]$ and $k\in\N$,
	and $\delta^k_G\to 0$ is fulfilled.
	
	Analogously, we may set, for all $i\in[q]$ and $k\in\N$,
	\[
		(\delta^k_H)_i
		\coloneqq
		\begin{cases}
			-H_i(x^k)	&	i\in I^H(\bar x),
			\\
			\max\{-H_i(x^k),1/k\}	&	i\in I_{0+}(\bar x),
		\end{cases}
		\qquad
		(\bar\lambda^k_H)_i
		\coloneqq
		\begin{cases}
			(\lambda^k_H)_i	&	i\in I^H(\bar x),
			\\
			0				&	i\in I_{0+}(\bar x)
		\end{cases}
	\]
	in order to guarantee that \eqref{eq:aw5} holds with $(\lambda^k_H)_i\coloneqq(\bar{\lambda}^k_H)_i$ 
	for all $i\in[q]$ and $k\in\N$ while $\delta^k_H\to 0$ is fulfilled.
	We also observe that this choice guarantees validity of \eqref{eq:aw6}
	for indices $i\in I^H(\bar x)$.
	Together with the above, \eqref{eq:aw6} holds.
	
	Recall that we verified the convergence $\delta^k\to 0$ above. 
	To finalize the proof, let us set
	\[
		\bar\varepsilon^k
		\coloneqq
		\nabla_x\mathcal L^\textup{cc}(x^k,\bar\lambda^k),
	\]
	where $\bar\lambda^k=(\bar\lambda^k_g,\bar\lambda^k_h,\bar\lambda^k_G,\bar\lambda^k_H)$,
	for each $k\in\N$.
	By construction and from \eqref{eq:aw1}, we find
	\begin{align*}
		\bar\varepsilon^k
		&=
		\varepsilon^k
		+
		\nabla_x\mathcal L^\textup{cc}(x^k,\bar\lambda^k)
		-
		\nabla_x\mathcal L^\textup{cc}(x^k,\lambda^k)
		\\
		&=
		\varepsilon^k
		-
		\sum_{i\in[m]\setminus I^g(\bar x)}(\lambda^k_g)_i\nabla g_i(x^k)
		+
		\sum_{i\in I_{+0}(\bar x)}(\lambda^k_G)_i\nabla G_i(x^k)
		+
		\sum_{i\in I_{0+}(\bar x)}(\lambda^k_H)_i\nabla H_i(x^k)
	\end{align*}	 	
	for all $k\in\N$. Due to \eqref{eq:def_andreani4}, $\varepsilon^k\to 0$, 
	continuous differentiability of all data functions, and, as shown above,
	$(\lambda^k_g)_i\to 0$ if $i\in[m]\setminus I^g(\bar x)$,
	$(\lambda^k_G)_i\to 0$ if $i\in I_{+0}(\bar x)$, and
	$(\lambda^k_H)_i\to 0$ if $i\in I_{0+}(\bar x)$,
	$\bar\varepsilon^k\to 0$ follows.
		
	To conclude, $\{(x^k,\bar\lambda^k,\delta^k,\bar\varepsilon^k)\}_{k=1}^\infty$
	is an approximately W-stationary sequence w.r.t.\ $\bar x$,
	i.e., $\bar x$ is AW-stationary.
\end{proof}

The next lemma shows that our concepts of AC- and AM-stationarity from \cref{def:aw_ac_am_as_stationary}
coincide with the ones from \cite[Definition~3.3]{AndreaniHaeserSeccinSilva2019}.

\begin{lemma}\label{lem:def_equi_c_m}
	Let $\bar x\in\R^n$ be feasible for \eqref{eq:MPCC}.
	Then the following assertions hold.
	\begin{enumerate}
		\item
			The point $\bar{x}$ is AC-stationary 
			if and only if there exists a sequence 
			$\{(x^k,\lambda^k) \}_{k=1}^\infty\subset\mathbb{R}^{n+\ell}$,
			where $\lambda^k=(\lambda^k_g,\lambda^k_h,\lambda^k_G,\lambda^k_H)$, 
			such that \eqref{eq:def_andreani} holds and, additionally, we have
			\begin{equation} \label{eq:def_andreani_ac}	
				\begin{aligned}
				\liminf_{k\to\infty} 
					\min \{
						\max&\{ (\lambda^k_G)_i, -(\lambda^k_H)_i \},
						\\
						\max&\{ -(\lambda^k_G)_i, (\lambda^k_H)_i \}
					\} 
				\geq 
				0, 
				\quad 
				\forall i \in [q].
				\end{aligned}
			\end{equation}
		\item
			The point $\bar{x}$ is AM-stationary 
			if and only if there exists a sequence 
			$\{(x^k,\lambda^k) \}_{k=1}^\infty\subset\mathbb{R}^{n+\ell}$,
			where $\lambda^k=(\lambda^k_g,\lambda^k_h,\lambda^k_G,\lambda^k_H)$, 
			such that \eqref{eq:def_andreani} holds and, additionally, we have
			\begin{equation*}% \label{eq:def_andreani_am}
				\begin{aligned}
				\liminf_{k\to\infty} 
					\min \{ 
						\max& \{ (\lambda^k_G)_i, -(\lambda^k_H)_i \}, 
				\\		\max& \{ -(\lambda^k_G)_i, (\lambda^k_H)_i \},
						\max \{ (\lambda^k_G)_i, (\lambda^k_H)_i \}
					\} 
				\geq 
				0, 
				\quad 
				\forall i \in [q].
				\end{aligned}
			\end{equation*}	
	\end{enumerate}
\end{lemma}
\begin{proof}
	\begin{enumerate}
		\item
			$ [ \Longrightarrow ]$:
			First, suppose that $\bar{x}$ is AC-stationary, i.e., there exists an
			approximately C-stationary sequence 
			$\{(x^k,\lambda^k,\delta^k, \varepsilon^k) \}_{k=1}^\infty$
			w.r.t.\ $\bar x$, where \eqref{eq:decoupling_lambda_and_delta}.
			Then $\{(x^k,\lambda^k,\delta^k,\varepsilon^k)\}_{k=1}^\infty$
			is an approximately W-stationary sequence w.r.t.\ $\bar x$ as well.
			We have already shown in the proof of \cref{lem:def_equi} that this implies 
			\eqref{eq:def_andreani} for the sequence 
			$\{(x^k,\lambda^k) \}_{k=1}^\infty$.
			Moreover, as \eqref{eq:ac} holds for all $k\in\N$, we find
			\[
				\min \{ 
						\max \{ (\lambda^k_G)_i, -(\lambda^k_H)_i \}, 
						\max \{ -(\lambda^k_G)_i, (\lambda^k_H)_i \} 
					\}
				\geq 0 
			\]
			for all $k\in\N$ and $i\in I_{00}(\bar x)$.
			For $i\in I_{+0}(\bar x)\cup I_{0+}(\bar x)$,
			\eqref{eq:def_andreani3} yields $(\lambda^k_G)_i\to 0$ or $(\lambda^k_H)_i\to 0$,
			and
			\[
				\lim\limits_{k\to\infty}
					\min \{ 
						\max \{ (\lambda^k_G)_i, -(\lambda^k_H)_i \}, 
						\max \{ -(\lambda^k_G)_i, (\lambda^k_H)_i \} 
					\}
				=
				0
			\]
			follows from \cref{lem:andreani_c_function}.
			Altogether, \eqref{eq:def_andreani_ac} is fulfilled.
	
			$ [ \Longleftarrow ]$:
			Assume there exists a sequence 
			$\{(x^k,\lambda^k) \}_{k=1}^\infty$, where $\lambda^k=(\lambda^k_g,\lambda^k_h,\lambda^k_G,\lambda^k_H)$,
			such that \eqref{eq:def_andreani} and \eqref{eq:def_andreani_ac} hold. We have seen
			in the proof of \cref{lem:def_equi} that this implies the existence of a sequence 
			$\{(x^k,\bar\lambda^k,\delta^k,\bar\varepsilon^k)\}_{k=1}^\infty$
			being approximately W-stationary w.r.t.\ $\bar x$ such that
			\[
				(\bar\lambda^k_G)_i
				=
				\begin{cases}
					(\lambda^k_G)_i	&	i\in I^G(\bar x),
					\\
					0		&	i\in I_{+0}(\bar x),
				\end{cases}
				\qquad
				(\bar\lambda^k_H)_i
				=
				\begin{cases}
					(\lambda^k_H)_i	&	i\in I^H(\bar x),
					\\
					0		&	i\in I_{0+}(\bar x)
				\end{cases}
			\]
			holds for all $k\in\N$ and $i\in[q]$.
			Thus, for the AC-stationarity of $\bar{x}$, it remains to make sure
			that these multipliers also satisfy
			$(\bar{\lambda}^k_G)_i(\bar{\lambda}^k_H)_i \geq 0$ for all $i\in I_{00}(\bar{x})$ and 
			$k\in\mathbb{N}$. Due to $i\in I_{00}(\bar{x})$ and the definitions of $\bar{\lambda}^k_G$ 
			and $\bar{\lambda}^k_H$, this is equivalent to showing $(\lambda^k_G)_i(\lambda^k_H)_i \geq 0$ for all 
			$i\in I_{00}(\bar{x})$ and $k\in\mathbb{N}$. 
			For $i\in I_{00}(\bar x)$ such that
			\[
				\liminf\limits_{k\to\infty}
				\min \{ 
						\max \{ (\lambda^k_G)_i, -(\lambda^k_H)_i \}, 
						\max \{ -(\lambda^k_G)_i, (\lambda^k_H)_i \} 
					\} 
				> 0,
			\]
			\cref{lem:andreani_c_function} indeed yields $(\lambda^k_G)_i(\lambda^k_H)_i \geq 0$ for all large enough $k\in\N$.
			Hence, let us assume that $i\in I_{00}(\bar x)$ is chosen such that
			\begin{equation}\label{eq:critical_C_multiplier}
				\liminf\limits_{k\to\infty}
				\min \{ 
						\max \{ (\lambda^k_G)_i, -(\lambda^k_H)_i \}, 
						\max \{ -(\lambda^k_G)_i, (\lambda^k_H)_i \} 
					\}
				= 0.
			\end{equation}
			Then \cref{lem:andreani_c_function} guarantees that $(\lambda^k_G)_i\to 0$ or $(\lambda^k_H)_i\to 0$
			holds along a subsequence (without relabeling).
			Without loss of generality, let us assume that $(\lambda^k_G)_i\to 0$.
			In this situation,
			we redefine $(\bar\lambda^k_G)_i\coloneqq 0$ for all $k\in\N$
			and append the null sequence $\{(\lambda^k_G)_i\nabla G_i(x^k)\}_{k=1}^\infty$ to $\{\bar\varepsilon^k\}_{k=1}^\infty$.
			Trivially, $(\bar{\lambda}^k_G)_i(\bar{\lambda}^k_H)_i \geq 0$ and 
			$\bar\varepsilon^k=\nabla_x\mathcal L^\textup{cc}(x^k,\bar\lambda^k)$
			hold for all $k\in\N$ after this procedure,
			and $\bar\varepsilon^k\to 0$ is still valid.
			Proceeding like this for each index $i\in I_{00}(\bar x)$
			with the property \eqref{eq:critical_C_multiplier}
			shows AC-stationarity of $\bar x$.
		\item 
			With the aid of \cref{lem:def_equi,lem:andreani_m_function},
			this statement can be proven analogously to the preceding one. 
			\qedhere
	\end{enumerate}
\end{proof}

Let us note that we can further show that AM-stationarity is equivalent to the 
fulfillment of the MPEC-AKKT conditions as introduced in \cite[Definition~3.2]{Ramos2021}. 
However, due to the equivalence shown in \cref{lem:def_equi_c_m}, this result can also 
be distilled following the rough arguments stated in 
\cite[p.\ 3211]{AndreaniHaeserSeccinSilva2019}, which is why
we omit a proof here.

It is further possible to show the equivalence of AW-stationarity
to the following formulation which is closely related
to the ones in
\cite[Proposition~4.1]{Mehlitz2023} and \cite[Definition~3.2]{Ramos2021},
see \cite[Theorem~3.2]{BirginMartinez2014}
and \cite[Lemma~3.3]{Mehlitz2023} for related results.
The proof is analogous to the one of \cref{lem:def_equi} and, thus, omitted.
We also would like to note that associated alternative formulations
of AC- and AM-stationarity can be derived.

\begin{lemma}\label{lem:def_equi_def2}
	A point $\bar{x}\in\mathbb{R}^n$ feasible for \eqref{eq:MPCC} is AW-stationary 
	if and only if there exists a sequence 
	$\{(x^k,\lambda^k) \}_{k=1}^\infty\subset \mathbb{R}^{n+\ell}$,
	where $\lambda^k=(\lambda^k_g,\lambda^k_h,\lambda^k_G,\lambda^k_H)$, 
	such that the conditions \eqref{eq:def_andreani4} and \eqref{eq:def_andreani1} hold
	while the relations
	\begin{align*}
		&(\lambda^k_g)_i\geq 0,\quad\min\{-g_i(\bar x),(\lambda^k_g)_i\}=0,
		\quad \forall i\in[m],
		\\
		&\min\{G_i(\bar x),|(\lambda^k_G)_i|\}=0,
		\quad
		\min\{H_i(\bar x),|(\lambda^k_H)_i|\}=0,
		\quad \forall i\in[q]
	\end{align*}
	are valid for each $k\in\N$.
\end{lemma}

Thanks to \cref{cor:local_minimizers_AM_stat} and \cref{lem:stat_DP_vs_MPCC},
we obtain that AM- (and, thus, AC- as well as AW-) stationarity
provides a necessary optimality condition for \eqref{eq:MPCC}.

\begin{theorem}
	If $\bar{x}\in\R^n$ is a local minimizer of \eqref{eq:MPCC}, 
	then $\bar{x}$ is AM-stationary.
\end{theorem}

Due to \cref{lem:def_equi,lem:def_equi_c_m}, 
the above result is also an immediate consequence of \cite[Theorem~3.6]{AndreaniHaeserSeccinSilva2019}.

On the contrary, the following example illustrates that local minimizers of \eqref{eq:MPCC}
do not need to be SAS-stationary in general.
\begin{example}
	We consider the MPCC
	\[
		\min\limits_x\quad -x_1\quad\textup{s.t.}\quad
		x_1-x_2 \leq 0,\quad 0\leq x_1 \perp x_2 \geq 0.
	\]
	The feasible set of this problem is $\{0\}\times\R_+$,
	and hence, the origin $\bar x\coloneqq(0,0)$ is a minimizer.
	To ensure SAS-stationarity of $\bar{x}$, we need to find a sequence 
	$\{(x^k,\lambda^k,\delta^k,\varepsilon^k)\}_{k=1}^\infty\subset\R^{2+3+3+2}$,
	where $\lambda^k = (\lambda_g^k, \lambda_G^k, \lambda_H^k)$ and 
	$\delta^k = (\delta_g^k, \delta_G^k, \delta_H^k)$,
	satisfying the convergences \eqref{eq:as_dp_conv_standard}
	as well as the conditions \eqref{eq:aw} and \eqref{eq:as} for each $k\in\N$.
	Due to \eqref{eq:aw1} and \eqref{eq:as}, we particularly need that the 
	sequence fulfills
	\[
		\varepsilon^k
		=
		\begin{pmatrix}
			-1\\0
		\end{pmatrix}
		+
		\lambda^k_g
		\begin{pmatrix}
			1\\-1
		\end{pmatrix}
		-
		\lambda^k_G
		\begin{pmatrix}
			1\\0
		\end{pmatrix}
		-
		\lambda^k_H
		\begin{pmatrix}
			0\\1
		\end{pmatrix}
	\]
	and $\min\{ \lambda_G^k, \lambda_H^k\} \geq 0$
	for each $k\in\N$, respectively.
	Due to $\varepsilon^k\to 0$,
	this requires $\lambda^k_g-\lambda^k_G\to 1$ and $-\lambda^k_g-\lambda^k_H\to 0$.
	Summing these relations up, we end up with $-\lambda^k_G-\lambda^k_H\to 1$,
	which is impossible due to $\min\{ \lambda_G^k, \lambda_H^k\} \geq 0$ for all $k\in\N$.
	Hence, we conclude that $\bar x$ is not SAS-stationary.
\end{example}

Nevertheless, in the presence of MPCC-LICQ, see \cite[Definition~2.8]{Ye2005},
which corresponds to LICQ from \cref{def:MPDC_LICQ} 
when applied to \eqref{eq:disjunctive_problem} in the setting \eqref{eq:dp_for_mpcc3},
see \cite[Section~5.1]{Mehlitz2020a},
each local minimizer of \eqref{eq:MPCC} is S- and, thus, SAS-stationary, as we
recall from \cref{lem:LICQ_stronger} and \cref{lem:stat_DP_vs_MPCC}, see also \cite[Theorem~3.8]{Mehlitz2020a}.
However, \cref{ex:LICQ_stronger}, which also considers a complementarity-constrained problem,
illustrates that local minimizers of \eqref{eq:MPCC} can be SAS-stationary
even in the absence of MPCC-LICQ, notice, once again, \cref{lem:stat_DP_vs_MPCC}.

\subsection{The subset Mangasarian--Fromovitz condition for complementarity-constrained problems}\label{sec:subMFC_MPCC}

In this subsection, we aim to apply ODP-subMFC 
from \cref{def:subMFC_ODP} to the special orthodisjunctive problem \eqref{eq:MPCC}.
Therefore, we reinspect \cref{def:subMFC_ODP} for \eqref{eq:orthodisjunctive_problem}
in the setting \eqref{eq:dp_for_mpcc3}.

Given $x\in\R^n$ and $\delta\in\R^\ell$, where $\delta=(\delta_g,\delta_h,\delta_G,\delta_H)$,
such that $F(x)-\delta\in\Gamma$, 
see \eqref{eq:dp_for_mpcc} for the definitions of $F$ and $\Gamma$,
we need to characterize the associated index set $I^\exists(x,\delta)$.
To start, we observe that a multi-index $j\in\JJ$ belongs to $J(x,\delta)$
if and only if
\[
	(-G_i(x)-(\delta_G)_i,-H_i(x)-(\delta_H)_i) \in C_{j_i}
\]
holds true for all $i\in[q]$. 
For indices $i\in I_{00}(x,\delta)$, the above is trivially satisfied,
and so we end up with the characterization
\[
	j\in J(x,\delta)
	\quad\Longleftrightarrow\quad
	\forall i\in I_{0+}(x,\delta)\colon\,j_i=1
	\;\text{and}\;
	\forall i\in I_{+0}(x,\delta)\colon\,j_i=2,
	\quad
	\forall j\in\JJ.
\]
Now, we are in position to characterize $I^\exists(x,\delta)$ by addressing the different
components in $F$ from \eqref{eq:dp_for_mpcc1} separately.
For the constraint functions $g_1,\ldots,g_m$,
precisely the indices from $I^g(x,\delta)$ belong to $I^\exists(x,\delta)$.
Addressing the constraint functions $h_1,\ldots,h_p$,
all the indices from $I^h$ belong to $I^\exists(x,\delta)$.
Let us now inspect all the constraint functions $G_1,\ldots,G_q$.
Then all indices from $I^G(x,\delta)$ are included in $I^\exists(x,\delta)$.
Similarly, when considering the component functions $H_1,\ldots,H_q$,
all indices from $I^H(x,\delta)$ are part of $I^\exists(x,\delta)$.

\begin{remark}\label{rem:seq_const}
	Given a feasible point $\bar x\in\R^n$ of \eqref{eq:MPCC}, 
	analogously to \cref{rem:dp_seq_const}, whenever we consider a sequence 
	$\{ (x^k,\lambda^k,\delta^k, \varepsilon^k) \}_{k=1}^\infty$ 
	being approximately W-stationary w.r.t.\ $\bar x$ in the following, 
	we may assume without loss of generality that the sets $I^g(x^k,\delta^k)$
	are the same for all $k\in\N$, and likewise, we assume the same for the
	sets $I^G(x^k,\delta^k)$ and $I^H(x^k,\delta^k)$,
	which is reasonable as it is always possible 
	to consider a suitable infinite subsequence.
	Similarly, we may also assume that each component of $\lambda^k$ possesses
	a constant sign for all $k\in\N$.
\end{remark}

Now, based on our considerations above and \cref{def:subMFC_ODP},
we are in position to introduce a qualification condition for
\eqref{eq:MPCC} which is closely related to ODP-subMFC when
applied to \eqref{eq:orthodisjunctive_problem} in the setting
\eqref{eq:dp_for_mpcc3}.

\begin{definition}\label{def:subMFC}
	Let $\bar{x}\in\R^n$ be an AW-stationary point of \eqref{eq:MPCC}. 
	We say that the \textit{MPCC Subset Mangasarian--Fromovitz Condition (MPCC-subMFC)}
	holds at $\bar{x}$ if there 
	exist index sets $I_1 \subseteq I^g(\bar{x})$, $I_2 \subseteq I^G(\bar{x})$, $I_3 \subseteq I^H(\bar{x})$, 
	and a sequence 
	$\{(x^k,\lambda^k,\delta^k, \varepsilon^k) \}_{k=1}^\infty\subset \mathbb{R}^{n+\ell+\ell+n}$,
	where \eqref{eq:decoupling_lambda_and_delta},
	such that the following conditions are satisfied.
	\begin{enumerate}[label=(\roman*)]
		\item \label{subMFCI} 
			Either $I_1=I^h=I_2=I_3=\emptyset$, or
			it holds for all 
			$(u,v,r,s)\in\R^{m+p+2q}\setminus\{(0,0,0,0)\}$
			with $u\geq 0$, $v\geq 0$, $r\geq 0$, $s \geq 0$ and
			$u_{[m]\setminus I_1}=0$, $r_{[q]\setminus I_2} = 0$, $s_{[q]\setminus I_3} = 0$
			that
			\begin{align*}
				0  
				\neq 
				\sum_{i\in I_1} u_i \nabla g_i(\bar{x})
				&+ 
				\sum_{i\in I^h} \sgn((\lambda^k_h)_i) v_i \nabla h_i(\bar{x}) 
				\\
				&- 
				\sum_{i\in I_2} \sgn((\lambda^k_G)_i) r_i \nabla G_i(\bar{x}) 
				- 
				\sum_{i\in I_3} \sgn((\lambda^k_H)_i) s_i \nabla H_i(\bar{x}).
			\end{align*}
		\item \label{subMFCII} 
			The sequence $\{ (x^k,\lambda^k,\delta^k,\varepsilon^k) \}_{k=1}^\infty$ 
			is approximately W-stationary w.r.t.\ $\bar x$, and
			$I_1 = I^g(x^k,\delta^k)$,
			$I_2 = I^G(x^k,\delta^k)$, and
			$I_3 = I^H(x^k,\delta^k)$ are valid for all $k\in\N$.
	\end{enumerate}
\end{definition}

It is important to observe that we do not use an approximately M-stationary sequence w.r.t.\ $\bar x$
in \cref{def:subMFC} but merely an approximately W-stationary sequence w.r.t.\ $\bar x$.
This will allow us to obtain results similar to \cref{thm:ODP_subMFC}
for all four types of stationarity from \cref{def:stat_MPCC}
which have been introduced for \eqref{eq:MPCC}.
Furthermore, let us emphasize that, 
for the treatment of the pure inequality constraints, we already removed
the sign condition in \cref{def:subMFC}~\ref{subMFCI},
as we already explained in \cref{sec:NLPs}.

Recall that, given a feasible point $\bar x\in\R^n$ for \eqref{eq:MPCC},
MPCC-MFCQ is said to hold at $\bar x$ if,
for all $(u,v,r,s)\in\R^{m+p+2q}\setminus\{(0,0,0,0)\}$ 
with $u\geq 0$, $u_{[m]\setminus I^g(\bar x)}=0$, $r_{[q]\setminus I^G(\bar x)}=0$,
$s_{[q]\setminus I^H(\bar x)}=0$,
\[ 
	0 \neq \sum_{i\in I^g(\bar{x})} u_i \nabla g_i(\bar{x}) 
	+ \sum_{i\in I^h} v_i \nabla h_i(\bar{x}) 
	- \sum_{i\in I^G(\bar{x})} r_i \nabla G_i(\bar{x}) 
	- \sum_{i\in I^H(\bar{x})} s_i \nabla H_i(\bar{x})
\]
is satisfied.
It is important to note that MPCC-subMFC from \cref{def:subMFC} relaxes MPCC-MFCQ in two different ways: 
First, the set of gradients for which the positive linear independence is demanded 
in item \ref{subMFCI} of \cref{def:subMFC} is a \textit{subset} of the gradients considered 
in MPCC-MFCQ as $I_1 \subseteq I^g(\bar{x})$,
$I_2 \subseteq I^G(\bar{x})$, and $I_3 \subseteq I^H(\bar{x})$. 
Moreover, in the positive linear independence condition in MPCC-MFCQ, 
all gradients besides the ones belonging to the 
functions $g_i$, $i\in I^g(\bar x)$, are considered with a multiplier of arbitrary sign. 
In MPCC-subMFC, for \textit{all} gradients, the sign of the multipliers required for the positive 
linear independence is fixed, see \cref{rem:seq_const} again, which makes MPCC-subMFC easier to fulfill.

The upcoming theorem illustrates the value of MPCC-subMFC and is closely related to 
\cref{thm:ODP_subMFC}.

\begin{theorem}\label{thm:submfc_implications}
	Let $\bar{x}\in\R^n$ be an AW-stationary point of \eqref{eq:MPCC},
	and let MPCC-subMFC be satisfied at $\bar{x}$ 
	with $\{ (x^k,\lambda^k,\delta^k,\varepsilon^k) \}_{k=1}^\infty\subset\R^{n+\ell+\ell+n}$,
	where \eqref{eq:decoupling_lambda_and_delta},
	being the involved approximately W-stationary sequence w.r.t.\ $\bar x$.
	Then the following assertions hold.
	\begin{enumerate}
		\item \label{thm:W} The point $\bar{x}$ is W-stationary.  
		\item \label{thm:C} Whenever $\{ (x^k,\lambda^k,\delta^k,\varepsilon^k) \}_{k=1}^\infty$
			is an approximately C-stationary sequence w.r.t.\ $\bar x$, then $\bar{x}$ is C-stationary.   
		\item \label{thm:M} Whenever $\{ (x^k,\lambda^k,\delta^k,\varepsilon^k) \}_{k=1}^\infty$
			is an approximately M-stationary sequence w.r.t.\ $\bar x$, then $\bar{x}$ is M-stationary.  
		\item \label{thm:S} Whenever $\{ (x^k,\lambda^k,\delta^k,\varepsilon^k) \}_{k=1}^\infty$
			is a strictly approximately S-stationary sequence w.r.t.\ $\bar x$, then $\bar{x}$ is S-stationary.   
	\end{enumerate}
\end{theorem}
\begin{proof}
	To start the proof, we elaborate on some general observations.
	Let $I_1\subset I^g(\bar x)$, $I_2\subset I^G(\bar x)$, and $I_3\subset I^H(\bar x)$
	be the index sets whose existence is guaranteed by \cref{def:subMFC}.
	By item \ref{subMFCII} of \cref{def:subMFC}, we find
	\begin{equation}\label{eq:proof1}
		\begin{aligned}   
			\varepsilon^k 
			= 
			\nabla f(x^k) + \sum_{i\in I_1} (\lambda^k_g)_i \nabla g_i(x^k) 
					&+ \sum_{i\in I^h} (\lambda^k_h)_i \nabla h_i(x^k) 
					\\
					&
					- \sum_{i\in I_2} (\lambda^k_G)_i \nabla G_i(x^k) 
					- \sum_{i\in I_3} (\lambda^k_H)_i \nabla H_i(x^k), 
		\end{aligned}
	\end{equation}
	as well as \eqref{eq:aw2} - \eqref{eq:aw6} for all $k\in\mathbb{N}$
	and the convergences \eqref{eq:as_dp_conv_standard}.
	If $I_1=I^h=I_2=I_3=\emptyset$, the continuity of $\nabla f$ implies 
	$0 = \nabla f(\bar{x})$	when taking the limit $x^k\to\bar{x}$ in \eqref{eq:proof1}.
	As $\bar{x}$ is feasible to \eqref{eq:MPCC}, it is already
	S-stationary with multipliers 
	$(\lambda_g, \lambda_h, \lambda_G,\lambda_H)\coloneqq (0,0,0,0)$.
	Hence, in the remainder	of the proof, we do not need to discuss this
	trivial situation anymore.		

	For the case where at least one $I_j$, $j\in\{1,2,3\}$, or $I^h$ is nonempty,
	we can show, similarly to the proof of \cref{thm:ODP_subMFC}, 
	that the sequence 
	$\{ ((\lambda^k_g)_{I_1}, \lambda^k_h, (\lambda^k_G)_{I_2}, (\lambda^k_H)_{I_3}) \}_{k=1}^\infty$
	is bounded by exploiting \cref{def:subMFC}~\ref{subMFCI}.
	Without loss of generality, 
	we may pick a quadruple	
	$
			(\widetilde{\lambda}_g, \widetilde{\lambda}_h, 
				\widetilde{\lambda}_G, \widetilde{\lambda}_H)
			\in 
			\R^{m+p+2q}
	$
	such that 
	\begin{equation}\label{eq:trivial_multipliers}
		(\widetilde\lambda_g)_{[m]\setminus I_1} = 0, \quad
		(\widetilde\lambda_G)_{[q]\setminus I_2} = 0, \quad
		(\widetilde\lambda_H)_{[q]\setminus I_3} = 0,
	\end{equation}
	and 
	$\{ ((\lambda^k_g)_{I_1}, \lambda^k_h, (\lambda^k_G)_{I_2}, (\lambda^k_H)_{I_3}) \}_{k=1}^\infty$
	converges to
	$
			((\widetilde{\lambda}_g)_{I_1}, \widetilde{\lambda}_h, 
				(\widetilde{\lambda}_G)_{I_2}, (\widetilde{\lambda}_H)_{I_3})
	$. 
	Taking the limit $k\to\infty$ in \eqref{eq:proof1}, it follows from 
	\eqref{eq:as_dp_conv_standard} and continuous differentiability of all involved functions that
	\[
		0
		=
		\nabla f(\bar x)
		+
		\sum_{i\in I_1}(\widetilde\lambda_g)_i\nabla g_i(\bar x)
		+
		\sum_{i\in I^h}(\widetilde\lambda_h)_i\nabla h_i(\bar x)
		-
		\sum_{i\in I_2}(\widetilde\lambda_G)_i\nabla G_i(\bar x)
		-
		\sum_{i\in I_3}(\widetilde\lambda_H)_i\nabla H_i(\bar x).
	\] 
	Furthermore, \eqref{eq:as_dp_conv_standard}, \eqref{eq:aw2}, and \eqref{eq:trivial_multipliers}  yield
	the complementarity slackness conditions 
	$(\widetilde\lambda_g)_i\geq 0$ and $(\widetilde\lambda_g)_i\,g_i(\bar x)=0$
	for all $i\in[m]$.
	Note that $\supp(\widetilde\lambda_G)\subset I_2\subset I^G(\bar x)$
	yields $(\widetilde \lambda_G)_i=0$ for all $i\in I_{+0}(\bar x)$,
	and $(\widetilde\lambda_H)_i=0$ for all $i\in I_{0+}(\bar x)$ is obtained similarly.
	Hence, we have shown that $\bar x$ is a W-stationary point of \eqref{eq:MPCC},
	which proves assertion~\ref{thm:W}.
	
	To show assertion~\ref{thm:C},
	assume that, additionally,
	\eqref{eq:ac} holds for all $k\in\N$.
	Then we also have 
	$(\widetilde{\lambda}_G)_i (\widetilde{\lambda}_H)_i = \lim_{k\to\infty} (\lambda^k_G)_i (\lambda^k_H)_i \geq 0$
	for all $i\in I_{00}(\bar{x})$, 
	such that $\bar{x}$ is a C-stationary point with multipliers 
	$(\widetilde{\lambda}_g, \widetilde{\lambda}_h, \widetilde{\lambda}_G, \widetilde{\lambda}_H)$.
		
	Let us proceed by verifying assertion~\ref{thm:M}.	
	To this end, we assume, additionally, that \eqref{eq:am} holds for all $k\in\N$.
	Fix $i\in I_{00}(\bar x)$.
	Whenever $\min\{(\lambda^k_G)_i,(\lambda^k_H)_i\}>0$ holds for
	infinitely many $k\in\N$, we find $\min\{(\widetilde\lambda_G)_i,(\widetilde\lambda_H)_i\}\geq 0$
	when taking the limit.
	Observe that this yields $\min\{(\widetilde\lambda_G)_i,(\widetilde\lambda_H)_i\}>0$
	or $(\widetilde\lambda_G)_i\,(\widetilde\lambda_H)_i = 0$.
	If $(\lambda^k_G)_i\,(\lambda^k_H)_i=0$ is valid for infinitely many $k\in\N$,
	$(\widetilde\lambda_G)_i\,(\widetilde\lambda_H)_i=0$ follows when taking the limit.
	Altogether, we find $\min\{(\widetilde\lambda_G)_i(\widetilde\lambda_H)_i\} > 0$
	or $(\widetilde\lambda_G)_i\,(\widetilde\lambda_H)_i = 0$,
	i.e., $\bar x$ is M-stationary with multipliers 
	$(\widetilde{\lambda}_g, \widetilde{\lambda}_h, \widetilde{\lambda}_G, \widetilde{\lambda}_H)$.
			
	Assertion~\ref{thm:S} follows in a similar way as the preceding two statements.
\end{proof}

Let us note that \cref{thm:submfc_implications}\,\ref{thm:M} can also be distilled
from \cref{thm:ODP_subMFC}\,\ref{item:ODP_subMFC_AM} directly, thanks to 
\cref{lem:stat_DP_vs_MPCC}\,\ref{lem:stat_DP_vs_MPCC_AM}, 
whereas \cref{thm:submfc_implications}\,\ref{thm:S} can be obtained by combining 
\cref{thm:ODP_subMFC}\,\ref{item:ODP_subMFC_AS}
and (the proof of) \cref{lem:stat_DP_vs_MPCC}\,\ref{lem:stat_DP_vs_MPCC_SAS}.

\begin{remark}\label{rem:multiplier_conditions}
	The additional conditions \eqref{eq:ac} and \eqref{eq:am} 
	on the multipliers $(\lambda^k_G)_i$ and $(\lambda^k_H)_i$ 
	in assertions~\ref{thm:C} and~\ref{thm:M} of \cref{thm:submfc_implications},
	hidden in the definition of an approximately C- and M-stationary sequence
	w.r.t.\ $\bar x$, respectively,
	only have to be imposed for all $i\in I_2 \cap I_3$ (instead of $i\in I_{00}(\bar{x})$). 
	Indeed, let us reinspect the proof of \cref{thm:submfc_implications}.
	For all $i\in I_{00}(\bar{x}) \setminus
	(I_2 \cap I_3)$, we have $i\in[q]\setminus I_2$ or $i\in[q]\setminus I_3$, 
	which, by \eqref{eq:trivial_multipliers}, 
	implies $(\widetilde{\lambda}_G)_i=0$ or $(\widetilde{\lambda}_H)_i=0$, respectively. 
	Hence, $(\widetilde{\lambda}_G)_i\, (\widetilde\lambda_H)_i = 0$ follows 
	for all $i\in I_{00}(\bar x)$.
	Note that this is not enough to distill S-stationarity of $\bar x$,
	i.e., a similar simplification of \cref{thm:submfc_implications}\,\ref{thm:S}
	is not possible.
\end{remark}

\section{Final remarks}\label{sec:final_rem}
In this paper, we investigated concepts of approximate stationarity and their 
connection to exact stationarity via regularity conditions and qualification 
conditions for problems with disjunctive constraints. Besides applying a 
tailored concept of AM-stationarity to infer M-stationarity in a limit point, 
we proposed the novel concept of SAS-stationarity, which is well-suited to 
guarantee even S-stationarity in a limit point. 
However, in contrast to AM-stationarity, SAS-stationarity is not a necessary 
optimality condition that holds without further assumptions, and it remains 
an open task to fully grasp the entirety of tractable conditions that 
guarantee SAS-stationarity of a local minimizer of 
\eqref{eq:disjunctive_problem}.
Building upon the considered approximate stationarity concepts, we introduced 
a new qualification condition that enables an alternative to using approximate 
regularity conditions in order to recognize an accumulation point's exact
stationary from approximate stationarity. Finally, we illustrated our findings
in the context of purely inequality-constrained problems and MPCCs, and
extended the available theory of approximate stationarity and 
qualification conditions for these problem classes.
Clearly, as \eqref{eq:orthodisjunctive_problem} covers many other prominent 
problem classes, such as problems with vanishing constraints, switching 
constraints, cardinality constraints, 
reformulated cardinality constraints of complementarity type, or 
relaxed probabilistic constraints, our results are likely to possess 
reasonable extensions to these problem classes as well.
From a practical viewpoint it can further be promising to explore the 
capability of the qualification conditions when numerically treating the 
aforementioned problem classes, 
see \cref{rem:algorithmic_consequences} as well.
Finally, following \cite{KaemingFischerZemkoho2024,Mehlitz2020b},
it is likely that many of our findings can be generalized to
(ortho)disjunctive optimization problems that are modeled with the
aid of nonsmooth functions. 

\subsection*{Acknowledgments}

We thank the reviewers for their constructive assessment that helped
to improve the presentation of the paper.
In particular, we acknowledge the comments of one reviewer which led
to the observations stated in \cref{rem:algorithmic_consequences}.

\appendix
\section{Appendix}\label{sec:appendix}

We consider the piecewise affine function $\varphi\colon\R^2\to\R$ given by
\[
	\varphi(a,b)
	\coloneqq
	\min\{\max\{a,-b\},\max\{-a,b\}\},
	\quad
	\forall (a,b)\in\R^2.
\]
A distinction of cases shows that
\begin{equation}\label{eq:explicit_rep_varphi}
	\varphi(a,b)
	=
	\begin{cases}
		a		&	|a|\leq b,
		\\
		-b		&	|b|\leq -a,
		\\
		-a		&	|a|\leq -b,
		\\
		b		&	|b|\leq a,
	\end{cases}
	\quad
	\forall (a,b)\in\R^2.
\end{equation}
In the subsequently stated lemma, we summarize some important properties of $\varphi$.

\begin{lemma}\label{lem:andreani_c_function}
	Let $\{a_k\}_{k=1}^\infty,\{b_k\}_{k=1}^\infty\subset\R$ be arbitrary sequences.
	\begin{enumerate}
		\item If $a_k\to 0$ or $b_k\to 0$, then $\varphi(a_k,b_k)\to 0$.
		\item If $\liminf_{k\to\infty}\varphi(a_k,b_k)>0$, then $a_kb_k\geq 0$ holds
			for all sufficiently large $k\in\N$.
		\item If $\liminf_{k\to\infty}\varphi(a_k,b_k)=0$,
			then $\{a_k\}_{k=1}^\infty$ or $\{b_k\}_{k=1}^\infty$ comprises a subsequence
			converging to $0$.
	\end{enumerate}
\end{lemma}
\begin{proof}
	The first and third assertion are obvious from \eqref{eq:explicit_rep_varphi}.
	In order to prove the second assertion,
	we assume on the contrary that there is an infinite set $K\subset\N$
	such that $a_kb_k<0$ for all $k\in K$.
	Then, $a_k<0$ and $b_k>0$, or $a_k>0$ and $b_k<0$ holds for each $k\in K$.
	Particularly, $\varphi(a_k,b_k)<0$ follows for each $k\in K$ 
	from \eqref{eq:explicit_rep_varphi}, and consequently,
	$\liminf_{k\to\infty}\varphi(a_k,b_k)\leq 0$,
	contradicting the assumptions.
\end{proof}

Let us also comment on the piecewise affine function $\psi\colon\R^2\to\R$ given by
\[
	\psi(a,b)
	\coloneqq
	\min\{\max\{a,-b\},\max\{-a,b\},\max\{a,b\}\},
	\quad
	\forall (a,b)\in\R^2.
\]
A distinction of cases shows that
\[
	\psi(a,b)
	=
	\begin{cases}
		a		&	|a|\leq b,
		\\
		-|b|	&	|b|\leq -a,
		\\
		-|a|	&	|a|\leq -b,
		\\
		b		&	|b|\leq a,
	\end{cases}
	\quad
	\forall (a,b)\in\R^2.
\]
The following result is very similar to \cref{lem:andreani_c_function}, which is why
a proof is omitted.
\begin{lemma}\label{lem:andreani_m_function}
	Let $\{a_k\}_{k=1}^\infty,\{b_k\}_{k=1}^\infty\subset\R$ be arbitrary sequences.
	\begin{enumerate}
		\item If $a_k\to 0$ or $b_k\to 0$, then $\psi(a_k,b_k)\to 0$.
		\item If $\liminf_{k\to\infty}\psi(a_k,b_k)>0$, then,
			for all sufficiently large $k\in\N$, 
			$\min\{a_k, b_k\} > 0$ or $a_kb_k=0$.
		\item If $\liminf_{k\to\infty}\psi(a_k,b_k)=0$,
			then $\{a_k\}_{k=1}^\infty$ or $\{b_k\}_{k=1}^\infty$ comprises a subsequence
			converging to $0$.
	\end{enumerate}
\end{lemma}

\end{document}